\documentclass[12pt]{article}

\usepackage{epic}
\usepackage{eepic}
\usepackage{epsf} 
\usepackage{ amssymb, amscd}
\usepackage{tikz} 
\usepackage{amsmath}
\def\cA            {{\mathcal{A}}}
\def\cB            {{\mathcal{B}}}
\def\cC            {{\mathcal{C}}}
\def\cD            {{\mathcal{D}}}

\def\cF            {{\mathcal{F}}}
\def\cG            {{\mathcal{G}}}
\def\cS            {{\mathcal{S}}}

\def\cJ            {{\mathcal{J}}}

\def\cV            {{\mathcal{V}}}

\def\cW            {{\mathcal{W}}}
\def\cT            {{\mathcal{T}}}
\def\cZ            {{\mathcal{Z}}}
\def\bbQ           {\mathbb{Q}}
\def\bbR           {\mathbb{R}}
\def\bbF           {\mathbb{F}}
\def\bbK           {\mathbb{K}}
\def\bbC           {\mathbb{C}}

\def\bbZ           {\mathbb{Z}}  
\def\lan{\langle}   \def\ran{\rangle}
\def\b1{\mathbf{1}}  \def\ri{\mathrm{i}}  \def\bo{\mathbf{0}} 
\def\fg{\ensuremath{\mathfrak{g}}}        \def\fC{\ensuremath{\mathfrak{C}}}
\def\la{\lambda}   \def\ka{\kappa}

\def\La{\Lambda}

\numberwithin{equation}{section}

\title{Exotic quantum subgroups and extensions of affine Lie algebra VOAs -- part I}

\author{Terry Gannon\\  {\footnotesize Department of Mathematics, University of Alberta,}\\
{\footnotesize Edmonton, Alberta, Canada T6G 2G1}\\
{\footnotesize e-mail: {\tt tjgannon@ualberta.ca}} }

\begin{document}

\maketitle \begin{abstract} Prototypical rational vertex operator algebras are  associated to affine Lie algebras $\fg^{(1)}$ at positive integer level $k$. They correspond physically to the Wess-Zumino-Witten theories, and their  representation theory can be captured by  quantum groups at roots of unity. One would like to identify the full (bulk) conformal field theories whose chiral halves are those VOAs. Mathematically, these correspond to their module categories. Until now, this has been done only for $\mathfrak{sl}(2)$ (famously, the A-D-E classification of Cappelli-Itzykson-Zuber) and $\mathfrak{sl}(3)$. The problem reduces to knowing the possible extensions of those VOAs, and the tensor equivalences between those extensions. Recent progress puts the tensor equivalences in good control, especially for $\mathfrak{sl}(n)$. This paper focuses on the extensions. We prove that, for any simple $\fg$, there is a bound $K(\fg)$  growing like rank$(\fg)^3$, such that for any  level $k>K(\fg)$, the only extensions are generic (i.e.\ simple-current ones). We use that bound to find all extensions for $\fg=\mathfrak{sl}(4)$ and $\mathfrak{sl}(5)$, at all levels, as well as all $\fg=\mathfrak{sl}(n)$ at levels $k\le 5$ (only those for $k=1$ had been classified before). In the sequel to this paper, we find all extensions for all simple $\fg$ of rank $\le 6$ (and the corresponding low level classifications). \end{abstract}

\section{Introduction}

The chiral halves of 2-dimensional conformal field theories (CFTs) are often identified with vertex operator algebras (VOAs) $\cV$. The friendliest examples, corresponding to so-called  rational theories, have a semisimple representation theory  Mod$(\cV)$ called a \textit{modular fusion category} (MFC) (also called a modular tensor category in the literature). 

The prototypical examples of rational VOAs and MFCs come from (finite-dimensional complex) simple Lie algebras $\fg$ at level $k\in\bbZ_{>0}$, which we denote by $\cV(\fg,k)$ and $\cC(\fg,k)=\mathrm{Mod}(\cV(\fg,k))$ respectively. These categories also  arise from quantum groups $U_q(\fg)$ at the appropriate root of unity $q$.  Their simple objects are parametrised by level $k$ highest weights $\la=(\la_1,...,\la_r)\in P_+^k(\fg)$ of the affine Lie algebra $\fg^{(1)}$ ($r$ is the rank of $\fg$; we'll usually drop as redundant the extended label $\la_0$); the tensor unit is $\bo=(0,0,\ldots,0)$.

The challenge of classifying all possible full (bulk) CFTs corresponding to a given pair of chiral halves is a fundamental one, going back to \cite{GW} and studied by many  in the 1990s and early 2000s (see e.g.\  \cite{EK,Oc,RFFS} and the recent review \cite{EKr}). The case of $\fg=A_1$ (for all levels $k$) falls into the famous A-D-E classification of Cappelli-Itzykson-Zuber \cite{CIZ}.  Following several insights by the subfactor community, a full CFT for a given rational VOA has been formulated categorically as a \textit{module category} \cite{Ost1} of the given MFC, which is the perspective we take. 

Given a  rational VOA $\cV$, each indecomposable module category for the MFC Mod$(\cV)$ is characterized by a triple $(\cV^e_L,\cV^e_R,\cF)$, where the VOAs $\cV^e_L$ and $\cV^e_R$ are both (conformal)  extensions of $\cV$, and $\cF:\mathrm{Mod}(\cV^e_L)\to\mathrm{Mod}(\cV^e_R)$ is a braided tensor equivalence  (see e.g.\ Corollary 3.8 of \cite{DNO}). 

In the physics literature this description of module categories appears as the \textit{naturality} of \cite{MS}. The extensions correspond to anyon condensations.
Physically, the trivial case $\cV^e_L=\cV^e_R=\cV(\fg,k)$ and $\cF=\,$Id corresponds to the Wess-Zumino-Witten model \cite{GW} of strings living on the compact simply connected Lie group $G$ corresponding to $\fg$. More generally, when the target $G$ is replaced with $G/Z$ for some subgroup $Z$ of the centre, the $\cV^e_{L,R}$ are so-called simple-current extensions of $\cV(\fg,k)$, and the equivalence $\cF$ is likewise  built from simple-currents. Simple-currents are the invertible  objects of $\cC(\fg,k)$, and the extensions and equivalences they generate are completely understood. In the VOA language, $\cV^e$ is a simple-current extension of $\cV$ iff $\cV$ is the orbifold $(\cV^e)^G$ (i.e.\ the fixed points of a $G$-action) by some (finite) abelian group $G$ of VOA automorphisms.

Braided tensor equivalences $\cF$ between MFC $\cC(\fg,k)$ of affine Lie algebra type are in relatively good control, thanks to Edie-Michell \cite{EM} (building on \cite{Gan3,GRW}), who completes this for $\fg=A_r,B_r,C_r$ and $G_2$ at all levels. The analogous work for the simple-current extensions of $\fg=A_r$ is completed in \cite{EM2}. Extending this to all $\fg$ should be accessible. The result is that almost all  equivalences come from simple-currents and outer automorphisms of $\fg$.

Thus the biggest challenge in the classification of module categories for the $\cC(\fg,k)$ is determining the possible extensions of $\cV(\fg,k)$. This is the question addressed in this paper. It is also the part of the module category classification which is of direct importance to the VOA community. Until this paper, the answer was known only for $\fg=A_1$ and $\fg=A_2$, though Ocneanu has announced \cite{Oc} progress for a couple other $A_r$ without supplying  proofs. 

We show that almost all extensions of the  VOAs $\cV(\fg,k)$ are simple-current extensions. By an \textit{exceptional extension} of $\cV(\fg,k)$ we mean one which is not a simple-current one. As we shall see, a relatively common source of exceptional  extensions are those of Lie type, where $\cV^e$ is also an affine Lie algebra VOA -- an example is $\cV(A_1,10)\subset\cV(C_2,1)$. These are also well-understood. 

There is a contravariant metaphor between finite groups and rational VOAs, where quotients of groups correspond to extensions of VOAs. In this metaphor, the results of this paper are analogous to the foundational group theoretic result that given a finite group of Lie type, e.g.\ $SL_n(\bbF_q)$, the quotient by the centre is almost always simple. The level corresponds to the finite field, and the simple-currents to the centre. Our result is analogous: after extending $\cV(\fg,k)$ by simple-currents, the resulting VOA is almost always \textit{completely anisotropic} (i.e.\ cannot be further extended).

Pushing this metaphor a little further, recall that the generic finite simple group is of Lie type. It is tempting then to guess that the MFC of a completely anisotropic VOA will generically be `of Lie type', whatever that may mean (presumably include arbitrary Deligne products of arbitrary Galois associates of $\cC(\fg,k)$). Compare this to Remark 6.6 in \cite{DMNO}.

The finite subgroups of SU(2) fall into an A-D-E classification, like the module categories for $A_1$. For this reason, Ocneanu \cite{Oc} introduced the term \textit{quantum subgroup}. More precisely, a \textit{quantum module} denotes any module category $(\cV_L^e,\cV_R^e,\cF)$ for the MFC $\cC(\fg,k)$, whilst a quantum subgroup denotes those of pure extension type $(\cV^e,\cV^e, \mathrm{Id})$. The latter are in natural bijection with the \textit{connected \'etale} (or \textit{condensable}) \textit{algebras} $\cA$ of $\cC(\fg,k)$.

The basic combinatorial invariant of a MFC is the \textit{modular data}, a unitary matrix representation $R$ of the modular group  SL$_2(\bbZ)$. It is generated by a symmetric matrix $S=R\left({0\atop 1}{-1\atop 0}\right)$ (coming from the Hopf link) and a diagonal matrix $T=R\left({1\atop 0}{1\atop 1}\right)$ (coming from the ribbon twist). For $\cC(\fg,k)$, the $S$-entries come from Lie group characters evaluated at elements of finite order, and the $T$-entries from the second Casimir.

\cite{NRW} classified the possible modular data coming from MFC of rank $\le12$. The most exotic of these (no.\ 139 in Appendix E for rank 11) is not yet known to be realized by a MFC. If it is, it would be,   to this author's knowledge,  the only known completely anisotropic MFC not of Lie type (recall the conjecture given earlier this section).

Fix a quantum subgroup of Mod$(\cV)$, i.e.\ an extension $\cV^e$ of $\cV$, and let $M$ be any simple $\cV^e$-module. Then $M$ restricts to a $\cV$-module
$\mathrm{Res}\,M\cong \oplus_\mu B_{M, \mu }\,\mu$,  where the $\mu$ are simple $\cV$-modules and the multiplicities $B_{M,\mu}\in\bbZ_{\ge 0}$ form the \textit{branching rule matrix}. Note that $B_{M,\bo}=\delta_{M,\bo^e}$, where $\bo=\cV$ and $\bo^e=\cV^e$ are the tensor units in the MFC Mod$(\cV)$ and Mod$(\cV^e)$ respectively. The branching rules intertwine the modular data of $\cV$ and $\cV^e$:
\begin{equation}\label{modinv}BS=S^eB\,,\ \ BT=T^eB\,.\end{equation}

Given any quantum module $(\cV_L^e,\cV_R^e,\cF)$, the \textit{modular invariant} is  $\cZ=B_R^t\Pi B_L$ where $\Pi$ is the permutation matrix implicit in $\cF$, and $B_L,B_R$ are the branching matrices for $\cV\subset\cV_L^e$ resp.\ $\cV\subset\cV_R^e$.
Then we see from \eqref{modinv} that $\cZ$ commutes with $S$ and $T$.

The physics literature (e.g.\ \cite{CIZ}) emphasized these modular invariants, expressed as the formal generating function (\textit{partition function}) $\sum_{\la,\mu}\cZ_{\la,\mu}\chi_\la\overline{\chi_\mu}$. But inequivalent quantum modules can have identical modular invariants (e.g.\ at $A_2$ level 3), and seemingly healthy modular invariants may not  be realized by a quantum module (e.g.\  at $B_4$ level 2).
 Nevertheless,  there is close to a bijection between the modular invariants and the module categories, at least for $\fg$ of small rank. This coincidence  though  plays no role in our analysis.

The restriction functor has an adjoint, called \textit{(alpha-)induction}. More precisely, each extension $\cV^e$ of $\cV$ corresponds to a commutative  (\textit{\'etale}) algebra object $\cA$ in $\cC=\mathrm{Mod}(\cV)$, which we'll also call a quantum subgroup: the object is $\cA=\mathrm{Res}\,\cV^e$ and multiplication $\mu\in\mathrm{Hom}_\cC(\cA\otimes\cA,\cA)$ comes from the vertex operator structure $Y^e(u,z)v$ on $\cV^e$. A  quantum subgroup $\cA$ of $\cC$ is exceptional iff  $\cA$ is not a direct sum of invertible objects (simple-currents) in $\cC$.  In terms of induction, the modular invariant can be written $\mathcal{Z}_{\la,\mu}=\mathrm{dim\,Hom}_{\mathrm{Mod}\,\cV^e}(\mathrm{Ind}\,\la,\mathrm{Ind}\,\mu)$. 
 Induction is a tensor functor from $\cC$ to the fusion category Rep$_\cC\,\cA$ of $\cA$-modules (by contrast, for groups the tensor functor is restriction). In terms of the contravariant metaphor, $\cA$-modules for $\cV^e$ are the analogues here of projective representations for a group $G$: just as projective representations become ordinary representations for some extension of $G$, so these $\cA$-modules become ordinary VOA modules for the subVOA $\cV$. Restriction Rep$_\cC\,\cA\to\cC$ is the forgetful functor, forgetting the multiplication by $\cA$.
The MFC Mod$(\cV^e)$ is the full subcategory of \textit{local} or \textit{dyslectic} $\cA$-modules, which we denote by (Rep$_\cC\,\cA)^{\mathrm{loc}}$. Under our metaphor,  local $\cA$-modules correspond to true $G$-representations. Restriction limited to (Rep$_\cC\,\cA)^{\mathrm{loc}}$ is what we call $B$. 

The following observation is the starting point of how induction constrains quantum subgroups, in the affine Lie algebra categories $\cC(\fg,k)$. Let $\la^*$ denote the dual of $\lambda$ (more precisely, $\la^*$ is the highest weight of the contragredient $\fg^{(1)}$-module $L(\la)^\vee$), $\rho$  the Weyl vector of $\fg$, and $\ka=k+h^\vee$ the shifted level. Write $(\la|\mu)$ for the usual invariant bilinear form of $\fg$.

\medskip\noindent\textbf{Ocneanu's Lemma.} (Lemma 4.1)  \textit{Let $\cA=\cV^e$ be a quantum subgroup of $\cC(\fg,k)$. Let $\la\in P_+^k(\fg)$  satisfy $(\la+\la^*+2\rho|\la+\la^*)<2\ka$. Then} Ind$\,\la$ \textit{is a simple $\cA$-module.}\medskip

It should be repeated that the area has been profoundly influenced by  the subfactor community, including Ocneanu. See e.g.\ \cite{EK,Oc,Oc2} for an abundance of ideas. In particular, 
Ocneanu realised that the Lemma implies that, for any $\fg$, there is a threshold $K(\fg)$ such that, for any level $k>K(\fg)$, no extension of $\cC(\fg,k)$ is exceptional \cite{Oc}. Unfortunately the details of his argument never appeared, nor did an estimate for $K(\fg)$.   To our knowledge, the first published proof of Lemma 4.1 appears implicitly in \cite{Sch} (see also \cite{Xu}). 
In that important paper, Schopieray used it to find explicit estimates for $K(\fg)$ when $\fg$ is rank 2: in particular, he obtains $19\,896$ for $C_2$ and $18\,271\,135$ for $G_2$. His method is expected to extend to all Lie algebras, with an estimate for $K(\fg)$  growing exponentially with the square of the rank of $\fg$, though a proof for $\fg$ with rank$\,>2$ is lacking. These bounds are far  too high however to be used in classification efforts, even for $G_2$. 

Fortunately Schopieray's estimates can be tightened significantly. Here we couple Ocneanu's Lemma to the Galois symmetry of modular data \cite{CG}, to find $\la$ such that Ind$\,\la$ is both simple and local. This gives much tighter thresholds.
Our main result, proved in Section 4, is:

 \medskip\noindent\textbf{Main Theorem.}   \textit{For any simple $\fg$ there exists a threshold $K(\fg)$ such that, for any level $k>K(\fg)$, any  quantum subgroup  of $\cC(\fg,k)$ must be a simple-current extension. As $r=\mathrm{rank}(\fg) \to\infty$, $K(\fg)\in O(r^{3+\epsilon})$ for all $\epsilon>0$.}\medskip
 
Actually, this argument shows that a quantum subgroup can be exceptional only if the level lies in a sparse set;  $K(\fg)$ is the maximum of that set. For example, for $\fg=A_r$ ($r$ even), the size of that sparse set is $\ll 5r^3/4$. A little more work (Step 2 of Section 4.2) reduces that number considerably. 
 Table 1.1 displays the resulting set of suspicious levels for each simple Lie algebra $\fg$ of rank at most 6.  Compare our thresholds $K(C_2)=12$ and $K(G_2)=4$ with those of Schopieray given above. The numbers in boldface are the levels which actually possess an exceptional quantum subgroup. In the table,  the dash denotes consecutive numbers whilst the dots denote an arithmetic sequence: e.g.\ `20--22' denotes `20,21,22',  and `21.23..55' denotes `21,23,25,27,...,55'.  We find that the number of such suspicious  levels for each $\fg$ actually grows like rank$(\fg)^2$.

 \medskip
\centerline{
{\scriptsize\vbox{\tabskip=0pt\offinterlineskip
  \def\tablerule{\noalign{\hrule}}
  \halign to 6.2in{
    \strut#&
    \hfil#&\vrule#&\hfil#&\vrule#&\vrule#&    
\hfil#&\vrule#&\hfil#&\vrule#&\hfil#&\vrule#&\hfil#&\vrule#&\hfil#
\tabskip=0pt\cr
&$\fg\,\,\,$&\,&&\hfill levels $k$\hfill&&\,total\,\cr
\tablerule&$A_1\,$&\,&&\hfill $\mathbf{10},\,\mathbf{28}$\hfill&&\hfill2\hfill\cr
\tablerule&$A_2\,$&\,&&\hfill $\mathbf{5},\,\mathbf{9},\,\mathbf{21},\,57$\hfill&&\hfill4\hfill\cr
\tablerule&$A_3\,$&\,&&\hfill $\mathbf{4},\,\mathbf{6},\,\mathbf{8},\,10,\,11,\,12,\,14,\,16,\,18,\,20,\,26,\,32,\,38,\,86$\hfill&&\hfill14\hfill\cr
\tablerule&$A_4$\,&\,&&\hfill $\mathbf{3},\,\mathbf{5},\,\mathbf{7},\,9,10,11,13,15,17,19,21,23,25,31,35,37,43,49,55,85,115$\hfill&&\hfill21\hfill\cr
\tablerule&$A_5\,$&\,&&\hfill $\,\mathbf{4},\,\mathbf{6},\,\mathbf{8},9,10,12,14\!-\!16,18,20\!-\!22,24,26,27,30,32,33,36,38,40,42,48,54,60,72,78,84\,$\hfill&&\hfill29\hfill\cr
\tablerule&$A_6\,$&\,&&\hfill $\mathbf{5},\mathbf{7},8,\mathbf{9},11,13,15,17\!-\!19,21.23..55,59,65,71,77,83,95,113,173$\hfill&&\hfill36\hfill\cr
\tablerule\tablerule&$B_3\,$&\,&&\hfill $\mathbf{5},\,\mathbf{9},\,19,\,25$\hfill&&\hfill4\hfill\cr
\tablerule&$B_4$\,&\,&&\hfill $\mathbf{2},\,\mathbf{7},\,8,\,\mathbf{11},\,13,\,14,\,17,\,21,\,22,\,23,\,26,\,29,\,53$\hfill&&\hfill13\hfill\cr
\tablerule&$B_5\,$&\,&&\hfill $\mathbf{9},11,12,\mathbf{13},15,19\!-\!21,24,27,30,33,39,51$\hfill&&\hfill14\hfill\cr
\tablerule&$B_6\,$&\,&&\hfill $4,7,9,\mathbf{11},13,\mathbf{15},17,19,22,24,25,28,29,31,34,37,43,49,55,67,79$\hfill&&\hfill22\hfill\cr
\tablerule\tablerule&$C_2\,$&\,&&\hfill $\mathbf{3},\,\mathbf{7},\,\mathbf{12}$\hfill&&\hfill3\hfill\cr
\tablerule&$C_3\,$&\,&&\hfill $\mathbf{2},\,\mathbf{4},\,\mathbf{5},\,6,\,7,\,9\!-\!11,\,14,\,17,\,26$\hfill&&\hfill11\hfill\cr
\tablerule&$C_4$\,&\,&&\hfill $\mathbf{3},\,\mathbf{5},\,6,\mathbf{7},\,8,9,\mathbf{10},11,13,16,18,19,22,23,25.28..43,55$\hfill&&\hfill22\hfill\cr
\tablerule&$C_5\,$&\,&&\hfill $\mathbf{3},\mathbf{4},\mathbf{6},8\!-\!12,14\!-\!16,18,19,21\!-\!24,27,33,36,39,42,54$\hfill&&\hfill23\hfill\cr
\tablerule&$C_6\,$&\,&&\hfill $3,4,\mathbf{5},\mathbf{7},8\!-\!15,17\!-\!23,25\!-\!29,32,35,38,41,44,48,53,56,59,68,98$\hfill&&\hfill35\hfill\cr
\tablerule\tablerule&$D_4$\,&\,&&\hfill $\mathbf{6},\,9,\,\mathbf{10},\,\mathbf{12},\,14,\,16,\,18,\,22,\,24,\,26,\,28,\,30,\,36,\,42,\,54,\,60$\hfill&&\hfill16\hfill\cr
\tablerule&$D_5\,$&\,&&\hfill $\mathbf{4},6,7,\mathbf{8},10,\mathbf{12},14.16..24,28,32,34,36,40,42,46,48,50,52.58..88,112$\hfill&&\hfill29\hfill\cr
\tablerule&$D_6\,$&\,&&\hfill $\,\,\mathbf{5},6,\mathbf{8},\mathbf{10},11,12,\mathbf{14},17,18,20,22,23,25,26,28-30,\!32,\!34,\!35,\!36,\!38,\!40,\!44.46..50,\!54,\!56,\!60,\!62.68..116,128,\!170,\!200\,\,$\hfill&&\hfill43\hfill\cr
\tablerule\tablerule&$E_6\,$&\,&&\hfill $\mathbf{4},\mathbf{6},8,10,\mathbf{12},14.16..48,54,58,60.66..120,138,168,198$\hfill&&\hfill39\hfill\cr
\tablerule&$F_4\,$&\,&&\hfill $\mathbf{3},\,\mathbf{6},\,\mathbf{9},11,12,13,15,17,19,20,21,24,26,30,51$\hfill&&\hfill15\hfill\cr
\tablerule&$G_2$\,&\,&&\hfill $\mathbf{3},\,\mathbf{4}$\hfill&&\hfill2\hfill\cr
\noalign{\smallskip}}}}}

\centerline {{\bf Table 1.1.} The possible  levels $k$ for exceptional quantum subgroups}\medskip

Most  levels appearing in Table 1.1 fail to possess an exceptional quantum subgroup. Our main strategy for eliminating a given suspicious level
is modular invariance \eqref{modinv}, which gives us the inequality
\begin{equation} \sum_{\la}Z_{\la}S_{\la,\mu} \ge 0\label{rhoineq}\end{equation}
valid  for all $\mu\in P_{+}^{\,k}(\fg)$, where we write $\cA=\bigoplus_\la Z_\la\la$ ($Z_\lambda=B_{\bo^e,\la}$ are the multiplicities).

 Building on this, in Section 5  we identify all exceptional quantum subgroups for $\fg=A_3$ and $A_4$, and work out their branching rules $B_{M,\mu}$. Exploiting level-rank duality this also yields the exceptional quantum subgroups for all $\fg=A_r$ when the level is $k\le5$ (see Section 5.4). The results are collected in Tables 2.1 and 2.2.  The main tools for  finding these branching rules are the formula
Res(Ind$(\nu))=\cA\otimes \nu$, and  modular invariance \eqref{modinv}.  In the follow-up paper \cite{Gan-ii}, we find all quantum subgroups and branching rules for all simple $\fg$ of rank $\le 6$. Combining this with the work of Edie-Michell \cite{EM,EM2,EM3}, this gives the module category (quantum module) classification for  $\fg=A_r$ \cite{EM3} for $r\le 6$
(see also Section 6).

We find in \cite{Gan-ii} that
there are precisely 75 exceptional quantum subgroups of the $\cC(\fg,k)$ when $\fg$ has  rank $\le 6$. Of these, 
59  are the aforementioned {Lie-type conformal extensions}.  In fact, Table 2.1 shows that all exceptional extensions for $\fg=A_1,A_2,A_3,A_4$ are of Lie type. 
Six of the 75 are simple-current extensions of ones of Lie type.
Three quantum subgroups can be interpreted using so-called  mirror extensions or level-rank duals of Lie-type extensions.
Only seven quantum subgroups of those $75$ (namely, at $(A_6,7),(C_4,10),(D_4,12)$, $(D_5,8),(D_6,5),(E_6,4),(F_4,6)$) are truly exotic and require ad hoc methods to construct. 


We also address uniqueness: for a given object $\cA$ in $\cC(\fg,k)$, how many inequivalent \'etale algebra structures can be placed on it? Most importantly, Theorem 5.2 addresses uniqueness for the $\cA$ appearing in the Lie-type extensions of arbitrary rank.

Incidentally, none of the MFC Mod$(\cV^e)=(\mathrm{Rep}_{\cC(\fg,k)}\cA)^{\mathrm{loc}}$ of these exceptional extensions (including the ones appearing in the $k\le 5$ classification) are in any way exotic. On the other hand, there are MFC we would call exotic (e.g.\ the doubles of the fusion categories of Haagerup type \cite{EG1}), but corresponding VOAs still haven't been found though are conjectured to exist.

 Surely, when rank $>6$, most  exceptional quantum subgroups will continue to come from Lie-type conformal embeddings. We expect that, given any simple $\fg$, there are at most four or five  levels where all their exceptional quantum subgroups reside. For example, $\fg=E_7$ appears to have exactly three exceptional quantum subgroups (namely at $k=3,12,18$, two of which are of Lie type) and $\fg=E_8$ appears to have only one (at $k=30$, of Lie type). In short, apart from simple-current extensions and Lie-type conformal embeddings (both of which are completely understood) and their level-rank duals, all other extensions of $\cV(\fg,k)$ appear to be extremely rare, with very few in each rank, perhaps only finitely many as one runs over all simple $\fg$. 

The quantum subgroups when $\fg$ is semisimple largely reduce to those of the simple summands, using a contravariant analogue of Goursat's Lemma. More precisely, Theorem 3.6 in \cite{DNO} says that the connected \'etale algebras of the Deligne product MFC $\cC\boxtimes\cD$ correspond naturally to choices of connected \'etale algebras $\cA_c,\cA_d$ of $\cC$ resp.\ $\cD$, fusion subcategories $\cF_c,\cF_d$ of  (Rep$_\cC\,\cA_c)^{\mathrm{loc}}$ and (Rep$_\cD\,\cA_d)^{\mathrm{loc}}$ respectively, and a braided tensor equivalence $\phi:\cF_c^{\mathrm{rev}}\to\cF_d$. We learn here that when $\cC=\cC(\fg,k)$, the categories  (Rep$_\cC\,\cA)^{\mathrm{loc}}$ for any connected \'etale algebra  $\cA$ are themselves usually (always?) minor tweaks of other $\cC(\fg',k')$; we learn in \cite{Sa} that with few exceptions the fusion subcategories of $\cC(\fg',k')$ are associated to simple-currents in standard ways; we know from \cite{EM,EM2} to expect that the braided equivalences between these subcategories are built from simple-currents and outer automorphisms of the Lie algebras. We return to semisimple $\fg$ in \cite{Gan-ii}.

The Galois associate $\cC^\sigma$ of a MFC $\cC$ is itself a MFC. The quantum subgroups (and quantum modules) of $\cC$ and $\cC^\sigma$ are in natural bijection. For example,  the MFC $\cC(F_4,3)$ is a Galois associate of $\cC(G_2,4)$, and they each have two quantum subgroups. A more important example of this phenomenon are the family of rational VOAs constructed in e.g.\ \cite{ACL}: most of those MFCs are  simple-current extensions of the Deligne product of Galois associates of MFC of affine Lie type, so their quantum subgroups are now also classifiable.

The connected \'etale algebras  (quantum subgroups)  are known for two other large classes of categories. 
When the MFC $\cC$ is pointed (i.e.\ all simples are simple-currents), these algebras obviously must  be simple-current extensions and so are under complete control. The MFC (Rep$_\cC\,\cA)^{\mathrm{loc}}$ of the extension will also be pointed. Much more interesting, the connected \'etale algebras of the twisted Drinfeld doubles $\cC=\cZ(\mathrm{Vec}_G^\omega)$ of finite groups are known \cite{DS,EG}. The MFC (Rep$_\cC\,\cA)^{\mathrm{loc}}$ of these extensions are those of another twisted double $\cZ(\mathrm{Vec}_K^{\omega'})$, where $K$ will be some subquotient of $G$, and the 3-cocycle $\omega'$ may be nontrivial even if $\omega$ is trivial \cite{EG}.

Incidentally, the Galois symmetry (more precisely, the condition $\epsilon_\sigma(\la)=\epsilon_\sigma(\mu)$ of Section 3.5) for $\fg=A_2$ also plays a fundamental role in the study of simple factors of the Jacobians of Fermat curves \cite{KR} -- see \cite{BCIR}
for a discussion. Recall also the A-D-E classification of $A_1$ quantum modules  \cite{CIZ,Ost1}, and the natural identification of module categories for lattice VOAs $\cV_L$ of rank $d=\mathrm{dim}\,L$ with rational points in the dual of the Grassmannian $G_{d,d}(\bbR)$. This suggests the intriguing possibility of other relations between quantum subgroups (and quantum modules) and geometry. We hope this paper, as well as  \cite{Gan-ii,EM3}, can inspire others to identify further relations.

 \subsection{Notation and terminology}{} 
 
 $\xi_n=e^{2\pi\ri/n}$
 
 
  \noindent $\overline{z}$ \quad complex-conjugate
 
 \noindent $\fg$ \quad simple finite-dimensional Lie algebra over $\bbC$

\noindent $h^\vee$ \ dual Coxeter number

 \noindent $k$ resp.\ $\kappa=k+h^\vee$ \  level resp.\ shifted level ($k\in\bbZ_{>0}$)
 
 \noindent MFC \ modular fusion (tensor) category (\S3.1)
 
 \noindent $\cC(\fg,k)$ resp.\ $\cV(\fg,k)$ \  the MFC resp.\ rational VOA associated to $\fg$ at level $k$
 
 \noindent $\cV$ \ a  rational VOA (\S3.2)
 
\noindent Mod$(\cV)$ \ the category of modules of $\cV$, a MFC

\noindent $\la,\mu,...$ \ highest weights, labelling simple objects in $\cC(\fg,k)$ 

\noindent $\la_0$ \ extended Dynkin label (eq.(2.3))
 
\noindent $\rho$ \ Weyl vector

\noindent $ P^{\,k}_{+}(\fg),P^{\,\ka}_{\!++}(\fg)$ \  level $k$ highest weights unshifted/shifted by $\rho$ (eq.(2.1)) 


 \noindent $\bo$ or $\b1$ \ tensor unit: $\bo=(0,0,...,0)\in P_+^k(\fg)$ or $\b1=(1,1,...,1)\in P_{\!++}^{\,\ka}(\fg)$
 
 \noindent $\la^*$ \ dual or contragredient
 
\noindent $C$ \ permutation $\la\mapsto\la^*$ of $P_{\!++}^{\,\ka}(\fg)$ or $P_+^k(\fg)$

\noindent $J$ \ simple-current, i.e.\ invertible simple object
 
\noindent $\cJ(\fg,k)$\ the group of simple-currents for $\cC(\fg,k)$

\noindent $\lan \la,\mu,...\ran_\tau$ \ orbit of $\la,\mu,...$ under permutation(s) $\tau$ (usually $\tau=C$ or $J$)
 
 \noindent $S,T$ \   modular data matrices normalised to get SL$_2(\bbZ)$-representation
 
\noindent $\theta(\la)=T_{\la,\la}\overline{T_{\b1,\b1}}$ \ ribbon twist

\noindent $h_\la$ conformal weight (\S2.1)
 
 \noindent qdim$(\la)$ \ quantum-dimension $S_{\la,\b1}/S_{\b1,\b1}$
 
\noindent $\varphi_J(\la)$ \ grading associated to simple-current $J$ (eq.\eqref{scform})

\noindent $(\la|\mu)$ \ standard inner product of (co)weights


 
 
 \noindent $Z_\cC(\la)$ \ centralizer (\S3.1)
 
 \noindent $\cC_0$ the adjoint subcategory $Z_\cC(\cJ(\cC))=\{\la\,|\,\varphi_J(\la)=1\,\forall J\}$ (\S4.1)
 
\noindent $\cV\subset\cV^e$ \ conformal extension of  rational VOAs (\S3.2)

 \noindent Lie-type conformal extension \ (\S3.2)

\noindent $\cA=\bigoplus_\la Z_\la\la$ \ quantum subgroup=connected \'etale algebra (\S3.3)

\noindent Rep$_\cC\,\cA$ \ fusion category of $\cA$-modules in $\cC$ (\S3.3)

\noindent $(\mathrm{Rep}_\cC\,\cA)^{\mathrm{loc}}=\mathrm{Mod}(\cV^e)$ \ MFC of local (dyslectic) $\cA$-modules (\S3.3)

\noindent $M$ \ an $\cA$-module, usually local and simple
 
\noindent $B_{M,\mu}$ \ branching rules=restriction $\mathrm{Mod}(\cV^e)\to\mathrm{Mod}(\cV)$ (\S3.3)

\noindent Res \ restriction, usually from Mod$(\cV^e)\to\mathrm{Mod}(\cV)$

\noindent Ind \ induction Mod$(\cV)\to\mathrm{Rep}_{\mathrm{Mod}\,\cV}\,\cA$

\noindent $\cZ=B^tB$ \ modular invariant for quantum subgroup

\noindent $\epsilon_\sigma(\la),\,\la^\sigma$ \  Galois parity $\pm 1$ and permutation (\S3.5)

\noindent $\fC_k(\fg)$ \ candidates $\la$  (Definition 3.10)


 \noindent $f_\fg\in\{1,2,3\}$ \ the smallest $f>0$  for which $f\rho$ is sum of coroots (\S2.1)
 
 \section{The exceptional  quantum subgroups}
 
 This section gives all quantum subgroups in the MFC $\cC(\fg,k)$, equivalently all conformal extensions of the rational VOA $\cV(\fg,k)$, when $\fg=A_1,A_2,A_3,A_4$ $\forall k\ge1$, as well as when $k\le5$ for $\fg$ being any $A_r$. In the process we fix the Lie theoretic notation used throughout this paper and its sequel \cite{Gan-ii}. The proof occupies Sections 3,4,5.
 
 \subsection{Lie theoretic background and notation}
 
 Let $\fg$ be a simple Lie algebra over $\bbC$. See \cite{Kac} for the basic theory of the affine Lie algebra $\fg^{(1)}$, which underlies much of the following material. The level $k$ highest weights of  $\fg^{(1)}$ (which parametrise up to equivalence the simple objects in the MFC $\cC(\fg,k)$, and equivalently the simple modules of the VOA $\cV(\fg,k)$) form the set 
 \begin{equation}\label{P+}P_{+}^{k}(\fg)=\left\{(\la_1,\la_2,\ldots,\la_r)\in\bbZ_{\ge0}^r\,|\, \sum_{i=1}^ra_i^\vee\la_i\le k\right\}\end{equation}
 where $a_i^\vee$ are the co-marks of $\fg$ and $r$ its rank. We also write  $\sum_{i=1}^r\la_i\La_i$ for $(\la_1,\la_2,\ldots,\la_r)$, where $\La_i$ are the fundamental weights. For large $k$, 
 \begin{equation}\|P_{+}^{k}(\fg)\|\approx \frac{k^r}{r!a_1^\vee\cdots a_r^\vee}\,.\label{growthP+}\end{equation}
Formulas below are often simplified by shifting $\la$ by  the Weyl vector $\rho=(1,1,...,1)$ and shifting the level by the dual Coxeter number, $\kappa=k+h^\vee$: hence  define
$P_{\!++}^{\,\ka}(\fg)=\rho+P_+^k(\fg)$.
 We denote the tensor unit in $\cC(\fg,k)$ by $\bo=(0,0,...,0)\in P_+^k(\fg)$ and equivalently $\b1=(1,1,...,1)\in P_{\!++}^{\,\ka}(\fg)$.
The dual (contragredient) of $\la$ is denoted $\la^*$ or $C\la$.

Scale the invariant bilinear form $(\la|\mu)$ on $\fg$ so the highest root has norm-squared 2, and use it to identify the Cartan subalgebra and its dual, as in Section 6.2 of \cite{Kac}. Then the roots and coroots are related by $a_i\alpha_i=a_i^\vee\alpha_i^\vee$, and the
 fundamental weights $\Lambda_i$ satisfy $(\Lambda_i|\alpha_j^\vee)=\delta_{i,j}$. For $\la\in P_{\!++}^{\,\ka}(\fg)$ write
 \begin{equation}\la_0=\kappa-\sum_{i=1}^ra_i^\vee\la_i\,.\label{lambda0}\end{equation}
 
The central charge of $\cC(\fg,k)$  is $c=k\,\mathrm{dim}(\fg)/\ka$. Its modular data is
\begin{align}T_{\la,\mu}=\,&\delta_{\la,\mu}e^{2\pi\ri(h_\la- c/24)} \cr
S_{\lambda,\mu}=\,&x\sum_{w\in W}\mathrm{det}(w)\exp\left(-2\pi \ri\,\frac{(\lambda| w(\mu))}{\kappa}\right)\label{Sgenform}\end{align}
for $\la,\mu\in P_{\!++}^{\,\kappa}(\fg)$, where  $h_\la=((\la|\la)-(\rho|\rho))/(2\ka)$ is  the \textit{conformal weight}, $W$ is the Weyl group of $\fg$ and $x$ is the square-root of a positive rational number independent of $\la,\mu$. Next subsection we spell this out in full detail for $\fg=A_r$. The ribbon-twist is $\theta(\la)=T_{\la,\la}\overline{T_{\b1,\b1}}=e^{2\pi\ri h_\la}$.
Useful is 
\begin{equation}\label{strange}(\rho|\rho)=h^\vee \mathrm{dim}(\fg)/12\,.\end{equation}
The $r^3$ growth in the Main Theorem ultimately comes from $(\rho|\rho)$, as we see in the proof of Proposition 4.9(b) below.

For any simple $\fg$, define $f_\fg$ to be the smallest positive integer $f$  for which $f\rho$ is a sum of coroots. In particular, $f_{G_2}=3$, $f_{E_6}=f_{E_8}=1$, $f_{A_r}=1$ when $r$ is even, $f_{D_r}=1$ when $r\equiv 0,1$ (mod 4), and $f_\fg=2$ otherwise.

A \textit{simple-current} $J$ is an invertible object in $\cC(\fg,k)$. They are classified in \cite{F}. 
With one exception ($E_8$ at $k=2$, which plays no role in this paper), the simple-currents form the group called $W_0^+\cong Q^*/Q^\vee$ in Section 1.3 of \cite{KW}, which is isomorphic to the centre of the simply-connected connected compact Lie group of $\fg$. 
$W_0^+$ can be identified with  a subgroup of  symmetries of the affine Dynkin diagram, and so acts on $P_{\!++}^{\,\ka}(\fg)$ by permuting the Dynkin labels (including $\lambda_0$) of $\la$. This coincides with 
 the fusion product $\la\mapsto J\la:=J\otimes \la$. We denote $W_0^+$ by $\cJ(\fg,k)$. Our Main Theorem shows that subgroups of $\cJ(\fg,k)$  are the source of almost all quantum subgroups of $\cC(\fg,k)$. 

In $\cC(\fg,k)$ as in any pseudo-unitary MFC, associated to a simple-current $J$  are roots of unity $\varphi_J:P_{\!++}^{\,\ka}(\fg)\to\bbC^\times$ satisfying 
\begin{equation}S_{\la,J\mu}=\varphi_J(\la)\,S_{\la,\mu}\,,\qquad \theta(J\la)=\overline{\varphi_J(\la)}\,\theta(J)\,\theta(\la)\,.\label{scform}\end{equation}
Each  $\varphi_J$ is a grading of the fusion ring of $\cC(\fg,k)$. Further background is provided in Section 3.

 \subsection{The $A$-series}
 
 Let $\fg=A_r$. Write $r'=r+1$ and $\ka=k+r'$. 
 The dual  of $\la$ is $\la^*=C\la=(\la_r,\la_{r-1},\ldots,\la_1)$.  We have   $\la_0=\ka-\sum_{i=1}^r\la_i$. The group $\cJ(A_r,k)$ of simple-currents is cyclic of order $r'$, generated by $J_a=k\Lambda_1+\rho\in P_{\!++}^{\,\ka}(A_r)$ with permutation $J_a\la=(\la_0,\la_1,...,\la_{r-1})$. The grading $\varphi_J$ in \eqref{scform} for $\la\in P_{\!++}^{\,\kappa}(\fg)$ is $\varphi_{J_a}(\la)=\xi_{r'}^{t(\la-\rho)}$  where  $t(\la)=\sum_{i=1}^ri\la_i$. Also, 
 \begin{align}\theta(J^\ell_a)&\,=\xi_{2r'}^{ka(r'-a)}\,,&\cr t(J_a^i\la)&\,\equiv ki+t(\la)\ (\mathrm{mod}\ r')\,.\label{tjmu}\end{align}
 
 By $\lan\la,\mu,\ldots\ran_d$ we mean the direct sum of all objects in the orbits of $\la,\mu,...$, over the subgroup $\lan J_a^{r'/d}\ran$ of $\cJ$ of order $d$. By $\lan\la,\mu,\ldots\ran_{dc}$ we mean the direct sum using instead the dihedral group $\lan J_a^{r'/d},C\ran$ of order $2d$. 

 For any divisor $d$ of $r'$, if $d$ is odd and $d^2$ divides $kr'$, or $d$ is even and $2d^2$ divides $kr'$, there is a simple-current extension with algebra $\cA$ generated by $J_a^{r'/d}$, i.e.\ $\cA=\lan \b1\ran_d$. See Section \ref{sect:sc} for more details. 
 
We have $(\La_i|\La_j)=\frac{i\,(r'-j)}{r'}$, provided $i\le j$. Hence \begin{equation}(\la|\mu)=\sum_i\la[i]\mu[i]-t(\la)t(\mu)/r'\,,\label{Anorm}\end{equation} if we write  $\la[i]=\sum_{j=i}^r\la_i$. An effective formula for the $S$-matrix is
\begin{equation}S_{\la,\mu}=\frac{\ri^{rr'/2}}{\ka^{r/2}\sqrt{r'}}\xi_{\ka r'}^{t(\la)\,t(\mu)}\,\mathrm{det}(\xi_{\ka}^{-\la[i]\,\mu[j]})_{1\le i,j\le r'}\,,\nonumber\end{equation}
in terms of the determinant of the $r'\times r'$ matrix with $(i,j)$-entry $\xi_{\ka}^{-\la[i]\,\mu[j]}$.

The exceptional quantum subgroups for $\cC(A_r,k)$ when $r\le 4$ are given in Table 2.1, using $P_+^k(A_r)$ rather than $P^{\,\ka}_{\!++}(A_r)$ for readability. We find that all exceptional quantum subgroups for these $A_r$ are of Lie type. All of these extensions are unique.

 \medskip
\centerline{\scriptsize{\vbox{\tabskip=0pt\offinterlineskip
  \def\tablerule{\noalign{\hrule}}
  \halign to 6.4in{
    \strut# &
    \tabskip=0em plus1em &    
    \hfil#&
    #&\hfil#&\vrule\vrule#&    
\hfil#&\vrule#&\hfil#&\vrule#
\tabskip=0pt\cr
&&$\fg$&&\hfill $k$\hfill && 
\hfill MFC\hfill&&\hfill branching rules\hfill&\cr
\tablerule\tablerule&&$A_1\,$&&\hfill$\,10$\hfill&&\hfill$\,C_{2,1}$\hfill&&\hfill$\,\bo^e\mapsto\bo\oplus(6)\,,\ \La_1^e\mapsto(3)\oplus(7)\,,\ \La_2^e\mapsto J_a(\bo\oplus(6))=(4)\oplus(10)$\hfill&\cr
\tablerule&&&&\hfill$\,28$\hfill&&\hfill$G_{2,1}$\hfill&&\hfill$\bo^e\mapsto\lan\bo,(10)\ran_2=(0)\oplus(10)\oplus(18)\oplus(28)\,,\ \La_2^e\mapsto\lan(6),(12)\ran_2=(6)\oplus(12)\oplus(16)\oplus(22)$
\hfill&\cr
\tablerule\tablerule&&$A_2\,$&&\hfill5\hfill&&\hfill$A_{5,1}$\hfill&&\hfill$\Lambda_{2j}^e\mapsto J_a^j(\bo\oplus(2\,2))\,,\ \Lambda_{1+2j}^e\mapsto J_a^j((0\,2)\!\oplus\!(3\,2)),\ \ j=0,1,2$\hfill&\cr
\tablerule&&&&\hfill9\hfill&&\hfill$E_{6,1}$\hfill&&\hfill$\bo^e\mapsto\lan\bo,(4\,4)\ran_3,\ \Lambda_1^e,\Lambda_5^e\mapsto\lan(2\,2)\ran_3$\hfill&\cr
\tablerule&&&&\hfill\,21\hfill&&\hfill $E_{7,1}$\hfill&&\hfill$\bo^e\mapsto\lan\bo,(4\,4),(6\,6),(10\,10)\ran_3, \ \Lambda_6^e\mapsto\lan(0\,6),(4\,7)\ran_{3c}$\hfill&\cr
\tablerule\tablerule&&$A_3$&&\hfill4\hfill&&\hfill$B_{7,1}$\hfill&&\hfill$\bo^e\mapsto \lan\bo,(012)\ran_2,\ 
\La_1^e\mapsto J_a \lan\bo,(012)\ran_2,\ \La_7^e\mapsto 2\cdot(111)$\hfill&
\cr\tablerule&&&&\hfill6\hfill&&\hfill$A_{9,1}$\hfill&&\hfill$\Lambda_{5j}^e\mapsto J_a^j\lan\bo,(202)\ran_2,\, \ \La_{5j+(-1)^i}^e\mapsto J_a^j(\lan C^i(200)\ran_2\!\oplus\!(212)),\, \ \La_{5j+2(-1)^i}^e\mapsto J_a^j(\lan C^i(210)\ran_2\!\oplus\!(303)),\ \ i,j=0,1$\hfill&\cr
\tablerule&&&&\hfill8\hfill&&\hfill$D_{10,1}$\hfill&&\hfill$\bo^e\mapsto 
\lan\bo,(121)\ran_4,\ \La^e_1\mapsto \lan(020),(303)\ran_4,\ \La^e_9,\La^e_{10}\mapsto\lan(113)\ran_4$\hfill&
\cr\tablerule\tablerule&&$A_4$&&\hfill3\hfill&&\hfill$A_{9,1}$\hfill&&\hfill$\La_{2j}^e\mapsto J_a^{2j}(\bo\!\oplus\!(0110)),\ \La_{1+2j}^e\mapsto J_a^{2j}((0010)\!\oplus\!(0201)),\ \ 0\le j\le 4$\hfill&
\cr\tablerule&&&&\hfill5\hfill&&\hfill$D_{12,1}$\hfill&&\hfill$\bo^e\mapsto\lan\bo,(0220)\ran_5,\ \La_1^e\mapsto(1111)\!\oplus\!\lan(1001)\ran_5,\ \La_{11}^e,\La_{12}^e\mapsto2\cdot(1111)$\hfill&
\cr\tablerule&&&&\hfill7\hfill&&\hfill$A_{14,1}$\hfill&&\hfill$\,\,\La_{3j}^e\mapsto J_a^{3j}(\bo\!\oplus\!(0330)\!\oplus\!(2002)\!\oplus\!(2112)\!\oplus\!\lan(0403)\ran_c),\ $\hfill&\cr&&&&&&&&\hfill$
\La^e_{3j+(-1)^i}\mapsto J_a^{3j}C^i((2000)\!\oplus\!(0312)\!\oplus\!(1240)\!\oplus\!(2102)\!\oplus\!(3022))\,,\ \ 0\le i\le1,\ 0\le j\le 4$\hfill&
\cr\noalign{\smallskip}}}}}
\noindent{\textbf{Table 2.1.}} The exceptional quantum subgroups for $A_1,A_2,A_3,A_4$. $X_{s,1}$ here means  a Lie-type embedding into $\cV(X_s,1)$.  Weights are not shifted by $\rho$.\medskip

For example we read that 
$A_1$ at level $k=10$ has a Lie-type conformal embedding $\cV(A_1,10)\subset\cV(C_2,1)$, with
branching rules (restrictions) $\mathbf{0}^e\mapsto\bo\oplus(6)$, $\La_1^e\mapsto(3)\oplus(7)$, $\La_2^e\mapsto(4)\oplus(10)$. As always, $\bo^e$ denotes the algebra object $\cA$. The corresponding modular invariant $\cZ=B^tB$ can be read off from the branching rules: the partition function $\sum_{\la,\mu}\cZ_{\la,\mu}\chi_\la\overline{\chi_\mu}$  is $$|\chi^e_{\bo^e}|^2+|\chi^e_{\La_1^e}|^2+|\chi^e_{\La_2^e}|^2=|\chi_{0}+\chi_6|^2+|\chi_3+\chi_7|^2+|\chi_4+\chi_{10}|^2\,.$$

Level-rank duality (see e.g.\ Theorem 5.1 in \cite{OsS})  relates some quantum modules and subgroups of $\cC(\mathfrak{sl}_n,k)$ and $\cC(\mathfrak{sl}_k,n)$. For example, duality relates the exceptional quantum subgroups of $\cC(A_2,5)$ and $\cC(A_4,3)$. We use it  in Sections 4.4 and 5.4 to find all exceptional quantum subgroups of $\cC(A_r,k)$ when $k\le 5$. The results are collected in Table 2.2  (to avoid redundancy with Table 2.1 we restrict to $k\le 5\le r$). The MFC (Rep$_{\cC(\fg,k)}\,\cA)^{\mathrm{loc}}$ is described in Section 4.4.

 \medskip
\centerline{\scriptsize{\vbox{\tabskip=0pt\offinterlineskip
  \def\tablerule{\noalign{\hrule}}
  \halign to 5.7in{
    \strut# &
    \tabskip=0em plus1em &    
    \hfil#&
    #&\hfil#&\vrule\vrule#&    
\hfil#&\vrule#
\tabskip=0pt\cr
&&$\fg$&&$\,k$&&\hfill branching rules\hfill&\cr
\tablerule\tablerule&&$A_5\,$&&\hfill$\,4$\hfill&&\hfill$J_e^j\bo^e\mapsto J_a^j\lan \bo,(01010)\ran_2,J_e^jC_e^i\Lambda^e_1\mapsto J_a^j(\lan C^i(30100)\ran_2\oplus(11011)),$\hfill&
\cr &&&&&&\hfill$J_e^jC_e^i\Lambda^e_2\mapsto J_a^j(\lan C^i(01002)\ran_2\oplus(00200)),\ 0\le i\le 1,0\le j\le2$\hfill&\cr
\tablerule\tablerule&&$A_6\,$&&\hfill$\,5$\hfill&&\hfill$J_e^j\bo^e\mapsto J_a^j( \bo\oplus(002200)\oplus(010010)\oplus(101101)\oplus\lan(003002)\ran_c),$\hfill&\cr 
&&&&&&\hfill$J_e^jC_e^i\La^e_1\mapsto J_a^jC^i((401000)\oplus(100210)\oplus(010130)\oplus(210101)\oplus(020020)),\ 0\le j\le 6$\hfill&
\cr
\tablerule\tablerule&&$A_7\,$&&\hfill$\,4$\hfill&&\hfill$J_e^j\bo^e\mapsto J_a^j\lan\bo,(1100011)\ran_4,J_e^j\La^e_1\mapsto J_a^j\lan (2000002),(0010100)\ran_4,J_e^j\La^e_3,J_e^j\La_4^e\mapsto J_a^j\lan (1100100)\ran_4,\ j=0,1$\hfill&
\cr\tablerule\tablerule&&$A_8\,$&&\hfill$\,3$\hfill&&\hfill$J_e^j\bo^e\mapsto J^j_a\lan \bo,\La_4+\La_5\ran_3, J_e^j\La_1^e,J_e^j\La^e_5\mapsto J^j_a\lan  \La_2+\La_7\ran_3,\ 0\le j\le 2$\hfill&
\cr\tablerule\tablerule&&$A_9\,$&&\hfill$\,2$\hfill&&\hfill$J_e^j\bo^e\mapsto J_a^j\bo\oplus J_a^j(\La_3+\La_7),  J_e^j\La_1^e\mapsto J_a^j\La_3\oplus J_a^j(\La_5+\La_8),\ 0\le i\le 9,\ 0\le j\le4$\hfill&
\cr\tablerule\tablerule&&$A_{20}\,$&&\hfill$\,3$\hfill&&\hfill$J_e^j\bo^e\mapsto J^j_a\lan \bo,\La_4+\La_{17},\La_6+\La_{15},\La_{10}+\La_{11}\ran_3,$\hfill&
\cr&&&&\hfill\hfill&&\hfill$J_e^j\La_6^e\mapsto J^j_a\lan \La_1+\La_8+\La_{12},\La_9+\La_{13}+\La_{20},2\La_2+\La_{17},\La_4+2\La_{19}\ran_3,\ 0\le j\le 6$\hfill&
\cr\tablerule\tablerule&&$A_{27}\,$&&\hfill$\,2$\hfill&&\hfill$J_e^j\bo^e\mapsto J^j_a\lan \bo,\La_5+\La_{23}\ran_2,J_e^j\La_2^e\mapsto J^j_a\lan \La_3+\La_{25},\La_6+\La_{22}\ran_2,\ 0\le j\le 13$\hfill&
\cr\noalign{\smallskip}}}}}
\noindent{\textbf{Table 2.2.}} The exceptional quantum subgroups of $\cC(A_r,k)$ when $k\le5\le r$.\medskip

\section{Background}

\subsection{Modular fusion categories}

The standard reference for  tensor categories is \cite{EGNO}.
A \textit{fusion category}  has direct sums $\la\oplus\mu$, tensor (fusion) products $\la\otimes\mu$, a tensor unit $\b1$, and duals $\la^*$. It  
is semisimple with finitely many isomorphism classes $[\la]$ of simple objects, the Hom-spaces are finite-dimensional vector spaces over $\bbC$,  and End$(\b1)=\bbC$. 
A  \textit{modular fusion category (MFC)}, also called a \textit{modular tensor category}, is a fusion category $\mathcal{C}$ with a braiding $c_{\la,\mu}\in\mathrm{Hom}_\cC(\la\otimes\mu,\mu\otimes\la)$ satisfying a certain nondegeneracy condition. Examples of  MFCs are $\cC(\fg,k)$, where $\fg$ is a simple Lie algebra and $k\in\bbZ_{>0}$.


Associated with any MFC is \textit{modular data},  an SL$_2(\bbZ)$-representation
$R$ on the formal $\bbC$-span of the (equivalence classes of) simple objects. The matrix $T=R\left({1\atop 0}{1\atop 1}\right)$ is diagonal and $S=R\left({0\atop 1}{-1\atop 0}\right)$  is unitary and symmetric. More precisely, given a MFC, $R$ is uniquely defined up to tensoring by a one-dimensional representation of PSL$_2(\bbZ)$: the $S$-matrix is uniquely determined up to a global sign (which we can naturally fix by requiring that $S$ has a strictly positive row) and $T$ is then well-defined up to an arbitrary third root of 1. The \textit{ribbon twist} $\theta(\la)=T_{[\la],[\la]}\overline{T_{[\b1],[\b1]}}$ is independent of these choices. Define $C_{[\la],[\mu]}:=\delta_{[\la^*],[\mu]}$. Then $S^2=C$ and 
$S_{[\lambda^*],[\mu]}=\overline{S_{[\la],[\mu]}}$. All examples we consider are \textit{pseudo-unitary}, i.e.\  $S_{[\la],[\b1]}>0$ for all $[\la]$ (see Corollary 3.6), which simplifies things somewhat. In fact $\cC(\fg,k)$ are \textit{unitary} (Theorem 5.5 of \cite{Ten}), though we don't need that.

Fix once and for all a representative $\la\in\mathrm{Irr}(\cC)$ for each  class $[\la]$. Let Fus$(\cC)$ denote the \textit{fusion ring} (Grothendieck ring) of $\cC$.
By semisimplicity, the product decomposes as $\la\otimes\mu\cong\bigoplus_\nu N_{\la,\mu}^\nu\nu$, where the multiplicities $N_{\la,\mu}^\nu\in\bbZ_{\ge0}$ are called \textit{fusion coefficients}.
\textit{Verlinde's formula} says
\begin{equation}\label{verl}N_{\la,\mu}^\nu=\sum_\phi S_{\la,\phi}S_{\mu,\phi}\overline{S_{\nu,\phi}}/S_{\b1,\phi}\,,\end{equation}
where we sum over the representatives $\phi$.
Combining this with the relation $(ST)^3=S^2$ yields \cite{BE}
\begin{equation}\label{verlmod}S_{\lambda,\mu}=\theta(\la)\,\theta(\mu)\sum_{\nu}\overline{\theta(\nu)}\,N_{\lambda,\mu}^\nu\, S_{\b1,\nu}\,.
\end{equation}

Let $\cC$ be any MFC.
A \textit{simple-current} is a simple object $\la$ in $\cC$ such that $\la\otimes \la^*\cong \b1$.  Using the tensor product, the (equivalence classes of) simple-currents in $\cC$ form an abelian group denoted $\cJ(\cC)$. Recall  \eqref{scform}.

By a \textit{fusion subring} $F$ of the fusion ring Fus$(\cC)$ of an MFC $\cC$, we mean the formal $\bbZ$-span of any subset of Irr$(\cC)$ which is closed under tensor product. For example, the $\bbZ$-span $\bbZ\cJ$ of any subgroup $\cJ\le\cJ(\cC)$ of simple-currents is one. By Corollary 4.11.4 of \cite{EGNO}, a fusion subring is also closed under taking duals. In fact, any fusion subring is the fusion ring of a fusion subcategory of $\cC$.

The \textit{centralizer} $Z_\cC(P)$ of a subset $P\subseteq\mathrm{Irr}(\cC)$ is all $\la\in\mathrm{Irr}(\cC)$ satisfying 
\begin{equation}\frac{S_{\la,\mu}}{S_{\b1,\mu}}=\frac{S_{\la,\b1}}{S_{\b1,\b1}}\ \ \forall\mu\in P\,.\end{equation} 
Then $\la\in Z_\cC(\mu)$ iff $\mu\in Z_\cC(\la)$.
An important example is, for any $\cJ\subseteq\cJ(\cC)$, \begin{equation}\label{ZJ}Z_\cC(\cJ)=\{\la\in\mathrm{ Irr}(\cC)\,|\,\varphi_J(\la)=1\ \forall J\in\cJ\}\,.\end{equation}

\medskip\noindent\textbf{Proposition 3.1.} \textit{Let $\cC$ be any pseudo-unitary MFC and  $P$ any subset of $\mathrm{Irr}(\cC)$. Then the $\bbZ$-span $\bbZ Z_\cC(P)$ is a fusion subring of $\mathrm{Fus}(\cC)$.}

\medskip\noindent\textit{Proof.}  For any $\nu\in \mathrm{Irr}(\cC)$, let $N_\nu$ denote the (non-negative) matrix $(N_\nu)_{\nu',\nu''}=N_{\nu,\nu'}^{\nu''}$ corresponding to multiplying in Fus$(\cC)$ by $\nu$. Then Verlinde's formula \eqref{verl} tells us the eigenvalues of $N_\nu$ are $S_{\nu,\nu'}/S_{\b1,\nu'}$ for $\nu'\in\mathrm{Irr}(\cC)$, with eigenvector the $\nu'$ column of $S$. Hence, $S_{\nu,\b1}/S_{\b1,\b1}$ is the Perron-Frobenius eigenvalue, corresponding as it does to a strictly positive eigenvector, so
\begin{equation}\frac{S_{\nu,\b1}}{S_{\b1,\b1}}\ge\left|\frac{S_{\nu,\nu'}}{S_{\b1,\nu'}}\right|\ \ \forall\nu'\in \mathrm{Irr}(\cC)\label{PFineq}\,.\end{equation}

 For any $\la,\la'\in Z_\cC(P)$ and $\mu\in P$, we find
$$\sum_\nu N_{\la,\la'}^\nu\frac{S_{\nu,\mu}}{S_{\b1,\mu}}=\frac{S_{\la,\mu}}{S_{\b1,\mu}}\frac{S_{\la',\mu}}{S_{\b1,\mu}}=\frac{S_{\la,\b1}}{S_{\b1,\b1}}\frac{S_{\la',\b1}}{S_{\b1,\b1}}=\sum_\nu N_{\la,\la'}^\nu\frac{S_{\nu,\b1}}{S_{\b1,\b1}}\,.$$
Applying the triangle inequality and \eqref{PFineq} to this, we find that $\nu\in Z_\cC(P)$ whenever the fusion coefficient  $N_{\la,\la'}^\nu\ne 0$. In other words, the span $\bbZ Z_\cC(P)$ is closed under fusion products. 
\textit{QED to Proposition 3.1}\medskip

All fusion subrings for  $\cC(\fg,k)$ were classified by Sawin \cite{Sa}. We need:

\medskip\noindent\textbf{Theorem 3.2.} \cite{Sa} \textit{Any fusion subring of $\cC(\fg,k)$ is of the form $\bbZ\cJ$ or $\bbZ Z_\cC(\cJ)$ for some $\cJ\le\cJ(\fg,k)$, unless $(\fg,k)=(B_r,2),(D_r,2)$ or $(E_7,2)$.}\medskip

\subsection{Vertex operator algebras}\label{sect:voa}

For the basic theory and terminology of vertex operator algebras (VOAs), see e.g.\ \cite{LeLi}.
All VOAs $\cV$ in this paper are \textit{(completely) rational}, i.e.\ $\cV$ is $C_2$-cofinite, regular, simple, self-dual and of CFT-type. Huang  \cite{H2} proved that the category Mod$(\cV)$ of modules of a  rational $\cV$ is a MFC. The tensor unit $\b1$ is $\cV$ itself.

We are interested in the VOAs  $\cV(\fg,k)$ where $\fg$ is a simple finite-dimensional Lie algebra over $\bbC$ and $k\in\bbZ_{>0}$ (see e.g.\ Sections 6.2,6.6 of  \cite{LeLi}). These are  rational. The simple $\cV(\fg,k)$-modules are labeled by highest-weights $\la\in P_+^k(\fg)$ or equivalently $P_{\!++}^{\,\ka}(\fg)$ (see \eqref{P+}) and are naturally identified with the level $k$ integrable modules of the affine Lie algebra $\fg^{(1)}$.

Unlike the modular data of a MFC, the modular data of a VOA is uniquely determined.
In any  $\cV$-module $M$, the Virasoro operator $L_0$ is diagonalizable with rational eigenvalues. The minimal eigenvalue exists and is called the \textit{conformal weight} $h_M$.  Then the $T$-matrix is $T_{M,N}=\delta_{M,N}e^{2\pi \ri \,(h_M-c/24)}$ where  $c\in\bbQ$ is the \textit{central charge} of $\cV$ ($c$ fixes the third root of 1 for $T$). The {one-point functions} of $\cV$ yield  vector-valued modular forms for the modular data $R$; most importantly, this holds for the graded-dimensions (characters) $\chi_M(\tau)$ of the $\cV$-modules.

Any CFT-type VOA  (e.g.\ $\cV(\fg,k)$ and their  extensions) have the property that their conformal weights $h_M$ for simple $M$ are $\ge 0$, and $h_M=0$ only for $M=\cV$. It is easy to show from $\chi_M(-1/\tau)=\sum_NS_{M,N}\chi_N(\tau)$ that Mod$(\cV)$ for such $\cV$ is pseudo-unitary:
$$0\le \lim_{q\to 1^-}\frac{\chi_M(q)}{\chi_\cV(q)}=\lim_{q\to 0^+}\frac{\sum_PS_{M,P}\chi_P(q)}{\sum_PS_{\cV,P}\chi_P(q)}=\frac{S_{M,\cV}}{S_{\cV,\cV}}\,.$$

When a  rational VOA $\cV$ is a subVOA of a rational $ \cV^e$ and shares the same conformal vector (hence the same central charge),  we call $\cV$ a \textit{conformal embedding} in $\cV^e$, and $\cV^e$ a \textit{conformal extension} of $\cV$. This paper studies the conformal extensions of $\cV=\cV(\fg,k)$.

The
\textit{Lie-type conformal embeddings} occur when the extended VOA $\cV^e$ is also of Lie algebra type, $\cV^e\cong\cV(\fg',k')$. Then  $k'=1$ and $\fg'$ will be simple if $\fg$ is. The complete list  is explicitly given in \cite{BB,SW}. When $\fg=A_1,A_2,A_3,A_4$, they occur at the levels $k$ listed in Table 2.1, together with $(A_1,4),(A_2,3),(A_3,2)$, where they are also simple-current extensions.

The following elementary observation is very helpful.

\medskip\noindent\textbf{Proposition 3.3.} \textit{Suppose $\cV=\cV(\fg,k)$  has no proper conformal extension of Lie type. If $\cV^e$ is any conformal extension of $\cV$, then the homogeneous spaces $(\cV^e)_1$ and $\cV_1=\fg$ coincide.}

\medskip\noindent\textit{Proof.} Recall that $(\cV^e)_1$ has the structure of a reductive Lie algebra (\cite{DM}, Theorem 1), with $\cV_1$ a Lie subalgebra. Clearly, to be conformal,  $(\cV^e)_1$ must in fact be semisimple (otherwise its central charge would be strictly greater). Let $\widetilde{\cV^e}$ be the subVOA of $\cV^e$ generated by $(\cV^e)_1$. Then $\widetilde{\cV^e}$ is a conformal extension of $\cV$ (necessarily of Lie type), since $\cV^e$ is.
\textit{QED to Proposition 3.3}

\subsection{Algebras in categories and conformal extensions}

One can think of a (modular) fusion category   as a categorification of  (commutative) finite-dimensional  rings. In the same sense,  a categorification of their modules is called
a \textit{module category}. The formal definition is given in \cite{Ost1}; see also Chapter 7 of \cite{EGNO}. There are obvious notions of equivalence, of direct sums, and of indecomposability of module categories. 

An algebra in say a MFC $\cC$ consists of an object $\cA\in\cC$, a multiplication $\mu\in\mathrm{Hom}_\cC(\cA\otimes\cA,\cA)$ satisfying the associativity constraint, and a unit $\iota\in\mathrm{Hom}_\cC(\b1,\cA)$. We can speak of commutative algebras, as well as (left) modules of $\cA$ in $\cC$, in the obvious way.  Again, see \cite{Ost1,EGNO} for formal definitions.
The Main Theorem of \cite{Ost1} says that each indecomposable module category of  $\cC$ is equivalent to the category Rep$_\cC\,\cA$ of left modules in $\cC$ of some associative algebra $\cA\in\cC$ with unit. We are interested though in module categories of pure extension type (recall the discussion in the Introduction). The corresponding algebras are \'etale:

\medskip\noindent\textbf{Definition 3.4.} \textit{Let $\cC$ be a MFC. We call $\cA$ an} \'etale algebra \textit{for $\cC$ if $\cA$ is a commutative associative algebra in $\cC$  and the category $\mathrm{Rep}_\cC\,\cA$ is semisimple. We call $\cA$} connected \textit{if} $\mathrm{dim\, Hom}_\cC(\b1,\cA)=1$. \textit{For $\cC=\cC(\fg,k)$, a connected \'etale algebra is also called a} quantum subgroup.\medskip

The quantum subgroups of $\cC(\fg,k)$  are the basic objects we wish to identify. In pseudo-unitary categories like $\cC(\fg,k)$, a connected \'etale algebra $\cA$ always has trivial ribbon twist, i.e.\ $\theta(\cA)=1$ (Lemma 2.2.4 in \cite{Sch}). 

\'Etale algebras go by other names in the literature -- e.g.\ \cite{RFFS} calls them commutative symmetric special Frobenius algebras and in physics they're called condensable.
Any MFC always has at least one connected \'etale algebra, namely  $\cA=\b1$. Any \'etale algebra is canonically a direct sum of connected ones. Any \'etale algebra is necessarily self-dual: $\cA=\cA^*$ (see Remark 3.4 of \cite{DMNO}). An intrinsic characterization of $\cA$ for which Rep$_\cC\,\cA$ is semisimple is that $\cA$ be \textit{separable} (see Proposition 2.7 in \cite{DMNO}). Section 3 of \cite{DMNO}, as well as \cite{KO}, collect many  properties of \'etale  $\cA$. 

When $\cA$ is connected \'etale, Rep$_\cC\,\cA$ is a fusion category (see Section 3.3 of \cite{DMNO}).  In Rep$_\cC\,\cA$, we write $M\otimes_\cA N$ and Hom$_\cA(M,N)$.

\textit{Restriction} $\mathrm{Res\,{:}\,Rep}_\cC\,\cA\to\cC$ is the forgetful functor sending an $\cA$-module to the   underlying object in $\cC$. It has an adjoint Ind$:\cC\to\mathrm{Rep}_\cC\,\cA$, called \textit{(alpha-)induction}, sending $\la\in\cC$ to a canonical $\cA$-module structure on $\cA\otimes \la$ (Section 1 of \cite{KO}). For all $\cA$-modules $M$ and $\la\in\cC$, we have 
\begin{equation}\mathrm{Hom}_{\cA}(\mathrm{Ind}\,\la,M)\cong\mathrm{Hom}_\cC(\la,\mathrm{Res}\,M)\,.\ 
\label{adjoint}\end{equation}
 Moreover, Ind is a tensor functor, Ind($\la^*)=(\mathrm{Ind}\,\la)^*$, and 
\begin{equation}\mathrm{Res(Ind}(\la)\otimes_{\cA}M)=\la\otimes\mathrm{Res}\,M\,.\label{ResIndTimes}\end{equation}

 Let (Rep$_\cC\,\cA)^{\mathrm{loc}}$ be the full subcategory of \textit{local} (dyslectic) $\cA$-modules \cite{KO}. A simple $\cA$-module $M$ is local iff the ribbon twist $\theta(M)$ is scalar (Theorem 3.2 of \cite{KO}), i.e.\ all subobjects $\la\in \mathrm{Res}\,M$ have the same value of $\theta(\la)$. Although Rep$_\cC\,\cA$ is usually not braided, the subcategory (Rep$_\cC\,\cA)^{\mathrm{loc}}$ is a MFC (Theorem 4.5 of \cite{KO}). For example, $\theta_{\mathrm{(Rep}_\cC\,\cA)^{\mathrm{loc}}}(M)=\theta_\cC(M)$.

Conformal embeddings $\cV\subset\cV^e$ of VOAs were defined  last subsection. 
Conformal extensions $\cV\subset \cV^e,\cV\subset\cV^{\prime e}$ are equivalent if there are VOA isomorphisms $\cV\cong\cV',\cV^e\cong \cV^{\prime e}$ commuting with the inclusions $\cV\subset\cV^e,\cV'\subset\cV^{\prime e}$.    For example, the extensions of $\cV=\cV'=\cV(D_8,1)$ given by $\bo\oplus\Lambda_8$ and $\bo\oplus\Lambda_7$ are inequivalent even though $\cV^e\cong\cV^{\prime e}\cong\cV(E_8,1)$ as VOAs.

\'Etale algebras $\cA,\cA'$  are equivalent if there is an invertible $\phi\in\mathrm{Hom}_\cC(\cA,\cA')$ intertwining in the obvious way the maps defining unit, multiplication etc. This requires Res$\,\cA\cong\mathrm{Res}\,\cA'$ as $\cV$-modules, but that is not sufficient (e.g.\ at the end of Section 3.4 we describe a situation where the inequivalent algebras with the same restriction are parametrized by $H^2(G;\bbC^\times)$ for some group $G$). 

A key result is:

\medskip\noindent\textbf{Theorem 3.5.} \cite{KO,HKL,CKM} \textit{Let $\cV$ be a  rational VOA with MFC $\cC=\mathrm{Mod}(\cV)$. Assume the conformal weights $h_\la$ are positive for all simple $\cV$-modules $\la\not\cong\b1=\cV$. Then the following are equivalent:}

\smallskip\noindent\textbf{(i)} \textit{a conformal extension $\cV^e$ of $\cV$;} 

\smallskip\noindent\textbf{(ii)} \textit{a connected \'etale algebra $\cA$  in $\cC$ such that $\cA\cong \cV^e$ as objects in $\cC$ (with multiplication coming from the vertex operators in $\cV^e$).}

\smallskip\noindent\textit{Moreover, $\mathrm{Mod}(\cV^e)$ is equivalent to  $(\mathrm{Rep}_\cC\,\cA)^{\mathrm{loc}}$ as  MFCs. The forgetful functor $(\mathrm{Rep}_\cC\,\cA)^{\mathrm{loc}}\to \cC$ restricts $\cV^e$-modules to $\cV$-modules. Algebras $\cA,\cA'$ are equivalent iff the corresponding $\cV^e,\cV^{e\prime}$ are equivalent as extensions of $\cV$.}

\medskip In the language of VOAs, $\cA$-modules are a solitonic generalization of the notion of VOA module -- see Remark 2.7 in \cite{CKM}.

For example, let $\cW$ be a holomorphic VOA (i.e.\ one with Mod$(\cW)\cong\mathrm{Vec}$), and $G$  a finite group of automorphisms of $\cW$. Then the fixed points $\cV=\cW^G$, called the \textit{orbifold} of $\cW$ by $G$, is also expected  to be  rational. When this is so, Mod$(\cV)$ will be the twisted Drinfeld double $\cZ(\mathrm{Vec}_G^\omega)$ for some 3-cocycle $\omega$ of $G$. Up to equivalence, the simple $\cV$-modules are parametrized by pairs $(g,\chi)$ where $g$ runs over a set of conjugacy class representatives in $G$, and $\chi$ is a simple projective representation of the centralizer $Z_G(g)$, with projective multiplier determined by $\omega$. Then $\cW$ is a conformal extension of $\cV$ with \'etale algebra $\cA\cong\bigoplus_\chi\mathrm{dim}(\chi)\,(e,\chi)$. There is a unique $\cA$-module $M_g$ for each $g\in G$, called the $g$-twisted $\cW$-module. Induction sends $(g,\chi)$ to $\bigoplus_h\mathrm{dim}(\chi)\,M_h$ where the sum is over the conjugacy class of $g$.

A more pertinent example is the extension $\cV(A_1,10)\subset\cV(C_2,1)$. It has three local $\cA$-modules (corresponding to $P_+^1(C_2)$) and three nonlocal $\cA$-modules. The inductions are worked out in Example 2.2 of \cite{BEii}.

 Thanks to restriction, the conformal weights of $\cV^e$-modules are a subset of those of $\cV$-modules, so Mod$(\cV^e)$ is pseudo-unitary if Mod$(\cV)$ is:

\medskip\noindent\textbf{Corollary 3.6.}  \textit{Let $\cV$ satisfy the conditions of Theorem 3.5 and let $\cV^e$ be any conformal extension. Then $\mathrm{Mod}(\cV^e)$ is pseudo-unitary. In particular, for any quantum subgroup $\cA$ of $\cC(\fg,k)$, $(\mathrm{Rep}_{\cC(\fg,k)}\,\cA)^{\mathrm{loc}}$ is pseudo-unitary.}

\subsection{Simple-current extensions}\label{sect:sc}

Recall the discussion of simple-currents in Sections 2.1 and 3.1.

The following has appeared in many different guises in the literature: 

\medskip\noindent\textbf{Proposition 3.7.} \textit{Let $\cA$ be a connected \'etale algebra in a pseudo-unitary MFC $\cC$. Let $B=(B_{M,\mu})$ be the branching rule matrix (where $M$ runs over the simple local $\cA$-modules and $\mu$ the simple objects in $\cC$).}

\smallskip\noindent\textbf{(a)} \textit{When $J\in \cC$ is a simple-current  and $M$ is arbitrary, $B_{M,J}\in\{0,1\}$.}

\smallskip\noindent\textbf{(b)} \textit{Write $\cJ_\cA$ for the set of simple-currents $J$ with $B_{\b1^e,J}=1$. 
Then $\cJ_\cA$ is a subgroup of $\cJ(\cC)$ and  $\theta(J)=1$  $\forall J\in\cJ_\cA$. For any local $M$ and any $\mu\in\cC$, $B_{M,J\mu}=B_{M,\mu}$ $\forall J \in \cJ_\cA$. When $B_{M,\mu}\ne 0$ for $M$  local, we have $\mu\in Z_\cC(\cJ_\cA)$.}

\medskip\noindent\textit{Proof.} Choose any simple-current $J$ with $B_{\b1^e,J}>0$. Then \eqref{adjoint} tells us that Ind$\,J$ contains $B_{\b1^e,J}$ copies of $\b1^e$. But $J\otimes J^*=\b1$ and Ind is a tensor functor, so  Ind$\,J=\b1^e$ and $B_{\b1^e,J}=1$.
Now apply \eqref{ResIndTimes} to $X=J$ to see that Res$\,M=J\otimes \mathrm{Res}\,M$, i.e.\ $B_{M,J\mu}=B_{M,\mu}$ for all $M,\mu$. In particular,  $\cJ_\cA$ is a group. Moreover, locality of $\b1^e$ resp.\ $M$  forces $\theta(J)=1$ and $\varphi_J(\mu)=1$ whenever  $\mu\in\mathrm{Res}\,M$, using \eqref{scform}.  \textit{QED to Proposition 3.7}\medskip

Most quantum subgroups $\cA$ for $\cC(\fg,k)$ are direct sums of simple-currents:  $\mathrm{Res}\,\cA=\bigoplus_{J\in \cJ_\cA}J$. We call such $\cA$    \textit{simple-current \'etale algebras}  (\cite{FRS} calls them \textit{commutative Schellekens algebras}). When $\cV^e\supset\cV$ has a simple-current \'etale algebra $\cA$, we say $\cV^e$ is a \textit{simple-current extension} of $\cV$.

\medskip\noindent\textbf{Proposition 3.8.} \textit{Let $\cC$ be a  pseudo-unitary MFC.}

\smallskip\noindent\textbf{(a)} \textit{Let $\cA$ be a  simple-current  \'etale algebra for $\cC$. If  $M$ is a simple $\cA$-module, then $\mathrm{Res}\,M=n_M\langle \la\rangle_{\cJ_\cA}$ for some $\la$ (all multiplicities equal $n_M$). Moreover, $M$ is local iff $\la\in Z_\cC(\cJ_\cA)$. If $\la\in\cC$ is simple, then $\mathrm{Ind}\,\la$ is a direct sum of $|\cJ_\cA|/(n^2\|\lan \la\ran_{\cJ_\cA}\|)$ distinct simple $\cA$-modules $M$, each with the same multiplicity $n_M=n$.}

\smallskip\noindent\textbf{(b)} \textit{Conversely, let $\cJ$ be a subgroup of simple-currents $\cJ(\cC)$, with ribbon twist $\theta(J)=1$ for all $J\in \cJ$. Then (up to equivalence) there is exactly one connected \'etale algebra structure on $\cA=\bigoplus_{J\in\cJ}J$.}

\medskip\noindent\textit{Proof.} For part (a), see Section 4.2 of \cite{FRS}.
For (b),  existence and uniqueness is established in Corollary 3.30 of \cite{FRS}. \textit{QED to Proposition 3.8}\medskip

More generally, \cite{FRS} studies the (not necessarily commutative) algebras of simple-current type. They find (see their Remark 3.19(i)) that the different algebraic structures on $\cA=\oplus_{J\in \cJ}J$ are parametrised by the Schur multiplier $H^2(\cJ;\bbC^\times)$ (so  it's unique  iff $\cJ$ is cyclic). But for any $\cJ$, requiring commutativity of $\cA=\oplus_{J\in \cJ}J$ forces uniqueness of that algebraic structure.

When $\cC=\mathrm{Mod}(\cV)$ (which always holds for us), the existence and uniqueness in Proposition 3.8(b) can be done by VOAs (Proposition 5.3 of \cite{DM}). 

The modular data $S,T$ for the simple-current extensions of $\cV(\fg,k)$ are given in Section 6 of \cite{FSS}. An example with $n_M>1$ is $\cC(D_4,2)$ with $\cA$ being the direct sum of all 4 simple currents:  the extended VOA is $\cV(E_7,1)$ with branching rules Res$\,\bo^e=\cA$ and Res$\,\Lambda_1^e=2\cdot\Lambda_2$ (so $n_{\Lambda^e_1}=2$).

\subsection{The Galois symmetry}\label{Sect:GaloisLie}

Recall that $\xi_n=\exp(2\pi \ri/n)$. In this subsection we restrict to the MFCs $\cC(\fg,k)$. Recall that the coroot lattice $Q^\vee$ is even (i.e.\ all $(\alpha|\alpha)\in2\bbZ$), and its dual lattice $Q^{\vee*}$ is the weight lattice of $\fg$. As in Section 2.1, ${W}$ is the (finite) Weyl group and $\rho$ is the Weyl vector. We see from  \eqref{Sgenform} that  the $S$-entries $S_{\la,\mu}$ for $\cC(\fg,k)$ manifestly lie in the cyclotomic field $\bbK=\bbQ[\xi_{N_\fg\kappa},x]$, where $x^2\in \bbQ$ and  $N_\fg\in\bbZ_{>0}$ is such that $N_\fg\,Q^{\vee*}\subseteq Q^\vee$. For example, $N_{A_r}=r+1$.

Identify the Galois group Gal($\bbQ[\xi_{N_\fg\ka}]/\bbQ)$ with $(\bbZ/N_\fg\ka)^\times$, by $\sigma(\xi_{N_\fg\ka})=\xi_{N_\fg\ka}^\ell$.  Galois theory says any such $\sigma$ lifts to an   automorphism of  $\bbK$.

Crucial for us is a Galois symmetry  \cite{dBG,CG} holding in any MFC because of Verlinde's formula \eqref{verl}. We need the following geometric interpretation valid for $\cC(\fg,k)$. 

\medskip\noindent\textbf{Lemma 3.9.} \cite{comm} \textit{Fix any MFC $\cC(\fg,k)$, and write $\bbK=\bbQ[\xi_{N_\fg\ka},x]$ as before.  Choose any $\sigma\in \mathrm{Gal}(\bbK/\bbQ)$ and corresponding integer $\ell$ coprime to $N_\fg\ka$.}

\begin{itemize}

\item[(a)] \textit{Given any $\la\in P^{\,\ka}_{\!++}(\fg)$, there is a unique $\alpha\in Q^\vee$ and $\omega\in {W}$ such that $\omega(\ell\la)+\kappa\alpha\in P_{\!++}^{\,\ka}(\fg)$. Denote by $\la^\sigma=\ell.\la$ that element of $P_{\!++}^{\,\ka}(\fg)$.}

\item[(b)] \textit{There are signs $\epsilon_\sigma:P_{\!++}^{\,\ka}(\fg)\to\{\pm1\}$ such that, for any $\la,\mu\in P_{\!++}^{\,\ka}(\fg)$, }
\begin{equation}\sigma(S_{\la,\mu})=\epsilon_\sigma(\la)\,S_{\ell.\la,\mu}=\epsilon_\sigma(\mu)\,S_{\la,\ell.\mu}\,,\label{Galrule}\end{equation}
\textit{where $\ell.\la$ and $\ell.\mu$ are as in} (a).

\end{itemize}

\noindent\textit{Proof.} By Proposition 6.6 of \cite{Kac}, given any weight vector $\beta\in Q^{\vee*}$, there is an element  $\omega\in {W}$ of the finite Weyl group and a coroot vector $\alpha\in Q^\vee$ such that exactly one of the following holds:

(i) $\beta=\omega(\beta)+\kappa\alpha$ and det$\,\omega=-1$, or 

(ii) $\omega(\beta)+\kappa\alpha\in P^{\,\ka}_{\!++}(\fg)$, in which case $\omega$ and $\alpha$ are unique. 

\noindent For any $\beta,\gamma\in Q^{\vee*}$, define 
$$S(\beta,\gamma)=\sum_{w\in{W}}\mathrm{det}(w)\,\exp(-2\pi \ri (\beta|w(\gamma))/\kappa)\,.$$
Then for $\la,\mu\in P_{\!++}^{\,\ka}(\fg)$, \eqref{Sgenform} says $S_{\la,\mu}=xS(\la,\mu)$. Note that $S(w\beta+\ka \alpha,\gamma)=\mathrm{det}(w)\,S(\beta,\gamma)$ for any $w\in W$ and $\alpha\in Q^\vee$, so $S(\beta,\gamma)=0$ $\forall\gamma\in Q^{\vee*}$ when $\beta$ is type (i).

 Note that $\sigma(S(\beta,\gamma))=S(\ell \beta,\gamma)$. Let $\la\in P_{\!++}^{\,\ka}(\fg)$.
  Then for  $\ell$ coprime to $N_\fg\ka$, $\beta=\ell\la$ is type (ii) (otherwise the $\la$-row of the invertible $S$-matrix would be identically 0). Together these imply both statements (a) and (b), where $\epsilon_\sigma(\la)=\epsilon'_\sigma\,\mathrm{det}(\omega)$ for $\epsilon'_\sigma=\sigma(x)/ x\in\{\pm1\}$. \textit{QED to Lemma 3.9}\medskip

Equation \eqref{Galrule} holds in any MFC, though the geometric interpretation of Lemma 3.9(a) is lost.
To see the value of \eqref{Galrule}, let $B$ be any $\bbQ$-matrix satisfying \eqref{modinv}. Choose any automorphism $\sigma$ of the field generated over $\bbQ$ by the entries of $S,S^e$.
Applying $\sigma$ to the $(M,\la)$-entry of $\overline{S^e}BS=B$  gives 
\begin{equation}\epsilon^e_\sigma(M)\,\epsilon_\sigma(\la)\sum_{N,\mu}\overline{S^e_{M^\sigma,N}}B_{N,\mu}S_{\mu,\la^\sigma}=\epsilon^e_\sigma(M)\,\epsilon_\sigma(\la)B_{M^\sigma,\la^\sigma}\,,\nonumber\end{equation}
since the field is cyclotomic, and hence
\begin{equation}B_{M,\la}=\epsilon_\sigma(\la)\,\epsilon_\sigma^e(M)\,B_{M^\sigma,\la^\sigma}\,.\label{GalBranch}\end{equation}
To our knowledge, this equation appears here for the first time. It implies that whenever a branching rule $B_{M,\la}\ne 0$, then $\epsilon_\sigma(\la)=\epsilon^e_\sigma(M)$ for all such $\sigma$ (since branching rules are non-negative). It also implies the weaker but older  $\cZ_{\la,\mu}=\epsilon_\sigma(\la)\,\epsilon_\sigma(\mu)\cZ_{\la^\sigma,\mu^\sigma}$, true for any modular invariant $\cZ$.

In particular, let $\cA=\bigoplus_\la Z_\la\la$ be \'etale. Then for any $\la\in P_{\!++}^{\,\ka}(\fg)$ with $Z_{\la}>0$, we have $\epsilon_\sigma(\b1)=\epsilon_\sigma(\la)$ $\forall\sigma$. Equivalently,  $\sigma(S_{\la,\b1}/S_{\b1,\b1})>0$ $\forall\sigma$, i.e.\   qdim$(\la)$ is totally positive. 

\medskip\noindent\textbf{Definition 3.10.} \textit{If $\cV(\fg,k)$ has no conformal extension of Lie type, call $\la\in P_{\!++}^{\,\ka}(\fg)$ a} candidate \textit{if $\theta(\la) =1$,  $\mathrm{qdim}(\la)$ is totally positive, and $h_\la\ne1$. If $\cV(\fg,k)$ has an extension of Lie type, use the same definition except drop the condition  $h_\la\ne1$. In either case, let $\fC_k(\fg)$ be the set of all {candidates}.}

\medskip  Using Proposition 3.3, we have proved:

\medskip\noindent\textbf{Corollary 3.11.} \textit{Let $\cA$ be any quantum subgroup of $\cC(\fg,k)$. Then any $\la\in\cA$ lies in $\fC_k(\fg)$.}\medskip

The set of level $k$ weights  $P_{\!++}^{\,\ka}(\fg)$  is almost always much larger than the subset $\fC_k(\fg)$ of candidates. 
For example, for $\fg=A_1$, we know \cite{su2}
\begin{equation}
\fC_k(A_1)=\left\{\begin{matrix}\{\b1\}&k\equiv 1,2,3\ (\mathrm{mod}\ 4),\ k\ne 10\cr \lan \b1\ran_2&k\equiv 0\ (\mathrm{mod}\ 4),\ k\ne 28\cr \{\b1,7\}&k=10\cr\lan\b1,11\ran_2&k=28\end{matrix}\right.\nonumber\end{equation}
Here, $\lan\cdot\ran_2$ means the orbit by the order-2 group of simple-currents. By comparison, the number of simple objects in $\cC(A_1,k)$ is $k+1$. $\fC_k(A_2)$ is also known for all $k$ \cite{su3}. We expect that max$_k\|\fC_k(\fg)\|$ is finite for any given $\fg$ (it is 12 for $\fg=A_2$), though $\|P_{+}^{k}(\fg)\|$ grows like \eqref{growthP+}. Every $\la\in\fC_k(A_1)$ appears in a quantum subgroup; for higher rank $\fg$ this often fails.

The signs $\epsilon_\sigma(\la)$ can be computed efficiently as follows. Since the sign of $\sigma(S_{\la,\b1}/S_{\b1,\b1})$ equals $\epsilon_\sigma(\la)\epsilon_\sigma(\b1)$, the Weyl denominator identity for $\fg$ gives us
\begin{equation}\epsilon_\sigma(\la)\,\epsilon_\sigma(\b1)=\prod_{\alpha>0}\mathrm{sign}\left(\sin(\pi\ell(\la|\alpha)/\ka)\right)\,\mathrm{sign}\left(\sin(\pi\ell(\rho|\alpha)/\ka)\right)\label{parityform}\end{equation}
where the product is over the positive roots of $\fg$.

\section{The underlying theorems}

\subsection{The algebra constraint}

Given any MFC $\cC$, let $\cC_0$ denote the \textit{adjoint subcategory} $Z_\cC(\cJ(\cC))$ (recall \eqref{ZJ}).  It is a braided fusion subcategory of $\cC$. For example, the simples of $\cC(A_r,k)_0$ are those $\la\in P_+^k(A_r)$ with $r'|t(\la) $. 

Fix any simple $\fg$ and any $k\in\bbZ_{>0}$. For any $\mu\in P_{\!++}^{\,\ka}(\fg)$,
 the condition  $\theta(\mu)=1$ is equivalent to the condition that $h_\mu:=\frac{1}{2\ka}(\mu-\rho|\mu+\rho)\in\bbZ_{\ge 0}$ ($h_\mu$ equals 0 only for $\mu=\b1$). Given a  quantum subgroup $\cA$ of $\cC(\fg,k)$, let $h(\cA)$ be the minimum of $h_\mu$ as we run over all $\mu\in\cA\cap\cC(\fg,k)_0$, $\mu\ne\b1$ (if no such $\mu$ exist, put $h(\cA)=\infty$).
Likewise, let $h_\fC(\fg,k)$ be the minimum of $h_\mu$ as we run over all candidates $\mu\in\fC_k(\fg)\cap\cC(\fg,k)_0$, $\mu\ne\b1$.  Clearly, $h_\fC(\fg,k)\le h(\cA)$.

The following key result was extracted from  the proof of Lemma 3.1.1 of \cite{Sch}. It was surely known to Ocneanu.

\medskip\noindent\textbf{Ocneanu's Lemma 4.1.} \textit{Let $\cA$ be a quantum subgroup of $\cC(\fg,k)$.  If \begin{equation}(\lambda+\lambda^*-2\rho|\lambda+\lambda^*)<2\kappa\, h(\cA)\label{OcneanusBound}\end{equation}
for some  $\la\in P^{\,\ka}_{\!++}(\fg)$, then $\mathrm{Ind}(\la)\in\mathrm{Rep}_{\cC(\fg,k)}\,\cA$ is simple.}\medskip

\noindent\textit{Proof.} Write $\cC=\cC(\fg,k)$. By Lemma 2.4 of \cite{Ost1},  
$$\mathrm{dim\,End}_{\mathrm{Rep}_\cC\cA}(\mathrm{Ind}\,\lambda)=\mathrm{dim\,Hom}_\cC(\lambda,\cA\otimes\lambda)=\mathrm{dim\,Hom}_\cC(\lambda\otimes\lambda^*,\cA)\,.$$
Hence if  Ind$\,\la$ is not simple, then there would exist some $\mu\in \cA$, $\mu\ne \b1$, such that $\mu$ appears in the decomposition of the fusion product $\lambda\otimes\lambda^*$. Of course, $\la\otimes\la^*\in\cC_0$. In particular, $\mu= \lambda+\lambda^*-\rho-\alpha$ for some root vector $\alpha\in\sum_i\bbZ_{\ge0}\alpha_i$. Since $\theta(\mu)=1$ and $\mu\ne\b1$, we have $(\mu|\mu)\ge 2\kappa+(\rho|\rho)$. But for $\nu\in P^{\,\kappa}_{\!++}(\fg)$ and any simple root $\alpha_j$ of $\fg$, $$(\nu|\nu)-(\nu-\alpha_j|\nu-\alpha_j)>2(\nu-\alpha_j|\alpha_j)\ge2(\rho-\alpha_j|\alpha_j)=2a_j^\vee/a_j-2a_j^\vee/a_j=0$$
The first inequality is because $(\alpha_j|\alpha_j)>0$; the second is because $\nu-\rho$ is a sum of fundamental weights.
So the norm of $\mu$ is bounded above by that of $\lambda+\lambda^*-\rho$. \textit{QED to Ocneanu's Lemma}\medskip

A key to this paper is blending Ocneanu's Lemma with the Galois symmetry of Section 3.5. Recall $f_\fg$ from Section 2.1.

\medskip\noindent\textbf{Lemma 4.2.} \textit{Let $\cA$ be a quantum subgroup of $\cC(\fg,k)$. Choose any $\ell\in\bbZ_{>1}$  coprime to $f_\fg\ka$. Suppose $\ell$ satisfies the bound $\ell<\frac{1}{2}+\sqrt{\frac{6L\kappa}{h^\vee \mathrm{dim}(\fg)}+\frac{1}{4}}$, where $L=\mathrm{min}\{h(\cA),h^\vee\mathrm{dim}\,\fg/(6(h^\vee-1))\}$. Then $\ell\rho\in P_{\!++}^{\, \kappa}(\fg)$, $\ell\rho$ is not a simple-current, and $\mathrm{Ind}(\ell\rho)$ is both simple and local.}\medskip

\noindent\textit{Proof.} 
Recall $N_\fg$ from Section 3.5. Galois automorphisms require $\ell$ coprime to $N_\fg\ka$, not merely $f_\fg\ka$. However,
 choose any $\ell'$ coprime to $N_\fg\kappa$ such that $\ell'\equiv \ell$ (mod $f_\fg\kappa$). Let $\ell'.\b1\in P_{\!++}^{\, \ka}(\fg)$ be as in Lemma 3.9(a). To show $\ell'.\b1=\ell\rho$, first note that $\ell'\rho\equiv \ell\rho$ (mod $\ka Q^\vee$), since by definition $f_\fg\rho$ is a sum of coroots. It suffices then to show that  $\ell\rho=\sum_{i=1}^r\ell\Lambda_i$ lies in $P_{\!++}^{\, \ka}(\fg)$. But 
 $$\sum_{i=1}^r\ell a_i^\vee=\ell\,(h^\vee-1)<\frac{h^\vee-1}{\ell-1}\frac{6\ka}{h^\vee\mathrm{dim}\,\fg}\frac{h^\vee\,\mathrm{dim}\,\fg}{6\,(h^\vee-1)} =\frac{\ka}{\ell-1}\le \ka\,.$$ Thus $\ell'.\b1=\ell\rho$.

None of the simple-currents of $\cC(\fg,k)$ are of the form $\ell\rho$ except for $\fg=A_1$ with $\ell=\ka-1$. However we learned in the previous paragraph that $\ell<\ka/(\ell-1)\le\ka/2<\ka-1$, since $f_{A_1}=2$.
 
Let $\cZ=B^tB$; then the Galois symmetry says $\cZ_{\ell\rho,\ell\rho}=\cZ_{\rho,\rho}=1$. This means there is exactly one local $\cA$-module $M$ such that $\ell\rho\in\mathrm{Res}\,M$. By Frobenius Reciprocity, $M\in\mathrm{Ind}\,(\ell\rho)$. However, Ind$(\ell\rho)$ is simple, by Ocneanu's Lemma: $(\ell\rho+(\ell\rho)^*-2\rho|\ell\rho+(\ell\rho)^*)=4(\ell^2-\ell)(\rho|\rho)$ satisfies \eqref{OcneanusBound}, using \eqref{strange}. Thus Ind$(\ell\rho)=M$.
\textit{QED to Lemma 4.2}\medskip

$L$ in Lemma 4.2 obeys $L\ge1$. We discuss local $\cA$-modules in Section 3.3. 

The next results  give implications of the existence of $\la$ for which Ind$\,\la$ is both simple and local. Recall the discussion of centralizer in Section 3.1. We say $\la\in Z_\cC(\cA)$ if $\la\in Z_\cC(\mu)$ for all $\mu\in\cA$.

\medskip\noindent\textbf{Proposition 4.3.} \textit{Let $\cC$ be a pseudo-unitary MFC with  connected \'etale  algebra  $\cA$, and $\mathrm{Ind}\,\la$ be simple  and local for some $\la\in\cC$. Then $\la\in Z_\cC(\cA)$.} 

\medskip\noindent\textit{Proof.} Since
Res(Ind$\,\la)=\cA\otimes \la$ and $\mathrm{Ind}\,\la$ is both simple  and local, 
 any $\nu$ with  fusion coefficient $N_{\la,\mu}^\nu>0$ for some $\mu\in\cA$ must have $\theta(\nu)=\theta(\la)$. Then \eqref{verlmod} becomes 
$$\frac{S_{\lambda,\mu}}{S_{\b1,\b1}}=\theta(\lambda)\,\theta(\mu)\sum_{\nu}\overline{\theta(\nu)}\,N_{\lambda,\mu}^\nu \frac{S_{\nu,\b1}}{S_{\b1,\b1}}=\frac{S_{\la,\b1}}{S_{\b1,\b1}}\frac{S_{\mu,\b1}}{S_{\b1,\b1}}$$
where we also use the fact that $\theta(\mu)=1$, and the fact that the quantum-dimensions $S_{\nu,\b1}/S_{\b1,\b1}$ give a 1-dimensional representation of the fusion ring Fus$(\cC)$. This says $\la\in Z_\cC(\cA)$, as desired. \textit{QED to Proposition 4.3}\medskip

\medskip\noindent\textbf{Proposition 4.4.} \textit{Let $\cC=\cC(\fg,k)$. Exclude here $(\fg,k)= (B_r,2)$ with $2r+1$ square-free, and $(\fg,k)=(D_r,2)$ with either $r$  or $r/2$ square-free. Let $\cA$ be an exceptional quantum subgroup in $\cC$. If  $\mathrm{Ind}\,\la$ is both simple  and local, for some $\la\in\cC$, then $\la$ is a simple-current}.

\medskip\noindent\textit{Proof.} Choose some $\mu\in \cA$ which is not a simple-current. Suppose for contradiction that there is a non-simple-current $\la\in \cC$  for which $\mathrm{Ind}\,\la$ is both simple and local. Then by Proposition 3.1, $\bbZ Z_\cC(\mu)$ is a fusion subring which contains non-simple-currents (e.g.\ $\la$ thanks to Proposition 4.3), as does its centralizer (which contains $\mu$). By Theorem 3.2, we must have $(\fg,k)=(B_r,2),(D_r,2)$ or $(E_7,2)$ for some $r$.

Consider first $(\fg,k)=(B_r,2)$, where $\ka=2r+1$ is square-free. Then the candidates $\fC_2(B_r)\subset P_+^2(B_r)$ cannot contain $\La_a$ or $2\La_r$ (since  $(\La_a|\La_a+2\rho)\equiv -a^2$ and $(2\La_r|2\La_2+2\rho)\equiv -r^2$ (mod $\ka$)), and $\fC_2(B_r)$ contains neither $\La_r$ nor $\La_1+\La_r$ (since both have quantum-dimension $\sqrt{\ka}$, which isn't totally positive). This excludes from $\fC_2(B_r)$ all non-simple-currents.

Consider next $(\fg,k)=(D_r,2)$, where $\kappa=2r$ is not a multiple of 16, nor has any divisor of the form $n^2$ for odd $n>1$.  Then $\fC_2(D_r)$ also contains only simple-currents: for $a<r-1$, $(\La_a|\La_a+2\rho)\equiv -a^2$ and  $(\La_{r-1}+\La_r|\La_{r-1}+\La_r+2\rho)\equiv -(r-1)^2$ (mod $\ka$),  
and the four spinors have quantum-dimension $\sqrt{r}$, which is not totally positive (unless $r=4$, when the spinors have the wrong ribbon twist $\theta$). 

Likewise,   $\fC_2(E_7)=\{\b1\}$. \textit{QED to Proposition 4.4}\medskip

For rank $r\le8$, the only $B_r$ and $D_r$ at level 2 excluded from Proposition 4.4 are $B_4$ and $D_8$. The modular invariants for all $B_r$ and $D_r$ at level 2 are classified in \cite{orth23}; many aren't realized by quantum subgroups and modules.

 \medskip\noindent\textbf{Corollary 4.5.} \textit{Let $\cC(\fg,k)$ and $\cA$ be as in Proposition 4.4. Suppose  some  non-simple-current $\mu\in P_{\!++}^{\,\ka}(\fg)$ has $(\mu+\mu^*-2\rho|\mu+\mu^*)<2\ka \,h(\cA)$. Then $\mu$ has zero multiplicity in any local  $\cA$-module, i.e.\ the $\mu$-column of the corresponding branching matrix $B$ is identically 0.}\medskip

Otherwise Ind$\,\mu$ would be both simple and local. We conclude:

 \medskip\noindent\textbf{Theorem 4.6.}   \textit{Let $\cC(\fg,k)$ be as in Proposition 4.4, and write $L=\mathrm{min}\{h_\fC(\fg,k),$ $h^\vee\mathrm{dim}\,\fg/(6(h^\vee-1))\}$. Let $p$ be the smallest prime coprime to $f_\fg\ka$. If}
\begin{equation}\label{keybound}
p\,(p-1)<\frac{6\kappa\, L}{h^\vee\,\mathrm{dim}(\fg)}\end{equation}
\textit{then $\cC(\fg,k)$ has no exceptional quantum subgroups.}

\medskip The following observation will prove very convenient:

\medskip\noindent\textbf{Theorem 4.7}. \textit{If $\cA$ is an exceptional quantum subgroup of $\cC(\fg,k)$, then $\cA\cap\cC(\fg,k)_0$ is also an exceptional quantum subgroup of $\cC(\fg,k)$.}

\medskip\noindent\textit{Proof.}  Write $\cC=\cC(\fg,k)$. Note that, by Corollary 3.3 of \cite{DNO}, $\cA\cap\cC_0$ is also  a quantum subgroup of $\cC$.

Assume first that  $S_{\b1,\b1}^2\not\in\bbQ$ for $\cC$. Then  (recall Lemma 3.9(b)) there is an integer $\ell$ coprime to $N_\fg\ka$ such that $\la:=\ell.\b1$ is not a simple-current (choose  $\ell$ so that $\sigma_\ell(S^2_{\b1,\b1})\ne S^2_{\b1,\b1}$).   Suppose for contradiction that $\cA\cap\cC_0$ is a simple-current \'etale algebra. Consider $M=\mathrm{ Ind}\,\la$ for $\cA$. As in the proof of Lemma 4.1, $M$ must be simple in Rep$_\cC\,\cA$: if dim$\,\mathrm{End}_{\mathrm{ Rep}_\cC\cA}(M)>1$ then $\la\in\mu\otimes\la$ for some $\mu\in\cA\cap\cC_0$ with $\mu\ne\b1$, but then $J\otimes\la=\la$ for some simple-current $J\ne\b1$ which is absurd ($J\otimes (\ell.\b1)=\ell.(J^{\ell^{-1}}\b1)$). By \eqref{GalBranch}, $M$ is also local. But then Proposition 4.4 contradicts the assumption that $\cA$ is exceptional. Thus $\cA\cap\cC_0$ is also exceptional, and we're done.

It remains to consider $\cC$ with $S_{\b1,\b1}^2\in\bbQ$. By Figure 2 of \cite{GS}, this happens only for the $A,B,D,E$ algebras at level 1, the $B$ and $D$ algebras at level 2, and  $(A_1,4),(A_2,3),(A_3,2),(C_4,1)$. The modular invariants for $\cC(\fg,1)$ are given in Theorem 5 of \cite{comm}; none of those relevant here have exceptional quantum subgroups. Table 2.1 shows none of $(A_1,4),(A_2,3),(A_3,2)$ possess exceptional quantum subgroups. The modular invariants of $B$ and $D$ at level 2 are classified in \cite{orth23}, and we see that only $B_r$ when $2r+1$ is a perfect square can have $\cA\not\in\cC(B_r,2)_0$ (namely $\cA=\bo\oplus\bigoplus_{i\ge1}\La_{i\sqrt{2r+1}}\oplus (\La_1+\La_r)$), but $\cA\cap \cC(B_r,2)_0=\bo\oplus\bigoplus_{i\ge1}\La_{i\sqrt{2r+1}}$ is also exceptional.  
 \textit{QED to Theorem 4.7}\medskip

$\cA\cap\cC_0$ in Theorem 4.7  means to retain only those subobjects of $\cA$  (together with their multiplicities) lying in $\cC_0$. Recall $\cJ_\cA$ from Proposition 3.6(b).

 \medskip\noindent\textbf{Corollary 4.8.} \textit{Suppose $\cC=\cC(\fg,k)$ has no exceptional quantum subgroup $\cA$ with $\cJ_\cA$ maximal (i.e.\ with $\cJ_\cA=\{J\in\cJ(\cC)\,|\,\theta(J)=1\}$). Then all quantum subgroups of $\cC$ are simple-current.}
 
 \medskip 
 \noindent\textit{Proof.} By Theorem 4.7, we know  there is an exceptional quantum subgroup $\cA'\in \cC_0$. Choose any simple-current $J\in\cJ(\cC)$. Then ${J}_e=\mathrm{Ind}(J)\in\mathrm{Rep}_{\cC}\,\cA'$ is simple and invertible. Since $\cA'\in \cC_0$, it is also local: by \eqref{scform}, for any $\la\in\cA'$, $\theta(J\la)=\theta(J)$ is independent of $\la$. Thus $J_e$ is a simple-current in the MFC (Rep$_\cC\,\cA')^{\mathrm{loc}}$. If $\theta(J)=1$, then $\theta(J^e)=1$. By inspection we see that for all $\cC(\fg,k)$ ($\fg$ simple) the set of all $J$ with $\theta(J)=1$ forms a group. Let $\cA''$ denote the direct sum of all such $J^e$ with $\theta(J^e)=1$ (to avoid repetition take at most one representative $J$ from each coset $\cJ(\cC)/\cJ_{\cA'}$). Then $\cA''$ will be \'etale in (Rep$_\cC\,\cA')^{\mathrm{loc}}$, and (see Section  3.6 of \cite{DMNO}) the restriction of $\cA''$ to $\cC$ is a quantum subgroup $\cA$ of $\cC$. By construction, $\cA$ will have $\cJ_\cA=\{J\in\cJ(\cC)\,|\,\theta(J)=1\}$. \textit{QED to Corollary 4.8}

\subsection{The strategy}

This subsection summarizes the steps used in identifying all exceptional quantum subgroups of  $\cC(\fg,k)$ when  $\fg=A_1,A_2,A_3,A_4$. In \cite{Gan-ii} we extend this to all $\fg$ of rank $\le 6$. The same methods would work for higher rank, though the analysis becomes increasingly computer intensive.

Consider for now any $\cC=\cC(\fg,k)$. Let $p$ be the smallest prime coprime to $f_\fg\ka$, where as always $\ka=k+h^\vee$ and $f_\fg$ is as in Section 2.1. Define \begin{equation}\label{Lmax}L_{\mathrm{max}}=L_{\mathrm{max}}(\fg,k)=\left\lfloor\frac{p(p-1)h^\vee\,\mathrm{dim}\,\fg}{6\ka}\right\rfloor\,.\end{equation}
Consider $\cV(\fg,k)$ with no conformal extension of Lie type (see \cite{BB,SW}). The previous subsection implies that if $\cA$ is exceptional, then there is some $\mu\in\cA$ with $2\le h_\mu\le L_{\mathrm{max}}$. In particular, if  $L>L_{\mathrm{max}}$, where $L$ is defined in Theorem 4.6, then $\cC$ has no exceptional quantum subgroups. 
 
\medskip\noindent\textbf{Step 1: First pass.} Retain any level $k$ where $\cV(\fg,k)$ possesses an exceptional conformal extension of Lie type. Of the remaining $k$, retain only those  with $L_{\mathrm{max}}(\fg,k)\ge2$. 

\medskip  As we explain next subsection, this is straightforward.  Only finitely many levels survive; Proposition 4.9 below estimates the number of survivors  and the size of the largest survivor  as rank$(\fg)\to\infty$. The survivors for $\fg=A_1,A_2,A_3,A_4$  are collected in Table 4.1. In \cite{Gan-ii} we list all survivors for all simple $\fg$ of rank $\le8$.

\medskip\noindent\textbf{Step 2: Second pass.} Eliminate any surviving $k$ with $\fC_k(\fg)\subseteq\cJ(\cC)$, or  with no  $\mu\in\fC_k(\fg)\cap\cC_0$ having $2\le h_\mu\le L_{\mathrm{max}}$.\medskip

The candidates $\fC_k(\fg)$ are defined in Definition 3.10. Call a level \textit{suspicious} if it survives Step 2.
 For $\fg$ of rank $\le 6$, the suspicious levels are collected in Table 1.1. We find that the numbers of suspicious levels  grows like $r^2$, where $r$ is the rank of $\fg$.

For $\fg=A_1,A_2,A_3,A_4$ it is elementary to complete Step 2. 
Suppose $(\fg,k)$ survives Step 1, and define $L_{\mathrm{max}}$ by \eqref{Lmax}. Run through the $\mu\in P_{\!++}^{\,\ka}(\fg)$ with $2\le h_\mu\le L_{\mathrm{max}}$, using \eqref{parityform} to see if any are candidates.  This bound on $h_\mu$ significantly limits the corner of $P_{\!++}^{\,\ka}(\fg)$ which must be searched.

The symmetries of the extended Dynkin diagram of $\fg$ also help here (and elsewhere): $\epsilon_\sigma(\la)=\epsilon_\sigma(C^iJ\la)$ for any simple-current $J$ and any $i=0,1$ and $\sigma$, and also $h_{C\la}=h_\la$ and $\theta({J\la})$ is related to $\theta(\la)$ through \eqref{scform}.

For example, the three worst cases in Table 4.1 are $A_4$ at levels 625, 835 and 1050, with $L_{\mathrm{max}}=3,2,2$ respectively. The cardinalities $\|P_{\!++}^{\,\ka}(A_4)\|$ range from around 6 billion to 50 billion for those three levels, but the numbers of $C_a$-orbits of $\la$ with integral conformal weight $h_\la$ between 2 and $L_{\mathrm{max}}$ are precisely 332, 181, and 212 respectively, permitting a fast search. Indeed, in \cite{Gan-ii} we go significantly further.

 \medskip\noindent\textbf{Step 3: For each suspicious $k$, identify the exceptional Res$\,\cA$.}\medskip

 By Corollary 3.11, Res$\,\cA$ will be a direct sum of candidates. We collect these
 candidates in Table 5.1. 
The largest part of  Step 3  involves eliminating most suspicious levels. Step 3 for our $\fg$ is described in detail in Section 5.1.

 \medskip\noindent\textbf{Step 4: Determine the branching rules for those Res$\,\cA$.}\medskip

At this point, we have two or three possibilities for Res$\,\cA$ for each $\fg$. In Section 5.2 we identify  how many simple local $\cA$-modules $M$ each has, and work out  Res$\,M$ for each of them. The results are collected in Table 2.1.

 \medskip\noindent\textbf{Step 5: Prove existence and uniqueness for those exceptional $\cA$.}\medskip

For $\fg=A_1,...,A_4$, existence is clear, since all $\cA$ are  of Lie type. Uniqueness for any $\cA$ whose restriction Res$\,\cA$ matches one of Lie type, is studied in Section 5.3 for any $\cC(\fg,k)$.

Section 5.4 applies the theory of Section 4.4 to Table 2.1 to classify all exceptional quantum subgroups for $\cC(A_r,k)$ when $k\le 5$ and $r$ is arbitrary. The results are in Table 2.2. 

\subsection{The first step and the proof of the Main Theorem}

In this subsection we explain how to use Section 4.1 to obtain for any $\fg$ a sparse set of levels where all exceptional extensions must lie. Specialized to $\fg=A_1,...,A_4$ the result is Table 4.1. In the process we prove the Main Theorem of the Introduction. 

Fix any simple Lie algebra $\fg$. Recall $f_\fg\in\{1,2,3\}$ from Section 2.1 and $L_{\mathrm{max}}$ from \eqref{Lmax}. 
Theorem 4.6  says   exceptional quantum subgroups of $\cC(\fg,k)$ occur  only if $L_{\mathrm{max}}\ge 2$ or $\cV(\fg,k)$ has a Lie-type conformal extension.

Table 4.1 collects all such $(\fg,k)$, for our  $\fg$. The total number for each $\fg$ is given in the final column. When the Lie algebra has $f_\fg>1$, the $f_\fg$th column won't contribute. Each entry $a..b$ in the table corresponds to a sequence $\{a,a+\Delta,a+2\Delta,\ldots,b\}$. For example, $4..10$ in the $A_1$ row  means $\{4,7,10\}$ (since $f_{A_1}=2$), whereas $625..1045$ in the $A_4$ row means $\{625,835,1045\}$ (since $f_{A_4}=1$). 


\medskip\centerline{
{\vbox{\tabskip=0pt\offinterlineskip
  \def\tablerule{\noalign{\hrule}}
  \halign to 4.5in{
    \strut#&
    \hfil#&\vrule#&\hfil#&\vrule#&\vrule#&    
\hfil#&\vrule#&\hfil#&\vrule#&\hfil#&\vrule#&\hfil#&\vrule#&\hfil#&\vrule#&\hfil#&\vrule#&\hfil#&\vrule#&\hfil#&\vrule#&\hfil#
\tabskip=0pt\cr
&$\fg\,\,\,$&&\,\,$f_\fg$\,\,&\,&&\,\,$\Delta=1$\,\,&&\,\,$\Delta=2$\,\,&&\,$\Delta=6/f_\fg$\,&&\hfill$\,\,\Delta=30/f_\fg$\,\hfill&&\hfill$\,\Delta=210/f_\fg$\,\hfill&&\,total\,\cr
\tablerule&$A_1\,$&&\hfill 2\hfill&\,&&\hfill$1$\hfill&&\hfill--\hfill&&\hfill$4..10$\hfill&&\hfill$13..28$\hfill&&&&\hfill6\hfill\cr
\tablerule&$A_2\,$&&\hfill 1\hfill&\,&&\hfill1\hfill&&\hfill$3..9$\hfill&&\hfill$15..33$\hfill&&\hfill$57$\hfill&&\hfill$207$\hfill&&\hfill11\hfill\cr
\tablerule&$A_3\,$&&\hfill 2\hfill&\,&&\hfill$1..26$\hfill&&\hfill--\hfill&&\hfill$29..95$\hfill&&\hfill$101..206$\hfill&&\hfill$311..521$\hfill&&\hfill60\hfill\cr
\tablerule&$A_4$\,&&\hfill 1\hfill&\,&&\hfill$1..15$\hfill&&\hfill$17..55$\hfill&&\hfill$61..193$\hfill&&\hfill$205..415$\hfill&&\hfill$625..1045$\hfill&&\hfill69\hfill\cr
\noalign{\smallskip}}}}}

\noindent{\bf Table 4.1.} The levels $k$ surviving Step 1, for $\fg=A_1,...,A_4$. $\Delta$ is the increment in the sequence $a..b$.\medskip

Let $p_n$ denote the $n$th prime (so $p_1=2$, $p_5=11$ etc). Write $A=h^\vee\mathrm{dim}\,\fg/12$ and $P_j=\prod_{{ i\le j\atop p_i\ne f_\fg}}p_i$. 

\medskip\noindent\textbf{Proposition 4.9.} \textit{Fix $\fg$. There is an  $l\ge 5$ such that $P_l> A\,(p_l^2-p_l)$. Write $\kappa_{max}(\fg)=P_l$ for the smallest such $l$.} 

\smallskip\noindent\textbf{(a)} \textit{If $(\fg,k)$ has $L_{\mathrm{max}}(\fg,k)\ge 2$, then $\kappa\le\ka_{max}(\fg)$. }

\smallskip\noindent\textbf{(b)} \textit{$\kappa_{max}(\fg)=O(r^{3+\epsilon})$ as $r=\mathrm{rank}\,\fg\to\infty$,  for any $\epsilon>0$.}

\medskip\noindent\textit{Proof.} Let $k$ have $L_{\mathrm{max}}\ge 2$. For each $j$ with  $p_j\ne f_\fg$, we get a threshold $\ka_j=A\,(p_j^2-p_j)$: if $\kappa>\ka_j$ then by definition of $L_{\mathrm{max}}$,  $p_j$ (hence $P_j$) must divide $\kappa$. For each $j=2,3,...,$ where $p_{j-1}\ne f_\fg$, the $j$th column of Table 4.1 contains all $k=\kappa-h^\vee$ for which $\ka$ is a multiple of $P_{j-1}$, and $\ka_{j-1}<\ka\le\ka_j$. In the first column  put all $k$ with $\ka\le\ka_1$ if $f_\fg\ne 2$ and  $\ka\le\ka_2$ otherwise. If $p_{j-1}=f_\fg$, put $-$ in the $j$th column. Table 4.1 then lists all such levels $k$ for those $A_r$. (We also include in Table 4.1 all $k$ for which $\cV(\fg,k)$ has $L_{\mathrm{max}}=1$ and an exceptional conformal extension of Lie type -- e.g.\ $\cV(A_1,10)$.)

From \cite{Loo} we know that $p_{i+1}/p_i\le 4/3$ for $i\ge 5$. From this we get $p_{i+1}-1\le 4p_i/3-4/3$,
and hence $\ka_{j+1}/\ka_j=(p_{j+1}/p_j)(p_{j+1}-1)/(p_j-1)<2$. However, $P_{j+1}/P_j=p_{j+1}>2$ so 
there must exist an $l\ge 5$ such that $P_l> \ka_l$. Then for all $m>l$, $P_m>\ka_m$.

Consider $\ka>P_m$ for some $m\ge l$, and suppose for contradiction that it has $L_{\mathrm{max}}\ge 2$. Then $\kappa$ must be a multiple of  $P_m$, so $\kappa\ge 2P_m\ge2\ka_m>\ka_{m+1}$, which means that 
 $ P_{m+1}|\ka$. But $\ka\ge P_{m+1}=p_{m+1}P_m\ge p_{m+1}\ka_m>4\ka_m>\ka_{m+2}$, which forces  $\kappa\ge P_{m+2}>P_{m+1}$. We have proved that if $\ka>P_m$, then also $\ka>P_{m+1}$. Continuing inductively, this impossibility gives us (a).

To see (b), it suffices to consider $f_\fg=2$ (since $f_\fg\le 2$ for rank $>2$, and $f_\fg=1$ grows slower). Note first that trivially $\sum_{3\le p\le x}\mathrm{log}\,p=(\pi(x)-1)\mathrm{log}\,x-\sum_{3\le p\le x }(\mathrm{log}\,x-\mathrm{log}\,p)$, where $p$ denotes primes and $\pi(x)$ is the number of primes $p\le x$.  By the Prime Number Theorem (see e.g.\ Chapter 9 of \cite{Hua}), $\pi(x)\,\mathrm{log}\,x=x+o(x)$, or equivalently for any $\epsilon>0$ there is an $X_\epsilon$ such that $(1-\epsilon)x<(\pi(x)-1)\mathrm{log}\,x<(1+\epsilon)x$ for all $x>X_\epsilon$; explicitly, we have for instance $(\pi(x)-1)\mathrm{log}\,x<1.3x$ for all $x\ge3$. Hence $\sum_{3\le p\le x}(\mathrm{log}\,x-\mathrm{log}\,p)=\sum_{3\le p\le x}\int^x_p t^{-1}\,dt=\int^x_3t^{-1}(\pi(t)-1)\,dt<1.3\int^x_3(\mathrm{log}\,t)^{-1}\,dt<16x/\mathrm{log}\,x$, since
$(\log\,t)^{-1}\le c\,((\log\,t)^{-1}-(\log\,t)^{-2})$ for $t\ge 3$ when $c\ge\log\,3/(\log\,3-1)$. Specializing  to $x=p_n$, what we've shown is that for any $\epsilon>0$,  there exists an $N_\epsilon$ such that $|1-\frac{1}{p_n}\sum_{i=2}^n\mathrm{log}\,p_i|<\epsilon $ for all $n>N_\epsilon$. The Prime Number Theorem  implies likewise that $p_n\sim n\,\mathrm{log}\,n$ for $n$ large. Together, these imply that
$\sum_{i=2}^n\log\,p_i\sim n\,\log\,n$.

Now, $\ka_{max}(\fg)=\prod_{2\le i\le l}p_i= A\,(p_l^2-p_l)e^b$ for some $0\le b<\log\,p_l$.
 Since $A$ like $(\rho|\rho)$ grows like $const\cdot r^3$ as $r\to\infty$, 
 we find $\log\,\kappa_{max}(\fg)\sim l\,\log\,l\sim 3\log(r)$. This gives (b). \textit{QED to Proposition 4.9}\medskip

Of course, the total number of $k$ with $L_{\mathrm{max}}\ge2$ is much smaller than $\ka_{max}$. An  upper estimate of that number when $f_\fg=1$ is $\ka_2+(\ka_3-\ka_2)/2+(\ka_4-\ka_3)/6+\cdots=A\,(p_1^2-p_1+(p_2^2-p_2-p_1^2-p_1)/2+\cdots)<7.5\,A$ (it is an upper estimate because, although we should add $l$ 1's to the expression, the first term $\ka_2$ has an extra $h^\vee$ which grows like $r$ or $2r$, much larger than $l\sim\log\,r/\log\,\log\,r$). Likewise, for $f_\fg=2$ this becomes $\ka_3+(\ka_4-\ka_3)/3+(\ka_5-\ka_4)/15+\cdots=A\,(p_2^2-p_2+(p_3^2-p_3-p_2^2-p_2)/2+\cdots)<13\,A$. Hence for sufficiently  large \textit{even} $r$,   the number of such $k$ for $\fg=A_r$ is less than $5r^3/4$, whilst that for large \textit{odd} $r$ is  less than $13r^3/6$. For $\fg=D_r$ when $r\equiv 0,1$ (mod 4), this number is $5r^3$, whilst for other $r$ is $26r^3/3$. For $B_r$ and $C_r$ at any $r$, the asymptotic upper bounds are $26r^3/3$ resp.\  $13r^3/3$.

\subsection{Level-rank duality and quantum subgroups}

In this subsection we explain how the quantum subgroups for some $\cC(A_r,k)$ are in natural bijection with those of $\cC(A_{k-1},r')$ where $r'=r+1$. 
An application of level-rank duality to quantum subgroups was first made in \cite{Wal}; we build instead on the formulation of \cite{OsS}.

Recall the adjoint subcategory $\cC(A_r,k)_0$ from Section 4.1.  
Then (Theorem 5.1 of \cite{OsS}) $\cC(A_r ,k)_0$ is braided tensor equivalent to $\cC(A_{k-1} ,r')_0^{\mathrm{rev}}$, using  $\la\leftrightarrow \tau_0(\la)^*$ where $\tau_0$ is defined in equation (3) of \cite{OsS}.  Hence their 
quantum subgroups are in natural bijection. The level-rank dual of a simple-current extension is always a simple-current extension, and the level-rank dual of a Lie-type conformal embedding is usually also of Lie-type.

For example, the level-rank duality  $\cC(A_2,5)_0\leftrightarrow\cC(A_4,3)_0$ sends e.g.\ $\mathbf{0} \leftrightarrow \mathbf{0}$, $(1,1) \leftrightarrow (1,0,0,1)$, $(2,2) \leftrightarrow(0,1,1,0)$, and $(3,0) \leftrightarrow (2,0,1,0)$. The exceptional quantum subgroup $\mathbf{0}\oplus(2,2)$ of $\cC(A_2,5)$ lies in $\cC(A_2,5)_0$, and corresponds  to the exceptional quantum subgroup $\mathbf{0}\oplus (0,1,1,0)$ of $\cC(A_4, 3)$.

Applying level-rank duality (Theorem 5.1 of \cite{OsS}) to Theorem 4.7 immediately gives the main result of this subsection:

\medskip\noindent\textbf{Corollary 4.10}. \textit{$\cC(A_r,k)$ has an exceptional quantum subgroup iff $\cC(A_{k-1}
,r')$ does.}\medskip

For most $r,k$ we can say more.
For any $\la\in P_{\!++}^{\,\kappa}(A_r)$, \eqref{Anorm} gives \begin{equation}(\la|\la)=-\frac{t(\la)^2}{r'}+\sum_{i=1}^ri\la_i^2+2\sum_{1\le i<j\le r}i\la_i\la_j\in-\frac{t(\la)^2}{r'}+t(\la)+2\bbZ\,.\label{Anormsimpl}\end{equation}
Using this we can prove:

\medskip\noindent\textbf{Proposition 4.11}. \textit{If either $r'$ or $r'/2$ is square-free, then the quantum subgroups of $\cC(A_r,k)$ are in natural bijection with a subset of those of $\cC(A_{k-1},r')$, the bijection given by $\tau_0$. If also $k$ or $k/2$ is square-free, then the quantum subgroups of $\cC(A_r,k)$ and $\cC(A_{k-1},r')$ are in natural bijection.}\medskip

\noindent\textit{Proof.} Let $\cA$ be any quantum subgroup of $\cC(A_r,k)$, and choose any $\la\in\cA$. Then $\theta(\la)=1$ so $(\la|\la)\equiv(\rho|\rho)$ (mod $2\ka$). Then \eqref{Anormsimpl} implies $-\frac{t(\la)^2}{r'}+t(\la)\equiv-\frac{t(\rho)^2}{r'}+t(\rho)$ (mod 2), or writing $\la=\la'+\rho$ this collapses to ${t(\la')^2}\equiv rr't(\la')$ (mod $2r'$). If  either $r'$ or $r'/2$ is square-free, this says $r'$ must divide $t(\la')$. In other words, $\cA\in\cC(A_r,k)_0$.  Hence $\tau_0$ identifies the connected \'etale algebras of $\cC(A_r,k)$ with those of $\cC(A_{k-1},r')_0$. \textit{QED to  Proposition 4.11}\medskip

We can identify the MFC (Rep$_{\cC(\fg,k)}\tau_0\cA)^{\mathrm{loc}}$ in terms of that of that of $\cA\in\cC(A_{k-1},r')_0$, as follows. By Remark 5.3 of \cite{OsS}, there is  a simple-current \'etale algebra  $\cB$ such that $\cC(A_r,k)$ is braided  equivalent to $(\mathrm{Rep}_{\cC(A_{k-1},r')^{\mathrm{rev}}\boxtimes \cC(\mathfrak{sl}_{r'k},1)}\cB)^{\mathrm{loc}}$. 
Proposition 3.16 of \cite{DMNO} tells us that $\tau_0(\cA)$ corresponds under the braided equivalence to some \'etale algebra $\tilde{\cA}$ of $\cC(A_{k-1},r')^{\mathrm{rev}}\boxtimes \cC(\mathfrak{sl}_{r'k},1)$ containing $\cB$ as a subalgebra. 
Theorem 3.6 of \cite{DNO} characterizes all \'etale algebras in a Deligne product: $\tilde{\cA}$  corresponds to  quantum subgroups $\tilde{\cA}_1\in\cC(A_{k-1},r')$ and $\tilde{\cA}_2\in\cC(\mathfrak{sl}_{r'k},1)$, fusion subcategories $\cD\subseteq (\mathrm{Rep}_{\cC(A_{k-1},r')}\,\tilde{\cA}_1)^{\mathrm{loc}}$ and $\cD'\subseteq (\mathrm{Rep}_{\cC(\mathfrak{sl}_{r'k},1)}\,\tilde{\cA}_2)^{\mathrm{loc}}$, and a braided equivalence $\cF:\cD\to\cD'$. Since
 $\cC(\mathfrak{sl}_{r'k},1)$ is pointed,   $\tilde{\cA}_2$ must be a sum of simple-currents; hence $(\mathrm{Rep}_{ \cC(\mathfrak{sl}_{r'k},1)}\,\tilde{\cA}_2)^{\mathrm{loc}}$, $\cD'$ and $\cD$ are all pointed. Putting  this together, the MFC (Rep$_{\cC(A_r,k)}\tau_0\cA)^{\mathrm{loc}}$ 
 is a simple-current extension of the Deligne product of $((\mathrm{Rep}_{\cC(A_{k-1},r')}\cA)^{\mathrm{loc}})^{\mathrm{rev}}$ with a pointed category. 
Thus when $\cA\in\cC(A_{k-1},r')_0$, the MFC (Rep$_{\cC(A_r,k)}\tau_0\cA)^{\mathrm{loc}}$ is no more exotic than is  $(\mathrm{Rep}_{\cC(A_{k-1},r')}\cA)^{\mathrm{loc}}$.

\section{Classifying quantum subgroups for $A_1,...,A_4$}

Fix $\cC=\cC(\fg,k)$. Let $\cA$ be a quantum subgroup: we write
\begin{equation}\mathrm{Res}\,\cA=\bigoplus_{\la\in \fC_k(\fg)} Z_\la\la\label{Arestr}\end{equation}
where $\fC_k(\fg)$ are the candidates (Definition 3.10). Recall that a quantum subgroup $\cA$ is \textit{exceptional} if $\la\in\cA$ (i.e.\  $Z_\la>0$) for some non-simple-current $\la$. This section eliminates most levels from Table 1.1 (Steps 2 and 3), and determines the branching matrices $B$ for the remainder (Steps 3 and 4). We also establish uniqueness when Res$\,\cA$ matches an extension of Lie type (Step 5). Much of this section applies to arbitrary $\fg$, but eventually we focus on $\fg=A_1,...,A_4$. Higher rank $\fg$ requires more sophistication \cite{Gan-ii}.

\subsection{Step 3: Identifying Res$\,\cA$}

Let's begin with some generalities.  Recall  $\cJ_\cA$ from Proposition 3.7(b), and $Z_\cC(\cJ_\cA)$ is all $\lambda\in P^{\,\kappa}_{\!++}(\fg)$ with $\varphi_J(\lambda)=1$ $\forall J\in\cJ_\cA$.

\medskip\noindent\textbf{Proposition 5.1.} \textit{Let $\cA$ in \eqref{Arestr} be a connected \'etale algebra for a pseudo-unitary MFC $\cC$. Then:}

\smallskip\noindent\textbf{(a)} \textit{$Z_\la\le \mathrm{qdim}(\la)\,$;} 

\smallskip\noindent\textbf{(b)} \textit{$Z_{\la^*}=Z_{\la}\,$;}

\smallskip\noindent\textbf{(c)} \textit{If $J\in\cJ_\cA$, then $Z_J=1$ and $Z_{J\la}=Z_\la$ for all $\la\in \mathrm{Irr}(\cC)$. Also, $Z_\la=0$ unless $\lambda\in Z_\cC(\cJ_\cA)$.}
\medskip

Part (a) follows from the Frobenius-Perron bound \eqref{PFineq} applied to modularity \eqref{modinv} (see (1.6) in \cite{BE}). Part (b) follows from modularity and $S^2=C$. Part (c) is Proposition 3.7.

 Let $\cA$ be any quantum subgroup of $\cC=\cC(\fg,k)$. Let $\cG_\cA$ be the group generated by $C$ and $\cJ_\cA$. Proposition 5.1 permits us to refine \eqref{Arestr}: 
  \begin{equation}\label{algshape}\cA=\lan\b1\ran\oplus\bigoplus_{\lan\la\ran\in\fC(\cJ_\cA)}Z_\la \lan\la\ran\,,\end{equation}
where $\lan\la\ran$ denotes the orbit of $\la$ with respect to $\cG_\cA$, and $\fC(\cJ_\cA)$ denotes the collection of all $\cG_\cA$-orbits of non-simple-currents in $\fC_k(\fg)\cap Z_\cC(\cJ_\cA)$.

Note that, for any $\mu\in Z_\cC(\cJ_\cA)$,  
\begin{equation}\|\cJ_\cA\|S_{\b1,\mu}+\sum_{\lan\la\ran\in \frak{C}(\cJ_\cA)}\|\lan\la\ran\|\,Z_\la\,\mathrm{Re}\,S_{\la,\mu}\ge 0\ ,\label{rhorow}\end{equation}
with equality iff the $\mu$-column of the branching matrix $B$ is identically 0. Equation \eqref{rhorow} is proved by evaluating
\eqref{modinv} at the entry $(\b1^e,\mu)$, and using $S^e_{\b1^e, \la}>0$.  In \eqref{rhorow},  Re$\,S_{\la,\mu}$ appears  because  $Z_\la=Z_{\la^*}$ and $S_{\la^*,\mu}=\overline{S_{\la,\mu}}$.

 As we will see, for most $(\fg,k)$  in Table 1.1, $\mu\in P_{\!++}^{\,\ka}(\fg)$ can be found for which Re$\,S_{\la,\mu}\le 0$ for all $\la\in\fC_k(\fg)$. Then for such $\mu$,  \begin{equation}Z_\la\le \left\lfloor \frac{\|\cJ_\cA\|}{\|\lan\la\ran\|}\frac{S_{\b1,\mu}}{|\mathrm{Re}\,S_{\la,\mu}|} \right\rfloor\,.\label{upperbd}\end{equation}
 Table 5.1 gives such $\mu$, when they exist, together with all candidates. When all but one orbit $\lan\la\ran\in\fC(\cJ_\cA)$ has Re$\,S_{\la,\mu}\le0$, then a similar argument gives a lower bound on $Z_\la$ if some of the other $Z_{\la'}$ are known to be nonzero.

Floating point calculation of $S_{\la,\mu}$ suffices because $Z_\la$ are bounded (Proposition 5.1(a)). For example, the condition Re$\,S_{\la,\mu}\le0$ in the previous paragraph can be replaced by  Re$\,S_{\la,\mu}<\epsilon$ where $\epsilon=S_{\b1,\mu}/(\|\fC(\cJ_\cA)\|S_{\la,\b1}/S_{\b1,\b1})$, since pseudo-unitary MFC obey the Perron-Frobenius inequality \eqref{PFineq}.

\medskip
\centerline{\scriptsize{\vbox{\tabskip=0pt\offinterlineskip
  \def\tablerule{\noalign{\hrule}}
  \halign to 7.15in{
    \strut# &
    \tabskip=0em plus1em &    
    \hfil#&
    #&\hfil#&\vrule\vrule#&    
\hfil#&\vrule#&\hfil#&\vrule#&\hfil#&\vrule#
\tabskip=0pt\cr
&&$\fg$&&\hfill $k$\hfill && 
\hfill candidates $\fC_k(\fg)$\hfill&&\hfill $\mu$ in  \eqref{rhorow}\hfill&&\hfill exceptional $\cA$ \hfill&\cr
\tablerule\tablerule&&$A_1\,$&&\hfill$\,10$\hfill&&\hfill$\,\bo,(6)$\hfill&&\hfill$\,(1)$\hfill&&\hfill$\bo\oplus(6)$\hfill&\cr
\tablerule&&&&\hfill$28$\hfill&&\hfill$\lan\bo,(10)\ran_2$\hfill&&\hfill$(2);\,(3);\,(1)$\hfill&&\hfill$\lan\bo, (10)\ran_2$\hfill&\cr
\tablerule\tablerule&&$A_2\,$&&\hfill5\hfill&&\hfill$\bo,(2\,2)$\hfill&&\hfill $(0\,1)$\hfill&&\hfill$\bo\oplus(2\,2)$\hfill&\cr
\tablerule&&&&\hfill9\hfill&&\hfill$\lan\bo,(4\,4)\ran_3$\hfill&&\hfill $(0\,3);\,(0\,1)$\hfill&&\hfill$\lan \bo,(4\,4)\ran_3$\hfill&\cr
\tablerule&&&&\hfill21\hfill&&\hfill$\lan\bo,(4\,4),(6\,6),(a\,a)\ran_3$\hfill&&\hfill $(03);\,(11);\,(22);\,(12)$\hfill&&\hfill$\,\lan\bo,(44),(66),(aa)\ran_3$\hfill&\cr
\tablerule&&&&\hfill57\hfill&&\hfill$\lan\bo,(a\,a),\lan(i\,i),(s\,s)\ran_3$\hfill&&\hfill $(1\,4)$\hfill&&\hfill$\emptyset$\hfill&\cr
\tablerule\tablerule&&$A_3$&&\hfill4\hfill&&\hfill$\lan\bo\ran_2,\lan(0\,1\,2)\ran_c$\hfill&&\hfill$(002)$\hfill&&\hfill$\lan\bo,(012)\ran_2$\hfill&
\cr\tablerule&&&&\hfill6\hfill&&\hfill$\lan\bo,(2\,0\,2)\ran_2$\hfill&&\hfill$(020);(011);(001)$\hfill&&\hfill$\lan\bo,(202)\ran_2$\hfill&
\cr\tablerule&&&&\hfill8\hfill&&\hfill$\lan\bo,(1\,2\,1)\ran_4$\hfill&&\hfill(see \S5.1)\hfill&&\hfill$\lan\bo,(121)\ran_4$\hfill&
\cr\tablerule&&&&\hfill10\hfill&&\hfill$\lan\bo,(4\,0\,4)\ran_2,\lan(0\,1\,6)\ran_{2c}$\hfill&&\hfill$(010)$\hfill&&\hfill$\emptyset$\hfill&
\cr\tablerule&&&&\hfill11\hfill&&\hfill$\bo,\lan(2\,0\,6)\ran_c$\hfill&&\hfill$(002)$\hfill&&\hfill$\emptyset$\hfill&
\cr\tablerule&&&&\hfill12\hfill&&\hfill$\lan\bo\ran_2,\lan(0\,5\,2)\ran_c$\hfill&&\hfill$(002)$\hfill&&\hfill$\emptyset$\hfill&
\cr\tablerule&&&&\hfill14\hfill&&\hfill$\lan\bo,(3\,0\,3),(6\,0\,6)\ran_2,\lan(4\,4\,0),(c\,0\,0)\ran_{2c}$\hfill&&\hfill$(030)$\hfill&&\hfill$\emptyset$\hfill&
\cr\tablerule&&&&\hfill16\hfill&&\hfill$\lan\bo,(1\,6\,1)\ran_4$\hfill&&\hfill$(012)$\hfill&&\hfill$\emptyset$\hfill&
\cr\tablerule&&&&\hfill18\hfill&&\hfill$\lan\bo,(8\,0\,8)\ran_2,\lan(2\,2\,6)\ran_{2c}$\hfill&&\hfill$(012)$\hfill&&\hfill$\emptyset$\hfill&
\cr\tablerule&&&&\hfill20\hfill&&\hfill$\lan\bo,(2\,9\,0)\ran_2,\lan(1\,0\,9)\ran_{2c}$\hfill&&\hfill$(006)$\hfill&&\hfill$\emptyset$\hfill&
\cr\tablerule&&&&\hfill26\hfill&&\hfill$\lan\bo,(3\,2\,3),(7\,2\,7),(c\,0\,c)\ran_2,\lan(0\,0\,k),(8\,1\,a)\ran_{2c}$\hfill&&\hfill$(01G)$\hfill&&\hfill$\emptyset$\hfill&
\cr\tablerule&&&&\hfill32\hfill&&\hfill$\lan\bo,(9\,0\,9)\ran_4$\hfill&&\hfill$(020)$\hfill&&\hfill$\emptyset$\hfill&
\cr\tablerule&&&&\hfill38\hfill&&\hfill$\lan\bo,(i\,0\,i)\ran_2,\lan(0\,3\,6),(3\,\,0\,z)\ran_{2c}$\hfill&&\hfill$(018)$\hfill&&\hfill$\emptyset$\hfill&
\cr\tablerule&&&&\hfill86\hfill&&\hfill$\lan\bo,(f\,0\,f),(r\,0\,r),(42\,\,0\,\,42)\ran_{2}$\hfill&&\hfill$(050)$\hfill&&\hfill$\emptyset$\hfill&
\cr\tablerule\tablerule&&$A_4$&&\hfill3\hfill&&\hfill$\bo,(0\,1\,1\,0)$\hfill&&\hfill$(0001)$\hfill&&\hfill$\bo\oplus(0110)$\hfill&
\cr\tablerule&&&&\hfill5\hfill&&\hfill$\lan\bo,(0\,2\,2\,0)\ran_5$\hfill&&\hfill$(0021);(0001)$\hfill&&\hfill$\lan\bo,(0220)\ran_5$\hfill&
\cr\tablerule&&&&\hfill7\hfill&&\hfill$\bo,(0\,3\,3\,0),(2\,0\,0\,2),(2\,1\,1\,2),\lan(0\,1\,6\,0),(0\,4\,0\,3),(1\,0\,1\,4)\ran_c$\hfill&&\hfill(see \S5.1)\hfill&&\hfill 1 (see Table 2.1)\hfill&
\cr\tablerule&&&&\hfill9\hfill&&\hfill$\bo,(0\,4\,4\,0),\lan(0\,0\,2\,6)\ran_c$\hfill&&\hfill$(0011)$\hfill&&\hfill$\emptyset$\hfill&
\cr\tablerule&&&&\hfill10\hfill&&\hfill$\lan\bo,(2\,3\,3\,2)\ran_5$\hfill&&\hfill$(0005)$\hfill&&\hfill$\emptyset$\hfill&
\cr\tablerule&&&&\hfill11\hfill&&\hfill$\bo,(0\,5\,5\,0),(4\,0\,0\,4),(4\,1\,1\,4),\lan(0\,3\,8\,0),(1\,2\,2\,6),(1\,4\,0\,4),(3\,0\,1\,6)\ran_c$\hfill&&\hfill$\,3(0010)\!\oplus\!(1314);\!(0401)$\hfill&&\hfill$\emptyset$\hfill&
\cr\tablerule&&&&\hfill13\hfill&&\hfill$\bo,\!(0660),\!\lan(0034),\!(0047),\!(010c),\!(0262),\!(0490),\!(0607),\!(1218),\!(3163)\ran_c$\hfill&&\hfill$(0413);(01a0)$\hfill&&\hfill$\emptyset$\hfill&
\cr\tablerule&&&&\hfill15\hfill&&\hfill$\lan\bo,(0\,7\,7\,0),(6\,0\,0\,6),(6\,1\,1\,6)\ran_5,\lan(1\,0\,2\,2),(3\,1\,0\,5)\ran_{5c}$\hfill&&\hfill$(0\,1\,4\,9)$\hfill&&\hfill$\emptyset$\hfill&
\cr\tablerule&&&&\hfill17\hfill&&\hfill$\bo,(0\,8\,8\,0),\lan(0\,6\,5\,2),(3\,6\,2\,6)\ran_c$\hfill&&\hfill$(0001)$\hfill&&\hfill$\emptyset$\hfill&
\cr\tablerule&&&&\hfill19\hfill&&\hfill$\bo,\!(0990),\!(2662),\!(4114),\!(4444),\!(8008),\!(8118),\!\lan(07c0),\!(2484),\!(4266),\!(701a)\ran_c$\hfill&&\hfill$(0303);(1118)$\hfill&&\hfill$\emptyset$\hfill&
\cr\tablerule&&&&\hfill21\hfill&&\hfill$\bo,(0\,a\,a\,0),\lan(2\,2\,6\,4),(3\,4\,0\,6)\ran_c$\hfill&&\hfill$(0007)$\hfill&&\hfill$\emptyset$\hfill&
\cr\tablerule&&&&\hfill23\hfill&&\hfill$\bo,(0bb0),(a00a),(a11a),\lan(09e0),(16d2),(2146),(2366),(2816),(2836),(901c)\ran_c$\hfill&&\hfill$(1212)$\hfill&&\hfill$\emptyset$\hfill&
\cr\tablerule&&&&\hfill25\hfill&&\hfill$\lan\bo,(0cc0),(2332),(2772)\ran_5,\lan(0505),(0607),(1006),(1077),(1161),(1365)\ran_{5c}$\hfill&&\hfill$\,(0055)\oplus(3082);(0055)$\hfill&&\hfill$\emptyset$\hfill&
\cr\tablerule&&&&\hfill31\hfill&&\hfill$\bo,(0ff0),(e00e),(e11e),(0cc0),(8338),(8448),$\hfill&&\hfill$(4505);(004b)$\hfill&&\hfill$\emptyset$\hfill&
\cr&&&&\hfill\hfill&&\hfill$\lan(00f0),(06ac),(0G0c),(37G5),(4h45),(604i),(635c),(364c)\ran_c$\hfill&&&&&
\cr\tablerule&&&&\hfill35\hfill&&\hfill$\lan\bo,(0\,h\,h\,0),(G\,0\,0\,G),(G\,1\,1\,G)\ran_5,\lan(0\,0\,0\,f),(1\,0\,f\,1)\ran_{5c}$\hfill&&\hfill$(0107)$\hfill&&\hfill$\emptyset$\hfill&
\cr\tablerule&&&&\hfill37\hfill&&\hfill$\bo,(0ii0),\lan(00aa),(00z0),(030G),(04GG),(080b),(090s),(0GL0)\ran_c,$\hfill&&\hfill$2(016L)\!\oplus\!(01w2)$\hfill&&\hfill$\emptyset$\hfill&
\cr&&&&\hfill\hfill&&\hfill$\lan(0n14),(1899),(262f),(274i),(453i),(637i),(71cf),(745a),(b02c)\ran_c$\hfill&&\hfill\hfill&&&
\cr\tablerule&&&&\hfill43\hfill&&\hfill$\bo,(0\,L\,L\,0),(8\,0\,0\,8),(8\,d\,d\,8),(c\,0\,0\,c),(c\,9\,9\,c),(k\,0\,0\,k),(k\,1\,1\,k)$\hfill&&\hfill$(0909)$\hfill&&\hfill$\emptyset$\hfill&
\cr\tablerule&&&&\hfill49\hfill&&\hfill$\bo,(0\,o\,o\,0),\lan(0\,5\,p\,A),(1\,4\,f\,J),(3\,8\,1\,7),(4\,c\,o\,5),(6\,3\,2\,s)\ran_c$\hfill&&\hfill$(001d)$\hfill&&\hfill$\emptyset$\hfill&
\cr\tablerule&&&&\hfill55\hfill&&\hfill$\lan\bo,\!(0rr0),\!(8118),\!(8ii8),\!(g11g),\!(gaag),\!(q00q),\!(q11q)\ran_5,\lan(000p),\!(006d),\!(10p1),\!(16d7)\ran_{5c}$\hfill&&\hfill$(0b07)$\hfill&&\hfill$\emptyset$\hfill&
\cr\tablerule&&&&\hfill85\hfill&&\hfill$\lan\bo,(0\,\,42\,\,42\,\,0)\ran_5,\lan(0\,c\,c\,A)\ran_{5c}$\hfill&&\hfill$(0005)$\hfill&&\hfill$\emptyset$\hfill&
\cr\tablerule&&&&\hfill115\hfill&&\hfill$\lan\bo,(0\,\,57\,\,57\,\,0),(k00k),(k\,37\,37\,k),(36\,\,0\,0\,\,36),(36\,LL\,36),(56\,\,0\,0\,\,56),(56\,\,1\,1\,\,56)\ran_5$\hfill&&\hfill$(0\,5\,0\,\,60)$\hfill&&\hfill$\emptyset$\hfill&
\cr\noalign{\smallskip}}}}}
  \noindent{\textbf{Table 5.1.} Step 3 data for $\fg=A_1,A_2,A_3,A_4$.}
  $\lan\cdot\ran_\alpha$ refers to the orbits  by permutations $\alpha$.  Letters a-z refer to numbers 10-35. Final column is Res$\, \cA$ for the exceptional quantum subgroups ($\emptyset$ when none exist).\medskip

 Table 1.1 lists all the levels $k$ (called \textit{suspicious} in Section 4.2) which remain after Step 2 for simple $\fg$ of rank $\le6$. Table 5.1 collects the data needed to complete Step 3 for the $2+4+14+21$ suspicious levels of $\fg=A_1,...,A_4$. Following the convention of Table 2.1,  it takes the weights from $P_+^k(\fg)$. To conserve space we use letters for numbers 10--35, and $\lan\cdot\ran$ notation for orbits. For example, for $\cC(A_3,26)$, $\lan(00k),(81a)\ran_{2c}$ denotes the orbits by $\lan J_a^2,C\ran$ of $(0,0,20)$ and $(8,1,10)$.

To see how to read the paper, suppose
for contradiction that $\cA$ is an  exceptional quantum subgroup of $\cC(A_2,57)$. By Corollary 4.8, it suffices to consider  $\cJ_\cA=\lan J_a\ran_3\cong\bbZ/3$.  Table 5.1 says that $\fC_{57}(A_2)=\langle\bo,\la_1,\la_2,\la_3\rangle_3$ where $\la_1=(10,10),\la_2=(18,18),\la_3=(28,28)$. Hence $\fC(\lan J_a\ran)=\{\lan\la_1\ran_3,\lan\la_2\ran_3,\lan\la_3\ran_3\}$. All quantum-dimensions $S_{\la,\b1}/S_{\b1,\b1}$ of candidates are between 1 and 2300, so accuracy of .01\% suffices. Table 5.1 gives  $\mu=(1,4)$; for it we compute 
$S_{\la,\mu}=.000746...,-.0141...,-.0427...,-.00725...$ respectively for $\la=\bo,\la_1,\la_2,\la_3$. Then  \eqref{upperbd} forces $Z_{\la_i}=0$  for all $i=1,2,3$, i.e.\ $\cA$ is not exceptional, so $\emptyset$ appears in the final column of Table 5.1.

When $\|\fC(\cJ)\|$ is relatively large, a $\mu\in Z_\cC(\cJ)$ satisfying Re$(S_{\la,\mu})\le0$ for all $\lan \la\ran\in\fC(\cJ)$ may not exist. But  \eqref{rhorow} implies
 \begin{equation} \|\cJ\|\sum_\mu x_\mu S_{\b1,\mu}+\sum_{\lan\la\ran\in\fC(\cJ)} \|\lan\la\ran\|\,Z_{\la}\sum_{\mu\in Z_\cC(\cJ)} x_\mu \mathrm{Re}\,S_{\la,\mu} \ge 0\,,\label{rhoineq2}\end{equation}
 whenever all $x_\mu\ge 0$.

 Semicolons in the $\mu$-column of Table 5.1 indicate that multiple $\mu$'s be considered, in the order given. When a linear combination of $\mu$'s is given,  use \eqref{rhoineq2}. For example, for $A_4$ level 11, first use $x_{(0010)}=3$ and $x_{(1314)}=1$ in  \eqref{rhoineq2} to show that five $Z_{\lan\la\ran}$ vanish, and then use $\mu=(0401)$ in  \eqref{upperbd} to show the remaining three $Z_{\lan\la\ran}$ vanish.

The higher rank $\fg$ considered in  \cite{Gan-ii}  require more sophisticated tools to do Step 3, and for even larger rank Step 3 will become a bottleneck.

We end this subsection by determining Res$\,\cA$ for the $2+3+3+3$ levels consistent with inequalities \eqref{rhorow},\eqref{rhoineq2}, i.e.\ those $k$ appearing in Table 2.1.

\medskip\noindent\textbf{$A_1$ level 10}\smallskip

\noindent By Table 5.1,  an exceptional quantum subgroup must look like $\cA=\bo\oplus Z_1(10)$ for some integer $Z_1>0$. Putting $\mu=(1)$ in \eqref{upperbd} forces $Z_1=1$ and we recover the  $\cA$ listed for $(A_1,10)$ in Table 5.1. The arguments for $(A_2,5)$ and $(A_4,3)$ are identical (use the $\mu$'s in the fourth column of  Table 5.1).

\medskip\noindent\textbf{$A_4$ level 5}\smallskip

\noindent An exceptional quantum subgroup $\cA$ looks like
 \begin{equation}\bo\oplus Z_1\lan(5000),(0500)\ran_c\oplus Z_2(0220) \oplus Z_2'\lan(1022)\ran_c\oplus Z_2''\lan(0102)\ran_c\nonumber\end{equation}
where $Z_2+Z_2'+Z_2''\ge 1$.   By Proposition 5.1(c), $Z_1\le1$.  Suppose first that $Z_1=1$. Then  $J_a\in\cJ_\cA$ so $Z_2=Z_2'=Z_2''\ge 1$.
Putting $\mu=(0021)$ in \eqref{upperbd} then gives $Z_2=1$ and we recover $\cA$ given in  Table 5.1.

Likewise, if $Z_1=0$, then $\mu=(0021)$  forces $Z_2=1$ and $Z_2'=Z_2''=0$. But $\cA=\bo\oplus (0220)$ cannot be an \'etale algebra, using $\mu=(0001)$.

The arguments for $(A_1,28)$, $(A_2,9)$, $(A_3,4)$ and $(A_3,6)$ are similar.

\medskip\noindent\textbf{$A_2$ level 21}\smallskip

\noindent Any quantum subgroup here has the shape
$$\bo\oplus Z_0\lan(21,\!0)\ran_c\oplus Z_1(4,\!4)\oplus Z_1'\lan(13,\!4)\ran_c\oplus Z_2(6,\!6)\oplus Z_2'\lan(9,\!6)\ran_c\oplus Z_3(10,\!10)\oplus Z_3'\lan(1,\!10)\ran_c$$
Because $(4,4),(6,6),(10,10)$ are Galois associates of $\bo$ (see Lemma 3.9(b)), $1=\cZ_{\bo,\bo}=\cZ_{(44),(44)}\ge Z_1^2$ etc, which requires $Z_1,Z_2,Z_3\le1$.

Consider first that $Z_0=1$, so $J_a\in\cJ_\cA$ and $Z_1=Z_1'$, $Z_2=Z_2'$, $Z_3=Z_3'$. Then $\mu=(0,3),(1,1)$ and (2,2) give $Z_2\ge Z_3$, $Z_1\ge Z_2$ and $Z_3\ge Z_1$ respectively, and we recover the $\cA$ of Table 5.1.

Otherwise, assume $\cA$ is an exceptional quantum subgroup with $Z_0=0$. Then we could extend $\cA$ by $\lan J_a\ran$ as in Corollary 4.8, and the result must be the algebra given in Table 5.1. This can only happen if $Z_1'=Z_2'=Z_3'=0$ and $Z_1=Z_2=Z_3=1$.  Now $\mu=(1,2)$ contradicts \eqref{rhorow}.

\medskip\noindent\textbf{$A_3$ level 8}\smallskip

\noindent An exceptional quantum subgroup looks like
$$\cA=\bo\oplus Z_1\lan(800)\ran_c\oplus Z_1'(080)\oplus Z_2(121)\oplus Z_2'(141)\oplus Z_3\lan(412)\ran_c$$
where $Z_1\le Z_1'\le 1$ and $Z_2+Z_2'+Z_3\ge 1$. Taking $\mu=(400)$, we find $2Z_1+Z_1'+1=Z_2+Z_2'+2Z_3$. 

Consider first $Z_1=1$, in which case $Z_1'=1$ and $Z_2=Z_2'=Z_3\ge 1$, and we recover  $\cA$ in Table 5.1. Next take $Z_1=0$ and $Z_1'=1$, so $Z_2=Z_2'$: comparing $\mu=(002),(010)$ gives  $Z_2=Z_3$, contradicting $2Z_1+Z_1'+1=Z_2+Z_2'+2Z_3$. Finally, if $Z_1=Z_1'=0$ then $Z_3=0$, which is eliminated by $\mu=(002)$.

\medskip\noindent\textbf{$A_4$ level 7}\smallskip

\noindent This is the most involved. A quantum subgroup looks like
$$\bo\oplus Z_1(0330)\oplus Z_2(2002)\oplus Z_3(2112)\oplus Z_4\lan(0160)\ran_c\oplus Z_5\lan(0403)\ran_c\oplus Z_6\lan(1014)\ran_c$$
 $\mu=(1112)$ implies $Z_4=0$ and $Z_1+Z_2\le2$. Comparing $\mu=(0310)$ and (1311), we now get $Z_2=Z_1\le1$, $Z_6=0$ and $Z_3\le1$. Then $\mu=(0003)$ requires $Z_2>0$, i.e.\ $Z_1=Z_2=1$, if $\cA$ is to be exceptional. $\mu=(0113)$ and (0016) then force $Z_5=Z_3=1$ as in Table 5.1.

\subsection{Step 4: Branching rules}

This subsection explains how to go from  Res$\,\cA$ to the modular invariant $\cZ=B^tB$, or more precisely the collection of branching rules $M\mapsto \mathrm{Res}\,M$ for all simple local $M$. The results for $\fg=A_1,...,A_4$ are collected in Table 2.1. This is done entirely at the combinatorial level -- next subsection we use this to identify the category (Rep$_\cC\,\cA)^{\mathrm{loc}}$. 

Some branching rules for Lie-type extensions  are known --- e.g.\  \cite{LeLib} gives those for $\cV(A_r,r-1)\subset\cV(A_{(r-1)(r+2)/2},1)$ and $\cV(A_r,r+3)\subset \cV(A_{r(r+3)/2},1)$. But we need to identify these branching rules using only Res$\,\cA$, since Theorem 5.2 identifies the extension using certain branching rules.

The branching rules only see the modules of $\cV^e$, i.e.\ the \textit{local} subcategory (Rep$_\cC\,\cA)^{\mathrm{loc}}$. Other aspects of  Rep$_\cC\,\cA$ aren't emphasized in this paper. For example, the fusion ring Fus$(\mathrm{Rep}_\cC\,\cA)$ is a module for Fus$(\cC(\fg,k))$, called the nim-rep, through the formula $\la.M=\mathrm{Ind}(\la)\otimes_\cA M$. 
The rich structure of many $\cC(A_r,k)$-module categories is studied e.g.\ in \cite{LR,CM,EMS}.

The MFC (Rep$_\cC\,\cA)^{\mathrm{loc}}$ is pseudo-unitary (see Corollary 3.6) so $S^e_{M,\b1^e}\ge S^e_{\b1^e,\b1^e}$. Unitarity of $S^e$ yields a test for completeness: if $M_i$ are pairwise inequivalent, they exhaust all  simple local $\cA$-modules (up to equivalence) iff
\begin{equation}\label{completeness}\sum(S^{e}_{M_i,\b1^e})^2=1\,.\end{equation}

Our main tool is modular invariance. Write \begin{equation}\cS_\nu(\mu):=\sum_{\la,M}B_{M,\nu}B_{M,\la} S_{\la,\mu}
=\sum_{\la,M}B_{M,\mu}B_{M,\la} S_{\la,\nu}\,,\label{SBBS}\end{equation} where the sum is over all simple $\cV^e$-modules $M$ (i.e.\ simples in (Rep$_\cC\,\cA)^{\mathrm{loc}}$), and over all simple $\cV$-modules $\la$.  The equality in \eqref{SBBS} is because $\cZ=B^tB$ commutes with $S$ (see \eqref{modinv}). Given $\cA$, call $\mu\in P^{\,\ka}_{\!++}(\fg)$ \textit{relevant} if $\cS_\b1(\mu)>0$ (this depends only on Res$\,\cA$).  Note that $\mu$ is relevant  iff it appears in the image Res$\,M$ for some simple local $M$, or equivalently iff the $\mu$-column of $B$ is not identically zero. 
Since $\cC(\fg,k)$ is pseudo-unitary, the definition of relevance is insulated from round-off error because $\cS_\b1(\mu)>0$ iff $\cS_\b1(\mu)\ge S_{\b1,\b1}$.

If $J\in\cJ(\cC)$, then Ind$(J)$ is invertible in Rep$_\cC\,\cA$. Given any simple $\cA$-module $M$, the simple $\cA$-module Ind$(J)\otimes_\cA M$ can be computed from \eqref{ResIndTimes}. $J$ is relevant iff $J\in\cJ_\cA$, in which case  Ind($J)\otimes_\cA M$ will be local iff $M$ is.

Other facts: $M$ is local iff the dual $M^*$ is. Moreover, Res$\,M^*=(\mathrm{Res}\,M)^*$.

If a simple $\cA$-module $M$ has $\la\in \mathrm{Res}\,M$, then Res$\,M$ is a subset of Ind(Res($\la$)). In particular, every simple $\mu\in\mathrm{Res}\,M$ lies in $\cA\otimes\la$. For local $M$, $\mu$ must in addition be relevant with  $\theta(\mu)=\theta(\la)$. For small $\la$, this yields an effective constraint on Res$\,M$.  

Call a relevant $\mu$ a \textit{singleton} if  there is a unique simple local $M$ such that $\mu\in\mathrm{Res}\,M$, and it occurs with multiplicity $B_{M,\mu}=1$. Equivalently,  $\cZ_{\mu,\mu}=1$.
 For a singleton $\la$, let $M(\la)$ denote the (unique) simple $\cV^e$-module with $\la$ in its restriction. If $M$ is any simple $\cV^e$-module, and $\la$ is a singleton, then \begin{equation}S^e_{M,M(\la)}=\sum_{\mu}B_{M,\mu}S_{\mu,\la}\,,\ \  \cS_\la(\mu)=\sum_MB_{M,\mu}S^e_{M,M(\la)}\,.\label{singleton}\end{equation}
If $\la,\mu$ are singletons, then $\cS_\la(\mu)=S^e_{M(\la),M(\mu)}$.

For example, $\b1$ is always a singleton, and hence $S^e_{\b1^e,\b1^e}=\cS_\b1(\b1)$. By pseudo-unitarity, if $\cS_\b1(\mu)<2\cS_\b1(\b1)$, then $\mu$ must also be a singleton.

By \eqref{GalBranch}, if $\la$ is a singleton, so is any Galois associate $\la^\sigma$, with $M(\la^\sigma)=M(\la)^\sigma$ and
Res$\,M(\la^\sigma)=\mathrm{Res}(M(\la))^\sigma$. For example, 
$M(\la^*)=M(\la)^*$. If $J\in\cJ_\cA$, then $J$ is a singleton with $M(J)=\mathrm{Ind}(J)$; in this case, $\la$ is a singleton iff $J\la$ is.

\medskip
\centerline{\scriptsize{\vbox{\tabskip=0pt\offinterlineskip
  \def\tablerule{\noalign{\hrule}}
  \halign to 5.65in{
    \strut# &
    \tabskip=0em plus1em &    
    \hfil#&
    #&\hfil#&\vrule\,\vrule#&    
\hfil#&\vrule#
\tabskip=0pt\cr
&&$\fg$&&\hfill $k$\hfill && \hfill relevant $\la$ grouped by $\theta\,[\cS_\b1]$\hfill&\cr
\tablerule\tablerule&&$A_1\,$&&\hfill$\,10$\hfill&&\hfill$1[.5]\!:\!\bo$,(6); \ $-1[.5]\!:\!(4),(10)$; \ $\xi_{16}^5[.707]\!:\!(3),(7)$\hfill&\cr
\tablerule&&&&\hfill$28$\hfill&&\hfill$1[.525]\!:\!\lan\bo,(10)\ran_2$; \ $\xi_5^2[.850]\!:\!\lan(6),(12)\ran_2$\hfill&\cr
\tablerule\tablerule&&$A_2\,$&&\hfill5\hfill&&\hfill$1[.408]\!:\!\bo,(22)$; \ $\xi_{12}^5[.408]\!:\!\lan(02),(23)\ran_c$; \ $\xi_3^2[.408]\!:\!\lan(05),(12)\ran_c$; \ $\xi_4^3[.408]\!:\!\lan(03)\ran_c$\hfill&\cr
\tablerule&&&&\hfill9\hfill&&\hfill$1[.577]\!:\!\lan\bo,(44)\ran_3$; \ $\xi_3^2[1.154]\!:\!\lan(22)\ran_3$\hfill&\cr
\tablerule&&&&\hfill21\hfill&&\hfill$1[.707]\!:\!\lan\bo,(44),(66),(10\,10)\ran_3$; \ $\xi_4^3[.707]\!:\!\lan(06),(47)\ran_{3c}$\hfill&\cr
\tablerule\tablerule&&$A_3$&&\hfill4\hfill&&\hfill$1[.5]\!:\!\lan\bo,(012)\ran_4$; \  $-1[.5]\!:\!\lan(004),(101)\ran_2$; \ $\xi_{16}^{15}[1.414]\!:\!(111)$\hfill&
\cr\tablerule&&&&\hfill6\hfill&&\hfill$1[.316]\!:\!\lan\bo,(202)\ran_2$; \ $\xi_{20}^9\!:\![.316]\lan(002)\ran_{2c},[.632](212)$; \ $\xi_5^4\!:\! [.316]\lan(012)\ran_{2c},[.632](303)$;\hfill\cr
&&&&&&\hfill$\,\xi_{20}\!:\![.316]\lan(123)\ran_{2c},[.632](030)$; \ $\xi_5\!:\![.316]\lan(004)\ran_{2c},[.632](121)$; \ $\xi_4[.316]\!:\!\lan(006),(022)\ran_2$\hfill&
\cr\tablerule&&&&\hfill8\hfill&&\hfill$1[.5]\!:\!\lan\bo,(121)\ran_4$; \ $-1[.5]\!:\!\lan(020),(303)\ran_4$; \ $\xi_4[1]\!:\!\lan(113)\ran_4$\hfill&
\cr\tablerule\tablerule&&$A_4$&&\hfill3\hfill&&\hfill$1[.316]\!:\!\bo,(0110)$; \ $\xi_{20}^9[.316]\!:\!\lan(0010),(0201)\ran_{c}$; \ $\xi_5^4[.316]\!:\!\lan(0030),(0101)\ran_{c}$;\hfill\cr
&&&&&&\hfill$\xi_{20}[.316]\!:\!\lan(0020),(1002)\ran_{c}$; \ $\xi_5[.316]\!:\!\lan(0003),(1011)\ran_{c}$; \ $\xi_4[.316]\!:\!\lan(0102)\ran_c$\hfill&
\cr\tablerule&&&&\hfill5\hfill&&\hfill$1[.5]\!:\!\lan\bo,(0220)\ran_5$; \ $-1\!:\![.5]\lan(1001)_5,[2.5](1111)$\hfill&
\cr\tablerule&&&&\hfill7\hfill&&\hfill$1[.258]\!:\! \bo,(0330),(2002),(2112),\lan(0403)\ran_c$; \ $\xi_5^4[.258]\!:\!(0007),(1033),(0023),(1121),(0040)\ran_c;$\hfill&\cr
&&&&&&\hfill$\,\xi_5[.258]\!:\!\lan(0070),(0103),(0232),(1212),(0030)\ran_c$; \ $\xi_3^2\!:\![.258]\lan(0005),(1213),(0042),(1022)\ran_c,[.516](0220)$;\hfill&\cr
&&&&&&\hfill$\xi_{15}^7\!:\! [.258]\lan(0002),(0052),(0312),(0121),(0421),(1004),(2012),(0222)\ran_c,[.516]\lan(2203)\ran_c$; \hfill&\cr
&&&&&&\hfill$\xi_{15}^{13}\!:\![.258]\lan(0250),(0025),(2103),(1031),(0012),(0124),(1222),(0122)\ran_c,[.516]\lan(0302)\ran_c$ \hfill&
\cr\noalign{\smallskip}}}}}
  \noindent{\textbf{Table 5.2.} Relevant $\la$ for $\fg=A_1,A_2,A_3,A_4$.} Grouped by common ribbon twist $\theta(\lambda)$ and
 $\cS_\b1(\la)$.  $\lan\cdot\ran_\alpha$ refers to the orbits  by permutations $\alpha$.\medskip

\medskip\noindent$\mathbf{A_1}$ \textbf{level 10}\smallskip

From Table 5.2 we see all $\cS_\bo(\la)<2\cS_\bo(\bo)$ and thus
all relevant $\la$ are singletons. We know $\cA=\bo\oplus(6)$. $(10)$ is a simple-current so Res$\,M(10)=J_a(\bo\oplus(6))=(10)\oplus(4)$ and $M(4)=M(10)$. Res$\,M(3)=(3)\oplus z(7)$ where $z=0,1$, but $1/\sqrt{2}=\cS_\bo(3)=1/(2\sqrt{2})+ z/(2\sqrt{2})$ so $z=1$ and $M(7)=M(3)$. 

 Theorem 5.2 below tells us (Rep$_{\cC(A_1,10)}\cA)^{\mathrm{loc}}\cong\cC(C_2,1)$. In order to use Theorem 5.2 to prove uniqueness of $\cA$ as an \'etale algebra, we need to identify $\Lambda_1^e$. But it is clear $M(10)$ is the simple-current $\Lambda^e_2$, so $M(3)$ is $\Lambda^e_1$. We recover the branching rules of Table 2.1.

\medskip\noindent$\mathbf{A_1}$ \textbf{ level 28}\smallskip


Again, all $\la$ are singletons. We know $\cA=\bo\oplus J_a\oplus (10)\oplus J_a(10)$. Since Ind$(J_a)=\bo^e$, $M(J_a(6))=M(6)$ and $M(J_a(12))=M(12)$. We have $M(6)=(6)\oplus J_a(6)\oplus z(12)\oplus zJ_a(12)$ where $z=0,1$. As in the argument for $(A_1,10)$, we get $z=1$ and $M(6)=M(12)$. Nothing more is needed for Theorem 5.2 to establish that Rep$_{\cC(A_2,28)}\,\cA\cong \cC(G_2,1)$ and that $\cA$ is unique. 

\medskip\noindent$\mathbf{A_2}$ \textbf{level 5}\smallskip


All $\la$ are singletons.  We know $\cA=\bo\oplus(2\,2)$ and thus $M(J_a)=(5\,0)\oplus(1\,2)$ and $M(J_a^2)=(0\,5)\oplus(2\,1)$. We get $M(0\,3)=(0\,3)\oplus(3\,0)=M(3\,0)$ by the usual argument.  Thus $M(2\,0)=J_a((0\,3)\oplus(3\,0))=(2\,0)\oplus(2\,3)$ and $M(0\,2)=((2\,0)\oplus(2\,3))^*=(0\,2)\oplus(3\,2)$.   For Theorem 5.2,  we can force $M(0\,2)\leftrightarrow\Lambda^e_1$ by hitting with the outer-automorphism $C$ if necessary. 

\medskip\noindent$\mathbf{A_2}$ \textbf{ level 9}\smallskip

Note that $S^e_{\bo^e,\bo^e}=\cS_{\bo}(\bo)=1/\sqrt{3}$. This means there are two remaining local simple $\cA$-modules (other than $\bo^e=\cA$) -- call them $M_1,M_2$ -- each with $S^e_{M_i,\bo^e}=1/\sqrt{3}$. After all, the MFC (Rep$_{\cC(A_2,9)}\cA)^{\mathrm{loc}}$ is pseudo-unitary (Corollary 3.6), so $S^e_{M_i,\bo^e}\ge S^e_{\bo^e,\bo^e}$ so its rank is at most 3; it can't be rank 2 because the nontrivial quantum-dimension $S^e_{M,\bo^e}/S^e_{\bo^e,\bo^e}$ would then have to be 1 or $(1+\sqrt{5})/2$ (see e.g.\ Theorem 3.1 in \cite{RSW}).

 Because $J_a\in\cA=\lan\bo,(4\,4)\ran_3$, we know $M_i=z_i(2\,2)\oplus z_i(5\,2)\oplus z_i(2\,5)$, as in $(A_1,28)$. As in $(A_1,10)$, we obtain $z_1=z_2=1$.  Res$\,\Lambda^e_1$ is now forced.

\medskip\noindent$\mathbf{A_2}$ \textbf{level 21}\smallskip

All relevant $\la$ are singletons. Since $S^e_{M(\la),\bo^e}=\cS_\bo(\la)=1/\sqrt{2}$ for any relevant $\la$, the pseudo-unitary MFC (Rep$_{\cC(A_2,21)}\cA)^{\mathrm{loc}}$ has rank 2, with simples $\cA=\bo^e$ and $M(0\,6)$. This suffices for Theorem 5.2.

\medskip\noindent$\mathbf{A_3}$ \textbf{level 4}\smallskip

The simple-current $J_a$ is a singleton, with $M(J_a)=J_a(\bo\oplus (040)\oplus (012)\oplus(210))=(400)\oplus(004)\oplus(101)\oplus(121)$. As $S^e_{\bo^e,\bo^e}=S^e_{M(J_a),\bo^e}=1/2$ and $\sqrt{2}=\sum_MB_{M,(111)}S^e_{M,\bo^e}$,  the pseudo-unitary MFC (Rep$_{\cC(A_3,4)}\cA)^{\mathrm{loc}}$ has rank 3. The missing simple local is $M(111)=z(111)$, with $S^e_{M(111),\bo^e}=1/\sqrt{2}$ and branching rule $B_{M(111),\bo}=z=2$. $\Lambda^e_1$ is the simple-current $M(J_a)$.

\medskip\noindent$\mathbf{A_3}$ \textbf{level 6}\smallskip

We know $M(\bo)=\cA$ and $M(J_a)=J_a\cA$ as usual. Consider next the relevant $\la$ with $\theta=\xi_{20}^9$. $(002)$ is a singleton and $J_a^2\in\cA$, so $M(002)=\lan(002)\ran_2\oplus z\lan(200)\ran_2\oplus z'(212)$, where $z=0,1$. We compute $S_{J_a^{2i}C^j(002),\bo}=.0437...$ and 
$S_{(212),\bo}=.228...$ so \eqref{singleton} becomes $.316...=(2+2z)\times.0437...+z'.228...$, and we obtain $z=0,z'=1$. Thus $M(200)=M(002)^*$. Note the Galois associates $7.(002)=(301)$ and $7.(212)=(030)$, so we also obtain $M(123)=M(301)=7.M(002)$ and $M(321)=M(123)^*$. Finally, $J_a$ applied to those give $M(012),M(210),M(004),M(400)$. The completeness test \eqref{completeness} then confirms that these 10 $\cA$-modules exhaust all simple local ones.  Without loss of generality (hitting with $C$ if necessary) we can assign $\Lambda^e_1$ to $M(200)$.

\medskip\noindent$\mathbf{A_3}$ \textbf{level 8}\smallskip

For $\theta=-1$, both $(020)$ and $(303)$ are singletons, and we obtain $M(020)=\lan(020)\ran_4\oplus\lan(303)\ran_4$ by the usual arguments. Analogously to $(A_3,4)$, the rank of (Rep$_{\cC(A_3,8)}\cA)^{\mathrm{loc}}$ must be 4; along with $M(\bo)=\cA$ and $M(020)$ we have $M_1,M_2$ with $M_i=z_i\lan(113)\ran_4$ and $S^e_{M_i,\bo}=1/2$; \eqref{singleton} then forces $z_1=z_2=1$. The connection to $\cC(D_{10},1)$ given in Table 2.1 is now forced.

\medskip\noindent$\mathbf{A_4}$ \textbf{level 3}\smallskip

All relevant $\la$ are singletons. We have $M(\bo)=\cA$ and hence $M(J_a^i)=J_a^i\cA$ for $i\le4$. Consider next $\theta=\xi_{20}^9$.  As explained earlier this subsection, if $B_{M(0010),\mu}>0$ then $\mu$ is relevant with $\theta(\mu)=\xi_{20}^9$ and $\mu\in\cA\otimes(0010)$. We find that the only such $\mu$ are  $(0010)$ and $(0201)$. Hence $M(0010)=(0010)\oplus z(0201)$, where $z=0,1$. We get $z=1$ as usual. The four remaining $M$ are of the form $J_a^iM(0010)$. As usual, we can force $\Lambda_1^e=M(0010)$. 

\medskip\noindent$\mathbf{A_4}$ \textbf{level 5}\smallskip

$(1001)$ is a singleton; then $M(1001)=\lan(1001)\ran_4\oplus z(1111)$. Using $S_{(1001),\bo}=1/20$ and $S_{(1111),\bo}=1/4$, we get $z=1$. As in $(A_3,8)$, the rank of (Rep$_{\cC(A_4,5)}\cA)^{\mathrm{loc}}$ must be 4; along with $M(\bo)=\cA$ and $M(1001)$ we have $M_1,M_2$ with $M_i=z_i(1111)$ and $S^e_{M_i,\bo}=1/2$; \eqref{singleton} then forces $z_1=z_2=2$. The assignment of $\Lambda^e_1\in\cC(D_{12},1)$ to $M(1001)$ is because $h^e_{\Lambda_1^e}=1/2=h_{(1001)}$ whilst $h^e_{\Lambda_{11}^e}=h^e_{\Lambda_{12}^e}=3/2=h_{(1111)}$. 

\medskip\noindent$\mathbf{A_4}$ \textbf{level 7}\smallskip

We obtain $M(J_a^i)=J_a^i\cA$ as usual. Consider the relevant $\la$ with $\theta=\xi_{15}^7$. All but $\lan(2203)\ran_c$ are singletons. The only singletons with $\theta=\xi_{15}^7$ in the fusion $(0002)\otimes\cA$ are $(0002),(0421),(2130),(2012),(2203)$, and thus $M(0002)=(0002)\oplus z(0421)\oplus z'(2130)\oplus z''(2012)\oplus z'''(2203)$ where $z,z',z''\le1$ and $z'''\le2$. Since $S_{\la,\bo}=.0057...,.0215...,.0803...,.0860...,.0645...$ respectively for those 5 relevant $\la$, we get that $z=z'=z''=z'''=1$. This gives us $M(J_a^iC^j(0002))$, and completeness confirms these 15 $\cA$-modules exhaust all that are simple and local. We can force $\Lambda_1^e=M(2000)$ as before.


\subsection{Step 5: Existence and uniqueness}

Section 5.1 shows that any exceptional quantum subgroup $\cA$ for  $\fg=A_1,...,A_4$ occurs at 11 possible pairs $(\fg,k)$, and that there is a unique Res$\,\cA$ for each pair.  Section 5.2  determines all possible branching rules for each $\cA$, finding that they are unique for each pair. Existence is clear from Table 2.1:  at least one exceptional quantum subgroup (namely one of Lie type) exists for each of those 11 pairs. Much more difficult is to establish uniqueness. This subsection addresses this, for any $\cC(\fg,k)$ when Res$\,\cA$ matches one  of Lie type (not just our 11 examples). As is permitted by Theorem 3.5, 
we work at the level of VOAs; equivalence of VOA extensions is defined in Section 3.3.

 {Lie-type conformal extensions} $\cV(\fg,k)\subset\cV(\fg^e,1)$ are defined in Section 3.2, and are classified in   \cite{BB,SW}.  We restrict to  simple $\fg,\fg^e$.  The homogenous spaces are $\cV(\fg,k)_1\cong\fg$ and $\cV(\fg^e,1)_1\cong\fg^e$ as Lie algebras. 

Embeddings (injective Lie algebra homomorphisms) $\iota_i:\fg\to\fg^e$ are studied in e.g.\ \cite{Dyn,Min}. We say $\iota_i$ are \textit{conjugate} if there is an inner automorphism $\alpha$ of $\fg^e$ such that $\iota_2=\alpha\circ\iota_1$. The notion relevant to us is slightly weaker: we call $\iota_i$ \textit{VOA-equivalent}  if  $\iota_2=\alpha\circ\iota_1$ for  any automorphism $\alpha$ of $\fg^e$. (Recall that the  automorphism group of the simple Lie algebra $\fg^e$ is the semi-direct product of its inner automorphisms with its outer ones. Its inner automorphisms are conjugations by the connected simply connected Lie group $G^e$ of $\fg^e$. The outer automorphisms correspond to the group of symmetries of the Dynkin diagram of $\fg^e$.)   

\medskip\noindent\textbf{Theorem 5.2.} \textit{Let $\cV(\fg,k)\subset \cV(\fg^e,1)$ be conformal  with $\fg,\fg^e$ both simple. Let $\cV(\fg,k)\subset \cW$ also be conformal, with $\cV(\fg^e,1)\cong\cW$ as $\cV(\fg,k)$-modules. Then $\cV(\fg^e,1)$ and $\cW$  are isomorphic as VOAs. Moreover, the extensions $\cV(\fg,k)\subset \cV(\fg^e,1)$ and $\cV(\fg,k)\subset \cW$ are equivalent, provided:}

\smallskip\noindent\textbf{(a)}  \textit{for $\fg^e\cong A_r,B_r,C_r,E_6$,  $\mathrm{Res}^{\cV(\fg^e,1)}_{\cV(\fg,k)}\,\Lambda^e_1\cong\mathrm{Res}^{\cW}_{\cV(\fg,k)}\,\Lambda^e_1$  as $\cV(\fg,k)$-modules;}

\smallskip\noindent\textbf{(b)} \textit{for $\fg^e\cong D_r$, $\mathrm{Res}^{\cV(\fg^e,1)}_{\cV(\fg,k)}\,\Lambda^e_1\cong\mathrm{Res}^{\cW}_{\cV(\fg,k)}\,\Lambda^e_1$  and $\mathrm{Res}^{\cV(\fg^e,1)}_{\cV(\fg,k)}\,\Lambda^e_r\cong\mathrm{Res}^{\cW}_{\cV(\fg,k)}\,\Lambda^e_r$   as $\cV(\fg,k)$-modules;}

\smallskip\noindent\textbf{(c)} \textit{for $\fg^e\cong E_7,E_8,G_2,F_4$, no other conditions are needed.}

\medskip\noindent\textit{Proof.} $\cV(\fg,k)$ is  generated as a VOA by its homogeneous space $\cV(\fg,k)_1\cong\fg$. $\cW_1$ is a reductive Lie algebra $\fg'$ (Theorem 1 of \cite{DM}). Let $\cW_{\mathrm{Lie}}$ be the subVOA of $\cW$ generated by  $\cW_1$. Since $\cW$ contains $\cV(\fg,k)$, $\fg$ is a Lie subalgebra of $\fg'$, so $\cV(\fg,k)\subset\cW_{{\mathrm{Lie}}}$. Thus the central charges satisfy $c(\cV(\fg,k))\le c(\cW_{{\mathrm{Lie}}})\le c(\cW)$. But $c(\cV(\fg,k))=c(\cW)$ by hypothesis, so $\cV(\fg,k)\subset\cW_{{\mathrm{Lie}}}\subset\cW$ are all conformal. The restriction of $\fg^e$ to $\fg$ is obtained from Res$\,\bo^e$ by deleting from it all simple $\cV(\fg,k)$-submodules with conformal weight $h_\la>1$, and replacing $\bo$ with the adjoint representation of $\fg$; as the restriction of $\fg'$ to $\fg$ is obtained similarly, the restrictions of $\fg'$ and $\fg^e$ agree. Hence the dimensions of $\fg^e$ and $\fg'$ are equal. Running through the list of \cite{BB,SW}, we see that only $(\fg,k)=(A_5,6),(A_7,10),(B_4,2)$ have multiple Lie-type conformal extensions, but  the dimensions of the extended Lie algebras disagree. Thus $\fg^e\cong\fg'$ as Lie algebras, and $\cW_{\mathrm{Lie}}\cong\cV(\fg^e,1)$ as VOAs. Since the restrictions of $\cW_{\mathrm{Lie}}$ and $\cW$ to $\cV(\fg,k)$ must agree,  $\cW=\cW_{\mathrm{Lie}}\cong\cV(\fg^e,1)$  as VOAs.

We need to show the associated Lie algebra embeddings $\iota_1:\fg\subset\fg^e$ and $\iota_2:\fg\subset\fg'$ are VOA-equivalent.
Following \cite{Dyn,Min}, we call embeddings $\iota_i$ \textit{linearly conjugate} if for any representation $\rho:\fg^e\to\mathfrak{gl}(V)$, the restrictions $\rho\circ\iota_i$ are equivalent as representations of $\fg$. By  Theorem 1.5 of \cite{Dyn} and Theorem 2 of \cite{Min}, our $\iota_i$ are linearly conjugate. More precisely, when $\fg^e\not\cong D_r$ for any $r$, it suffices to confirm this for $\rho$ being the first fundamental representation $\Lambda_1^e$. For $\fg^e\cong A_r,B_r,C_r,E_6$, $\Lambda_1^e\in P_{+}^{1}(\fg^e)$ so is included in the branching rules. For $\fg^e\cong E_7,E_8,G_2,F_4$, $\Lambda_1^e$ is the adjoint representation of $\fg^e$ and the restriction of that to $\fg$ equals the $h_\la=1$ part of Res$\,\bo^e$, as described last paragraph. When $\fg^e\cong D_r$, we see that  $\Lambda_1^e,\Lambda_r^e\in P_{+}^{1}(\fg^e)$ suffices.

Theorem 3 of \cite{Min} now tells us $\iota_i$ are conjugate (hence VOA-equivalent) when either $\fg\cong A_1$, or Im$\,\iota_i$ are both 
regular subalgebras of $\fg^e$, or $\fg^e\cong A_r,B_r,C_r,G_2,F_4$ for some $r$. In fact, if one of Im$\,\iota_i$ is a regular subalgebra, so must be the other, since $\fg\subset\fg^e$ is regular for a conformal extension $\cV(\fg,k)\subset\cV(\fg^e,1)$ iff rank$\,\fg=\mathrm{rank}\,\fg^e$. 

For $\fg^e=D_r$, Theorem 4 of \cite{Min} says that, because $\iota_i$ are linearly conjugate, if they are not conjugate then $\iota_2=\alpha\circ\iota_1$ where $\alpha$ is outer (i.e.\ not inner). But in this case $\iota_i$ are clearly VOA-equivalent (necessarily $\alpha$ here would permute $P_+^k(\fg^e)$ but preserve restriction).

The only conformal embeddings into $E_6$ for which linear conjugates aren't necessarily conjugates, are $\cV(A_2,9)\subset \cV(E_6,1)$ and $\cV(G_2,3)\subset \cV(E_6,1)$, but in both cases the embeddings coincide up to an outer automorphism of $E_6$ and so are VOA-equivalent (Theorem 11.1 in \cite{Dyn}, Corollary 3 of \cite{Min}). For conformal embeddings into $E_7$ and $E_8$, linear conjugacy implies conjugacy and hence VOA-equivalence. 
\textit{QED to Theorem 5.2}\medskip

For example, the two extensions of the form $\cV(D_8,1)\subset\cV(E_8,1)$ are distinguished by Res$\,\bo^e$.  The requirement that $\fg$ be simple is only used in showing that $\cW\cong\cV(\fg^e,1)$.

Thanks to Step 4, Theorem 5.2 proves that all exceptional quantum subgroups for $\fg=A_1,...,A_4$ are of Lie type, and are unique. Compare our argument to that of \cite{KO} for $\fg=A_1$, which established uniqueness by explicitly solving the associativity constraint in the category. The papers \cite{AL,DL} established the weaker statement for $\fg=A_1,A_2$ that these extensions define isomorphic VOAs.

\subsection{The level $k\le 5$ classification}

In this subsection we find all exceptional quantum subgroups of $\cC(A_r,k)$ for all $r$, when $k\le 5$. To avoid redundancy with Table 2.1, it suffices to consider $r\ge 5\ge k$. Again write $r'=r+1$.
 
 \medskip\noindent\textbf{Theorem 5.3.} \textit{The list of all exceptional quantum subgroups of $\cC(A_r,k)$, up to equivalence, when $k\le 5\le r$ is provided by Table 2.2.}\medskip
 
\noindent\textit{Proof.} First note that since $\cC(A_r,1)$ is pointed,  its quantum subgroups are all necessarily simple-current ones, so we need consider only $k>1$.
 
 By Corollary 4.10, the only pairs $(r,k)$ which can appear in Table 2.2 are those where $\cC(A_{k-1},r')$ appears in Table 2.1. Those  seven pairs $(r,k)$ are listed  in Table 2.2.  We claim that for these $\cC(A_r,k)$, all quantum subgroups lie in $\cC(A_r,k)_0$.  Since this is also true for their level-rank duals, thanks to the Proposition 4.11 argument, this would establish by Theorem  5.1 of \cite{OsS} the existence and uniqueness of their exceptional quantum subgroups. 
 
 By Proposition 4.11, all this is manifest  for $\cC(A_5,4)$, $\cC(A_6,5)$, $\cC(A_9,2)$, $\cC(A_{20},3)$, and $\cC(A_{27},2)$. On the other hand,  for $\cC(A_7,4)$ and $\cC(A_8,3)$,  \eqref{Anormsimpl} leaves open the possibility that $\cA$ can have subobjects  $\la$ with $t(\la-\rho)\equiv 4$ (mod $8$) resp.\ $t(\la-\rho)\equiv \pm3$ (mod $9$). However, direct inspection verifies that such $\la$ in these two categories cannot have $\theta(\la)=1$. This verification is easy: there are only 13 resp.\ 8 orbits of these $\la$ by the groups $\lan J_a^2,C\ran$ resp.\ $\lan J_a^3,C\ran$ (those groups preserve the ribbon twist here). \textit{QED to Theorem 5.3}\medskip

Thus it suffices to work out their branching rules. Those at levels 2 and 3 were determined in \cite{su23} using modularity so it is unnecessary to revisit them.

$\cC(A_r,k)_0$ is a braided fusion subcategory of the MFC $\cC(A_r,k)$. When $\cA\in \cC(A_r,k)_0$, we see from \eqref{adjoint} that whenever $M\in\mathrm{Rep}_{\cC(A_r,k)}\cA$ is indecomposable and Res$\,M$ has a nontrivial intersection with $\cC(A_r,k)_0$, then Res$\,M\in \cC(A_r,k)_0$. Hence such an $M$ will also lie in $\mathrm{Rep}_{\cC(A_r,k)_0}\cA$. Furthermore, by definition $M$ is local in $\mathrm{Rep}_{\cC(A_r,k)}\cA$ iff it is local in $\mathrm{Rep}_{\cC(A_r,k)_0}\cA$. The branching rules $B_{M,\la}$ for such $M$ are the same whether we regard them in $\mathrm{Rep}_{\cC(A_r,k)}\cA$ or in $\mathrm{Rep}_{\cC(A_r,k)_0}\cA$. But we know $\tau_0^*$ is a braided tensor equivalence between  ${\cC(A_r,k)_0}$ and  $(\cC(A_{k-1},r')_0)^{\mathrm{rev}}$, so the indecomposable local modules in $\mathrm{Rep}_{\cC(A_r,k)_0}\cA$ are in natural bijection with those of  $\mathrm{Rep}_{(\cC(A_{k-1},r')_0)^{\mathrm{rev}}}\cA$, and $\tau_0^*$ intertwines their branching rules.

Using Table 2.1, this gives us the branching rules in Table 2.2 for $(r,k)=(5,4),(6,5),(7,4)$, when $j=0$. To get the rest, consider first
$A_5$ at level 4. 
As in the proof of Corollary 4.8, $J^e=\mathrm{Ind}\,J_a$ is  a simple-current $J_e$ in the MFC (Rep$_{\cC(A_5,4)}\cA)^{\mathrm{loc}}$. Since $\cJ_\cA=\lan J_a^3\ran$, any $\la\in P_{\!++}^{\,\ka}(A_5)$ appearing in the restriction of a local $\cA$-module $M$ obeys $2|t(\la-\rho)$. Using \eqref{tjmu}, we have $J_a^{-t(\la-\rho)}\la\in\cC(A_5,4)_0$ for any such $\la$. In other words, for any indecomposable local $M\in \mathrm{Rep}_{\cC(A_5,4)}\cA$, there will be a $j$ for which $J_e^jM\in \mathrm{Rep}_{\cC(A_5,4)_0}\cA$. This gives the $j\ne0$ entries in that row of Table 2.2.  The arguments for 
$\cC (A_6,5)$ and $\cC(A_7,4)$ are similar: for $\cC (A_6,5)$ use $\la\mapsto J_a^{4t(\la-\rho)}\la$ and for $\cC (A_7,4)$ use $\la\mapsto J_a^{t(\la-\rho)}\la$.

 \section{Module category classifications}\label{sect:modcat}
 
 For completeness, this section lists all module categories for $\fg=A_1,...,A_4$. These for $A_1$ are classified e.g.\  in \cite{KO}, based on \cite{CIZ}, and fall into the ADE pattern. Those for $A_2$ are classified in \cite{EP}, based on \cite{su3}, and have connections to Jacobians of Fermat curves \cite{BCIR}. Recently, those for $A_3,A_4,A_5,A_6,G_2$ are classified in \cite{EM3,EM2,EM}, assuming  results from this paper and \cite{Gan-ii}. The module categories for $\fg=A_3$ are studied in full detail in \cite{CM}. Ocneanu announced \cite{Oc} the $A_3$  classification, though the arguments never appeared.
 
Recall  that a module category for Mod$(\cV)$ can be thought of as a triple $(\cV^e_L,\cV^e_R,\cF)$ where $\cV^e_{L,R}$ are conformal extensions of $\cV$ and $\cF:\mathrm{Mod}(\cV^e_L)\to\mathrm{Mod}(\cV^e_R)$ is a braided tensor equivalence. Associated to a module category is its modular invariant $\cZ=B_R^t\Pi B_R$, 
where $B_{L,R}$ are the branching (restriction) matrices associated to $\cV^e_{L,R}\supset\cV$ and $\Pi$ is the bijection Irr$(\mathrm{Mod}(\cV^e_L))\to\mathrm{Irr(Mod}(\cV^e_R))$. We identify the modular invariant $\cZ$ with the formal combination $\sum_{\la,\mu}\cZ_{\la,\mu}\chi_\la\overline{\chi_{\mu}}$.

Consider any $\cC(A_r,k)$. Write ${r}'=r+1$, and put ${k}'=\kappa$ if $k{r}'$ is odd, or  ${k}'=k$ otherwise. Choose any divisor $d$ of ${r}'$; if $d$ is even we require as well that ${r}'{k}'/d$ be even. Then we get a module category with modular invariant \begin{equation}\cZ[\cJ_d]_{\la,\mu}=\sum_{j=1}^d\delta^d(t(\la)+j{r}'\,{k}'/(2d))\,\delta_{\mu,J^{j{r}'/d}\la}\end{equation}
where $\delta^y(x)=1$ or 0 depending, respectively, on whether or not $x/y\in\bbZ$. It defines an auto-equivalence of the simple-current extension of $\cV(A_r,k)$ by $\langle J^{d'}\rangle$ where $d'=\mathrm{gcd}(d,{r}'/d)$. Note that $\cZ[\cJ_1]=I$.

When $r>1$ and $k\le 2$,  the products $C\cZ[\cJ_d]$ should be dropped.

 \medskip\noindent\textbf{Theorem A1.} \textit{The complete list of  module
categories for $\cC(A_1,k)$  is:}

\noindent{$\bullet$} \textit{those of simple-current type $\cZ[\cJ_d]$;} 

 \noindent{$\bullet$} \textit{at $k=16$, the exceptional auto-equivalence of the simple-current extension $\cZ[\cJ_2]$, whose modular invariant has the exceptional terms}
 $$\cdots 
+ (\chi_2+\chi_{14})\overline{\chi_8}+\chi_8(\overline{\chi_2}+\overline{\chi_{14}})+|\chi_8|^2\,;$$
 
  \noindent{$\bullet$} \textit{the conformal embeddings $ \cV(A_1,10)\subset\cV(C_2,1)$ and $\cV(A_1,28)\subset\cV(G_2,1)$:}
  \begin{align}&|\chi_{\bo}+\chi_{6}|^2+|\chi_{3}+\chi_7|^2+|\chi_4+\chi_{10}|^2 \,,\nonumber\\
&|\chi_{\bo}+\chi_{10}+\chi_{18}+\chi_{28}|^2+|\chi_{6}+\chi_{12}+\chi_{16}+\chi_{22}|^2\,.\nonumber\end{align}

 \medskip\noindent\textbf{Theorem A2.} \textit{The complete list of  module
categories for $\cC(A_2,k)$  is:}
 
 \noindent{$\bullet$} \textit{those of simple-current type $\cZ[\cJ_d]$ and (for $k>2$) their products $C\cZ[\cJ_d]$;}
 
  \noindent{$\bullet$} \textit{at $k=9$, the exceptional auto-equivalence of the simple-current extension $\cZ[\cJ_3]$, whose modular invariant has the exceptional terms}
 $$\cdots
 +(\chi_{11}+\chi_{71}+\chi_{17})\overline{\chi_{33}}+\chi_{33}(\overline{\chi_{11}}+\overline{\chi_{71}}+\overline{\chi_{17}})+|\chi_{33}|^2$$
\textit{ as well as its product with $C$;}

  \noindent{$\bullet$} \textit{the conformal embeddings $ \cV(A_2,5)\subset\cV(A_5,1)$, $\cV(A_2,9)\subset\cV(E_6,1)$, and $\cV(A_2,21)\subset\cV(E_7,1)$, with modular invariants}
  \begin{align}&|\chi_{\bo}+\chi_{22}|^2+|\chi_{02}+\chi_{32}|^2+|\chi_{20}+\chi_{23}|^2+|\chi_{21}+\chi_{05}|^2+|\chi_{30}+\chi_{03}|^2+|\chi_{12}+\chi_{50}|^2 \,,\nonumber\\
&|\chi_{\bo}+\chi_{09}+\chi_{90}+\chi_{44}+\chi_{41}+\chi_{14}|^2+2|\chi_{22}+\chi_{25}+\chi_{52}|^2\,,\nonumber\\
&|\lan\chi_{\bo}\ran_3+\lan\chi_{44}\ran_3+\lan\chi_{66}\ran_3+\lan\chi_{10,10}\ran_3|^2+|\lan\chi_{15,6}\ran_{3c}+\lan\chi_{10,7}\ran_{3c}|^2\,,\nonumber\end{align}
\textit{as well as products of the first two with $C$.}

 \medskip\noindent\textbf{Theorem A3.} \textit{The complete list of  module
categories for $\cC(A_3,k)$  is:}

\noindent{$\bullet$} \textit{those of simple-current type $\cZ[\cJ_d]$, and (for $k>2$) their products $C\cZ[\cJ_d]$;}

 \noindent{$\bullet$} \textit{at $k=8$, the exceptional auto-equivalence of the simple-current extension $\cZ[\cJ_4]$, whose modular invariant has the exceptional terms }
 $$\cdots +2|\chi_{222}|^2+\lan\chi_{012}\ran_{4c}\overline{\chi_{222}}+ \chi_{222}\lan\overline{\chi_{012}}\ran_{4c}\,,$$
\textit{together with $C$ times it;}

 \noindent{$\bullet$} \textit{the conformal embeddings $\cV(A_3,4)\subset\cV(B_7,1)$, $\cV(A_3,6)\subset\cV(A_9,1)$, and $\cV(A_3,8)\subset\cV(D_{10},1)$, with modular invariants}
   \begin{align}&|\lan\chi_\bo\ran_2+\lan\chi_{012}\ran_2|^2+|\lan\chi_{400}\ran_2+\lan\chi_{101}\ran_2|^2+|2\chi_{111}|^2\,,\cr
   &|\lan\chi_{\bo}\ran_2+\lan\chi_{202}\ran_2|^2+|\lan\chi_{600}\ran_{2}+\lan\chi_{220}\ran_2|^2+|\lan\chi_{200}\ran_2+\chi_{212}|^2+|\lan\chi_{002}\ran_2+\chi_{212}|^2+|\lan\chi_{420}\ran_2+\chi_{121}|^2\cr&\ \ +|\lan\chi_{024}\ran_2+\chi_{121}|^2+|\lan\chi_{210}\ran_2+\chi_{303}|^2+|\lan\chi_{012}\ran_2+\chi_{303}|^2+|\lan\chi_{321}\ran_2+\chi_{030}|^2+|\lan\chi_{123}\ran_2+\chi_{030}|^2\,,\cr
 &  |\lan\chi_\bo\ran_4+\lan\chi_{121}\ran_4|^2+|\lan\chi_{020}\ran_4+\lan\chi_{303}\ran_4|^2+2|\lan\chi_{113}\ran_4|^2\,,\nonumber\end{align}
\textit{together with $C$ times the second and an auto-equivalence of the third leaving the modular invariant unchanged.}

\medskip\noindent\textbf{Theorem A4.} \textit{The complete list of  module
categories for $\cC(A_4,k)$  is:}

\noindent{$\bullet$} \textit{those of simple-current type $\cZ[\cJ_d]$, and (for $k>2$) their products $C\cZ[\cJ_d]$;}

 \noindent{$\bullet$} \textit{at $k=5$, the exceptional auto-equivalence of the simple-current extension $\cZ[\cJ_5]$, whose modular invariant has the exceptional terms }
 $$\cdots +4|\chi_{1111}|^2+\lan\chi_{1001}\ran_5\overline{\chi_{1111}}+ \chi_{1111}\lan\overline{\chi_{1001}}\ran_5\,,$$
\textit{together with $C$ times it;}

  \noindent{$\bullet$} \textit{the conformal embeddings $\cV(A_4,3)\subset\cV(A_9,1)$, $\cV(A_4,5)\subset\cV(D_{12},1)$, and $\cV(A_4,7)\subset\cV(A_{14},1)$, with modular invariants}
   \begin{align}&|\chi_\bo+\chi_{0110}|^2+|\chi_{3000}+\chi_{1011}|^2+|\chi_{0300}+\chi_{0101}|^2+|\chi_{0030}+\chi_{1010}|^2+|\chi_{0003}+\chi_{1010}|^2\cr&\ \ +|\chi_{0010}+\chi_{0201}|^2+|\chi_{2001}+\chi_{0020}|^2+|\chi_{0200}+\chi_{1002}|^2+|\chi_{1020}+\chi_{0100}|^2+|\chi_{0102}+\chi_{2010}|^2\,,\cr
&|\lan\chi_\bo\ran_5+\lan\chi_{0220}\ran_5|^2+|\lan\chi_{1001}\ran_5+\chi_{1111}|^2+2|2\chi_{1111}|^2\,,\cr
&|\chi_{\bo}+\chi_{0330}+\chi_{2002}+\chi_{2112}+\lan\chi_{0403}\ran_c|^2+|\chi_{7000}+\chi_{1033}+\chi_{3200}+\chi_{1211}+\lan\chi_{0040}\ran_c|^2\cr&\ \ +|\chi_{0700}+\chi_{0103}+\chi_{2320}+\chi_{2121}+\lan\chi_{3004}\ran_c|^2+|\chi_{0070}+\chi_{3010}+\chi_{0232}+\chi_{1212}+\lan\chi_{0300}\ran_c|^2\cr&\ \ +|\chi_{0007}+\chi_{3301}+\chi_{0023}+\chi_{1121}+\lan\chi_{4030}\ran_c|^2
+|\chi_{2000}+\chi_{0312}+\chi_{1240}+\chi_{2102}+\chi_{3022}|^2\cr&\ \ +
|\chi_{5200}+\chi_{1031}+\chi_{0124}+\chi_{2210}+\chi_{0302}|^2+|\chi_{0520}+\chi_{2103}+\chi_{0012}+\chi_{2221}+\chi_{2030}|^2\cr&\ \ +|\chi_{0052}+\chi_{1210}+\chi_{4001}+\chi_{0222}+\chi_{2203}|^2+|\chi_{0005}+\chi_{3121}+\chi_{2400}+\chi_{1022}+\chi_{0220}|^2
\cr&\ \ +|\chi_{0002}+\chi_{2130}+\chi_{0421}+\chi_{2012}+\chi_{2203}|^2+|\chi_{0025}+\chi_{1301}+\chi_{4210}+\chi_{0122}+\chi_{2030}|^2\cr&\ \ +|\chi_{0250}+\chi_{3012}+\chi_{2100}+\chi_{1222}+\chi_{0302}|^2+|\chi_{2500}+\chi_{0121}+\chi_{1004}+\chi_{2220}+\chi_{3022}|^2\cr&\ \ +|\chi_{5000}+\chi_{1213}+\chi_{0042}+\chi_{2201}+\chi_{0220}|^2\,,\nonumber\end{align}
\textit{together with $C$ times the first, five other auto-equivalences of the second corresponding to any permutation of $\Lambda^e_1,\Lambda^e_{11},\Lambda^e_{12}$, and three other auto-equivalences of the third corresponding to $\cZ^e[\cJ^e_{3}],\cZ^e[\cJ^e_{5}],\cZ^e[\cJ^e_{15}]=C^e$.}

\bigskip\noindent\textbf{Acknowledgements.} This manuscript was largely prepared at Max-Planck-Institute for Mathematics in Bonn, who we thank for generous support. We learned of Schopieray's crucial work \cite{Sch} in the conference ``Subfactors in Maui'' in July 2017, without which this paper may not exist. We thank Cain Edie-Michell, David Evans, Liang Kong, Andrew Schopieray and Mark Walton for comments on a preliminary draft. This research was supported in part by NSERC.

  \newcommand\bibx[4]   {\bibitem{#1} {#2:} {\sl #3} {\rm #4}}

\end{document}